\documentclass[12pt]{amsart}
\usepackage{indentfirst,amsmath,amssymb,amsthm,mathrsfs}
\usepackage[colorlinks,linkcolor=blue,citecolor=red]{hyperref}

\numberwithin{equation}{section}
\theoremstyle{plain}
\newtheorem{theorem}{Theorem}[section]
\newtheorem{lemma}{Lemma}[section]
\newtheorem{corollary}{Corollary}[section]
\newtheorem{proposition}{Proposition}[section]
\theoremstyle{definition}
\newtheorem{definition}{Definition}[section]
\newtheorem{example}{Example}[section]
\newtheorem{problem}{Problem}[section]
\theoremstyle{remark}
\newtheorem{remark}{Remark}[section]
\newtheorem*{remark*}{Remark}

\providecommand{\Iref}[2]{#2}

\begin{document}
\date{}
\title[Shrinking Sasaki--Ricci Solitons II]{On the Structure of Gradient
Shrinking Sasaki--Ricci Solitons of Dimension Five II: Global Quotient
Classification and Irregular Reeb Foliation}
\author{Shu-Cheng Chang}
\address{Shanghai Institute for Mathematics and Interdisciplinary Sciences,
Shanghai 200433, China}
\email{scchang@math.ntu.edu.tw}
\author{Chien Lin}
\address{Mathematical Science Research Center, Chongqing University of
Technology, Chongqing 400054, China}
\email{chienlin@cqut.edu.cn}
\author{Hongbing Qiu}
\address{School of Mathematics and Statistics, Wuhan University, Wuhan
430072, China}
\email{hbqiu@whu.edu.cn}

\begin{abstract}
This is the second of two papers on complete noncompact gradient shrinking
Sasaki--Ricci solitons of dimension five. Building on the transverse
curvature estimates and the local splitting at transverse infinity
established in the first paper, we study the global structure of these
solitons according to the type of their Reeb foliation. For regular Reeb
foliations the leaf space is one of the four noncompact K\"{a}hler--Ricci
shrinker surfaces (Gaussian, Feldman--Ilmanen--Knopf, round cylinder and
Bamler--Cifarelli--Conlon--Deruelle), and we determine all regular lifts and
their topology. The quasi-regular case with developable leaf space reduces
to the regular one by a finite covering. If the Reeb foliation is not
quasi-regular, it is globally developable if and only if the soliton is the
Heisenberg model or an irrational Gaussian suspension, which happens if and
only if $R^{T}\equiv 0$. In the remaining nondevelopable case we show that,
up to finite quotients and the transversely flat models, the solitons of a
one-parameter family $\mathcal{S}_{a}$ on $\mathbb{C}\times \mathbb{S}^{3}$
are the only globally split shrinkers, and that global splitting is
equivalent to nonnegative transverse holomorphic bisectional curvature. This
serves as a Sasakian analogue of Naber's theorem on four-dimensional
shrinkers with nonnegative curvature operator.
\end{abstract}

\maketitle
\tableofcontents

\section{Introduction}

\label{sec:intro}

This is the second of two papers on complete noncompact gradient shrinking
Sasaki--Ricci solitons of dimension five; we refer to the introduction of
\cite{PartI} for background on the Sasaki--Ricci flow \cite{SWZ,Co1,He} and
its solitons \cite{FOW,CLL,CHW}. Throughout, a gradient shrinking
Sasaki--Ricci soliton $(M,\xi ,\eta ,\Phi ,g,f)$ is a Sasakian manifold
together with a basic function $f$ satisfying the full real transverse
equation
\begin{equation}
\mathrm{Ric}^{T}+\operatorname{Hess}^{T}f=\lambda g^{T},\qquad \lambda >0.
\label{1}
\end{equation}
Up to a transverse $D$-homothety, one may take $\lambda =2n+2$; the
additive constant in the potential is fixed by the Hamilton identity
\eqref{eq:pre-normal}. In dimension five the default scale in
Sections~\ref{sec:local} and~\ref{sec:structure} is $\lambda=6$, unless
another scale is explicitly stated. Formulas retaining $\lambda$ apply
at the indicated positive scale, and the temporary normalization
$\lambda=\tfrac12$ used for Ricci-flow limits is stated when invoked.

An isomorphism of Sasakian structures is a diffeomorphism preserving
$\xi,\eta,\Phi$ and $g$. An isomorphism of normalized solitons also
preserves the normalized potential $f$; an automorphism is such an
isomorphism from a soliton to itself. These identifications are at the
fixed value of $\lambda$, unless a $D$-homothety is explicitly specified.
A classification of quotient metrics up to biholomorphic pullback means
that a biholomorphism pulls back the stated model metric to $h_Z$.
In statements about the topology of total spaces, $\cong$ denotes a
diffeomorphism of the underlying smooth manifolds, not a product
isometry of their Sasakian metrics.

In \cite{PartI} we proved, without any regularity assumption on the Reeb
foliation, that an upper bound for the transverse scalar curvature $R^{T}$
implies
\begin{equation*}
|\mathrm{Rm}^{T}|\leq CR^{T}
\end{equation*}
(Theorem~\ref{thm:intro-curvature}), which is the Sasakian analogue of the
estimate of Munteanu--Wang \cite{MW1}, and that at transverse infinity,
i.e.\ as $f\rightarrow \infty $, every developed local limit splits off a
flat complex line (Theorem~\ref{thm:splitting}) as in the splitting theorem
of Naber \cite{Nab}.

Here and below, transverse infinity means $f\to\infty$. A gradient
trajectory is called escaping if $f(\gamma(s))\to\infty$ in its forward
direction. Uniform scalar decay means that for every $\varepsilon>0$
there is $A$ such that $|R^T|<\varepsilon$ on $\{f\geq A\}$.
Accordingly, either $R^T$ decays uniformly at transverse infinity or
there are $\delta>0$ and points $p_i$ with $f(p_i)\to\infty$ and
$R^T(p_i)\geq\delta$. These local results do not by themselves give a global
classification. The purpose of the present paper is therefore to derive the
global structure according to the type of the Reeb foliation
(Section~\ref{sec:types}): regular, quasi-regular, a free proper
$\mathbb{R}$-action, or
irregular.

\textbf{Regular Reeb foliations:} If the Reeb field generates a free circle
action of period $L>0$, then
\begin{equation*}
\pi :M\rightarrow Z=M/\mathbb{S}^{1}
\end{equation*}
is a principal circle bundle over a complete smooth K\"{a}hler surface
$(Z,h_{Z},J_{Z})$ with $\frac{1}{2}d\eta =\pi ^{\ast }\omega _{Z}$, $f=\pi
^{\ast }f_{Z}$, and the soliton descends to the gradient K\"{a}hler--Ricci
shrinker equation
\begin{equation*}
\mathrm{Ric}_{Z}+\operatorname{Hess}_{Z}f_{Z}=\lambda h_{Z}.
\end{equation*}
By Li--Wang \cite{LW} the curvature of $Z$ is then automatically bounded
(Theorem~\ref{thm:quotient}), and the classification of the noncompact
K\"{a}hler--Ricci shrinker surfaces by Conlon--Deruelle--Sun \cite{CDS},
Cifarelli--Conlon--Deruelle \cite{CCD} and
Bamler--Cifarelli--Conlon--Deruelle \cite{BCCD} applies (see also \cite{LZ}).
In the regular case our results are therefore consequences of this
classification together with an analysis of the circle bundles. We separate
the two behaviors at infinity. In the conical case:

\begin{theorem}
\label{thm:intro-decay} Let $(M^{5},\xi ,\eta ,\Phi ,g,f)$ be a complete
noncompact regular gradient shrinking Sasaki--Ricci soliton. Assume that
\begin{equation*}
R^{T}(p)\longrightarrow 0\qquad \text{uniformly as}\qquad
f(p)\longrightarrow \infty .
\end{equation*}
Then, up to biholomorphism, the quotient $(Z,h_{Z},f_{Z})$ is either the
Gaussian shrinker on $\mathbb{C}^{2}$ or the FIK shrinker \cite{FIK} on
\begin{equation*}
\operatorname{Bl}_{0}\mathbb{C}^{2}=\operatorname{Tot}\bigl(\mathcal{O}(-1)\rightarrow
\mathbb{CP}^{1}\bigr).
\end{equation*}
Moreover, by \cite{CDS}, $Z$ has a unique tangent K\"{a}hler cone at
infinity, biholomorphic to $\mathbb{C}^{2}$.
\end{theorem}

In the cylindrical case:

\begin{theorem}
\label{thm:intro-cyl} Let $(M^{5},\xi ,\eta ,\Phi ,g,f)$ be a complete
noncompact regular gradient shrinking Sasaki--Ricci soliton. Assume that
$\nabla ^{T}f$ has an escaping integral curve $\gamma $ along which
$R^{T}\nrightarrow 0$. Then, up to biholomorphism, either

\begin{enumerate}
\item
\begin{equation*}
Z\cong \mathbb{C}\times \mathbb{CP}^{1},
\end{equation*}
carrying the product of the Gaussian and the round shrinker, or

\item
\begin{equation*}
Z\cong \operatorname{Bl}_{p}(\mathbb{C}\times \mathbb{CP}^{1}),
\end{equation*}
carrying the BCCD shrinker \cite{BCCD}.
\end{enumerate}

In both cases, along the projected trajectory $\pi (\gamma )$ the pointed
quotients converge in the smooth pointed Cheeger--Gromov sense to the round
K\"{a}hler cylinder $(\mathbb{C}\times
\mathbb{S}^{2},h_{\mathrm{E}}+h_{\lambda })$, where $h_{\lambda }$ is the
round metric of Gaussian
curvature $\lambda $, and $R^{T}\rightarrow 2\lambda $ along $\gamma $. The
limit is independent of the sequence. The full Sasakian structures converge
along $\gamma $, after gauge transformations, to the cylindrical connection
lift of Remark~\ref{rem:cylmodel} with $d=8\pi /(\lambda L)$, where $L$ is
the Reeb period (see the paragraph following Remark~\ref{rem:cylmodel}).
\end{theorem}

\begin{corollary}
\label{cor:four-models} Every complete noncompact regular five-dimensional
gradient shrinking Sasaki--Ricci soliton has, up to biholomorphic pullback,
exactly one of the smooth quotient metrics
\begin{equation}
\mathbb{C}^{2}\ \text{(Gaussian)},\qquad \operatorname{Bl}_{0}\mathbb{C}^{2}\
\text{(FIK)},\qquad \mathbb{C}\times \mathbb{CP}^{1}\ \text{(product)},\qquad
\operatorname{Bl}_{p}(\mathbb{C}\times \mathbb{CP}^{1})\ \text{(BCCD)}.
\label{eq:fourmodels}
\end{equation}
\end{corollary}

These are the four transverse models for the regular Reeb foliation. The
total space is determined by the quotient with the Reeb period $L$ and the
integral Euler class of the circle bundle. A regular lift of $(Z,h_{Z},f_{Z})
$ with Reeb period $L$ exists if and only if
\begin{equation}
\frac{2}{L}[\omega _{Z}]=\frac{4\pi }{\lambda L}\,c_{1}(Z)_{\mathbb{R}}
\label{eq:intro-euler}
\end{equation}
is the real image of an integral class $e\in H^{2}(Z;\mathbb{Z})$
(Proposition~\ref{prop:lifting}). This leads to the following topological
description (Theorems~\ref{thm:regular} and~\ref{thm:decay-lifts}).

We fix the circle-bundle notation and signs before stating the result.
Give $\mathbb{CP}^{1}$ its complex orientation and put
$u=c_{1}(\mathcal O(1))$, so
$\langle u,[\mathbb{CP}^{1}]\rangle=1$. For a principal circle bundle
$p:P\to\mathbb{CP}^{1}$, our convention is
$e(p)_{\mathbb R}=[F_{\theta}/(2\pi)]$, where
$d\theta=p^{*}F_{\theta}$ and $\theta(U)=1$ for the positive circle
generator $U$ of period $2\pi$. For each integer $k\geq1$, let
$\mu_k=\{\zeta\in\mathbb C:\zeta^k=1\}$ act diagonally on
$\mathbb S^3\subset\mathbb C^2$, and define
\[
L(k,1):=\mathbb S^3/\mu_k,\qquad
p_k:L(k,1)\longrightarrow\mathbb{CP}^{1},\qquad
p_k([z_0,z_1]_{\mu_k})=[z_0:z_1].
\]
The effective circle action is
$e^{it}\cdot[z]_{\mu_k}=[e^{it/k}z]_{\mu_k}$; the formula is
independent of the choice of $t$ modulo $2\pi$. In the convention just
specified,
\[
e(p_k)=k u,\qquad
\langle e(p_k),[\mathbb{CP}^{1}]\rangle=k.
\]
Thus the integer is the Euler number of the bundle, whereas $p_k$ is
its circle-bundle projection, not a finite covering. In particular,
$p_1$ is the Hopf fibration and $L(1,1)=\mathbb S^3$. The notation
$L(k,1)$ denotes a lens space, while $L$ alone denotes the Reeb period.
For an oriented compact curve $C\subset Z$, we abbreviate
$e(C):=\langle e,[C]\rangle$.

\begin{theorem}
\label{thm:intro-topology}Let $(M^{5},\xi ,\eta ,\Phi ,g,f)$ be a complete
noncompact regular gradient shrinking Sasaki--Ricci soliton with Reeb period
$L$ and quotient $Z$. Let $E$ denote the exceptional curve and $F$ a generic
sphere fiber whenever they exist.

\begin{enumerate}
\item The conical case.

\begin{enumerate}
\item If $Z$ is the Gaussian shrinker, every $L>0$ occurs,
\begin{equation*}
M\cong \mathbb{C}^{2}\times \mathbb{S}^{1}
\end{equation*}
with $e=0$, $\pi _{1}(M)=\mathbb{Z}$ and $H^{2}(M;\mathbb{Z})=0$.

\item If $Z$ is the FIK shrinker, then $q:=4\pi /(\lambda L)\in
\mathbb{Z}_{>0}$,
\begin{equation*}
M\cong \operatorname{Tot}\bigl(p_{q}^{\ast }\mathcal{O}(-1)\rightarrow L(q,1)\bigr)
\end{equation*}
with $e(E)=q$ and $\pi _{1}(M)=H^{2}(M;\mathbb{Z})=\mathbb{Z}/q\mathbb{Z}$,
where $p_q$ is the principal circle bundle defined above, with
$e(p_q)=q u$. Here $\mathcal O(-1)\to\mathbb{CP}^{1}$ is the tautological
complex line bundle, with $c_1(\mathcal O(-1))=-u$, and
$p_q^{*}\mathcal O(-1)$ is its pullback to $L(q,1)$. For $q=1$,
\begin{equation*}
M\cong \mathbb{C}\times \mathbb{S}^{3}.
\end{equation*}
\end{enumerate}

\item The cylindrical case.

\begin{enumerate}
\item If $Z=\mathbb{C}\times \mathbb{CP}^{1}$, then $d:=8\pi /(\lambda L)\in
\mathbb{Z}_{>0}$,
\begin{equation*}
M\cong \mathbb{C}\times L(d,1)
\end{equation*}
with $e(\mathbb{CP}^{1})=d$ and $\pi
_{1}(M)=H^{2}(M;\mathbb{Z})=\mathbb{Z}/d\mathbb{Z}$.

\item If $Z=\operatorname{Bl}_{p}(\mathbb{C}\times \mathbb{CP}^{1})$, then $q:=4\pi
/(\lambda L)\in \mathbb{Z}_{>0}$,
\begin{equation*}
\pi ^{-1}(E)\cong L(q,1),
\end{equation*}
with $e(F)=2q$, $e(E)=q$, $\pi _{1}(M)=\mathbb{Z}/q\mathbb{Z}$ and
$H^{2}(M;\mathbb{Z})=\mathbb{Z}\oplus \mathbb{Z}/q\mathbb{Z}$.
\end{enumerate}
\end{enumerate}

In each case every positive integer $q$, resp.\ $d$, is realized.
\end{theorem}

\textbf{Quasi-regular Reeb foliations: }When the Reeb foliation is only
quasi-regular, the leaf space $\mathcal{Z}=M/\mathbb{S}^{1}$ is a K\"{a}hler
orbifold and $M\rightarrow \mathcal{Z}$ is a circle orbibundle \cite{BG,Kol1}.
Quasi-regularity alone does not provide a smooth manifold cover of
$\mathcal{Z}$, and the regular results do not descend in general
(Remark~\ref{rem:weighted-example}). In the developable case we obtain a
reduction to
the regular one.

\begin{definition}
\label{def:developable}An effective orbifold is \emph{developable} if its
universal orbifold cover is a manifold. Equivalently, it can be written
$[\widetilde{Z}/\Gamma ]$ with $\widetilde{Z}$ a simply connected manifold and
$\Gamma $ a discrete group acting properly by diffeomorphisms with finite
stabilizers. The action need not be free; for instance, a quotient by a
finite group with fixed points is allowed.
\end{definition}

\begin{remark}
Let $\mathcal{Z}$ be a complete noncompact gradient shrinking K\"{a}hler
orbifold whose potential is globally defined. If $\mathcal{Z}$ is
developable, its universal orbifold cover is a finite cover
(Lemma~\ref{lem:finite}).
\end{remark}

\begin{theorem}
\label{thm:good-intro} Let $(M^{5},\xi ,\eta ,\Phi ,g,f)$ be a complete
noncompact quasi-regular gradient shrinking Sasaki--Ricci soliton with
generic period $L$ and developable quotient $\mathcal{Z}$. There is a finite
ordinary covering $\widetilde{M}\rightarrow M$ preserving the pulled-back
Sasakian structure and potential, such that the Reeb flow on $\widetilde{M}$
is regular with period $L$. Its smooth K\"{a}hler quotient $\widetilde{Z}$ is
one of the four smooth models of Corollary~\ref{cor:four-models}. More
precisely,

\begin{enumerate}
\item if $R^{T}\rightarrow 0$ uniformly at transverse infinity,
$\widetilde{Z}$ is Gaussian or FIK;

\item if an escaping gradient trajectory has $R^{T}\nrightarrow 0$,
$\widetilde{Z}$ is the round K\"{a}hler cylinder or the BCCD
surface;

\item the original manifold is a finite free Sasakian quotient
\begin{equation}
M=\widetilde{M}/\Gamma ,  \label{eq:finitequotient-intro}
\end{equation}
where $\Gamma =\pi _{1}^{\mathrm{orb}}(\mathcal{Z})$ acts on $\widetilde{M}$
preserving $\eta ,\xi ,\Phi ,g,\widetilde{f}$ and commuting with the circle
action.
\end{enumerate}
\end{theorem}

\begin{corollary}
Let $(M^{5},\xi ,\eta ,\Phi ,g,f)$ be a complete noncompact gradient
shrinking Sasaki--Ricci soliton. If the Reeb action is closed, free and
proper with group $\mathbb{R}$, its smooth quotient is the Gaussian shrinker
on $\mathbb{C}^{2}$, and the total Sasakian metric is, in compatible
coordinates, the Heisenberg model (Corollary~\ref{cor:Rcase}).
\end{corollary}

\begin{theorem}[Corollary~\protect\ref{cor:nondev-qr}]
Let $M$ be a complete noncompact quasi-regular five-dimensional gradient
shrinking Sasaki--Ricci soliton whose orbifold leaf space is not
developable. Then

\begin{enumerate}
\item $R^{T}>0$, $\pi _{1}(M)$ is finite, and $M=\widehat{M}/\Gamma $, where
the universal cover $\widehat{M}$ is quasi-regular with a nondevelopable
leaf space.

\item $\Gamma \cong \pi _{1}(M)$ is a finite group of automorphisms of the
lifted soliton which acts freely and commutes with the circle action. In
particular $\widehat{M}$ is not regular.
\end{enumerate}
\end{theorem}

\begin{example}[A nondevelopable quasi-regular soliton]
\label{ex:weighted-example}Let $m_{1},m_{2}$
be coprime positive integers with $(m_{1},m_{2})\neq (1,1)$. The soliton of
Remark~\ref{rem:weighted-example} is the circle lift, with generic Reeb
period $L=4\pi (m_{1}+m_{2})/\lambda $, of the orbifold shrinker
\begin{equation*}
\mathcal{Z}=\mathbb{C}\times \mathbb{CP}_{[m_{1},m_{2}]}^{1}
\end{equation*}
with $h_{\mathcal{Z}}=|dz|^{2}+h_{m}$ and $f_{\mathcal{Z}}=\frac{\lambda
}{2}|z|^{2}+f_{m}$, where $h_{m}$ is the shrinking soliton of Wu \cite{Wu} on
the weighted projective line. It is a complete gradient shrinking
Sasaki--Ricci soliton on
\begin{equation*}
M=\mathbb{C}\times \mathbb{S}^{3}
\end{equation*}
with bounded $R^{T}=\pi ^{\ast }R_{h_{m}}$. It is quasi-regular but not
regular, the orbits over the two orbifold points have minimal periods
$L/m_{1}$ and $L/m_{2}$, its leaf space is not developable, and it is not one
of the four models of Corollary~\ref{cor:four-models}.
\end{example}

\textbf{Reeb foliations which are not quasi-regular: }For an irregular
foliation, there is no Hausdorff leaf space, and the relevant notion is the
following.

\begin{definition}
\label{def:globdev}The Reeb foliation of a connected complete
five-dimensional Sasakian manifold $M$ is \emph{globally developable} if,
on the universal cover
$\widehat{M}$ with the lifted structure, there is a surjective submersion
\begin{equation*}
\mathcal{D}:\widehat{M}\rightarrow \widehat{Z}
\end{equation*}
onto a connected smooth Hausdorff four-manifold whose fibers are
exactly the leaves of the lifted Reeb foliation. The target is equipped
with the descended transverse metric and complex structure. This is a
global quotient of the universal cover; the existence of a compatible
extension of a local transverse chart alone is not this condition.
\end{definition}

\begin{example}
The \emph{Gaussian suspension} $M_{A,\ell }$ is the quotient of the
Heisenberg model by $(z,t)\mapsto (Az,t+\ell )$, with $A\in U(2)$ and $\ell
>0$ (Example~\ref{ex:suspension}). It is irregular exactly when $A$ has
infinite order; we then call it an \emph{irrational} Gaussian suspension.
The fundamental group is finite unless the soliton is a Gaussian suspension
(Theorem~\ref{thm:pi1}).
\end{example}

\begin{theorem}
\label{cor:globdev-class-intro} A complete noncompact five-dimensional
gradient shrinking Sasaki--Ricci soliton with globally developable Reeb
foliation is isomorphic to exactly one of the following:

\begin{enumerate}
\item a finite free quotient $\widetilde{M}/\Gamma $, as in
\eqref{eq:finitequotient}, of a regular lift $\widetilde{M}$ of one of the
four models \eqref{eq:fourmodels};

\item the Heisenberg model $\mathbb{H}_{2}$; or

\item an irrational Gaussian suspension $M_{A,\ell }$, determined by $\ell $
and the conjugacy class of $A$.
\end{enumerate}

\noindent The FIK, round-cylinder and BCCD models occur only in $(1)$, with
generic Reeb period at most $8\pi /\lambda $.
\end{theorem}

\textbf{Quasi-regular soliton approximation: }For irregular Reeb foliations,
we follow the strategy used for compact transverse Fano Sasakian
five-manifolds in \cite[Section~4]{CHLWHT}. By the second structure theorem
of Sasakian geometry \cite{Ru,BG}, whose complete-manifold version with
explicit moment-map hypotheses is Proposition~\ref{prop:second-structure},
rational generators approximate the irregular generator in its compact torus.

However, the approximation has to be made by solitons: a deformation of the
Reeb field does not preserve the soliton equation
(Section~\ref{sec:qr-approx}), and the compactness theorem is formulated for
the
five-dimensional total spaces (Theorem~\ref{thm:compactness-total}). The
leaf spaces of quasi-regular approximations may degenerate with the
possibility that their generic Reeb periods tend to infinity
(Lemma~\ref{lem:period-blowup}). See Problem \ref{prob:open}.

The main structure theorem for Reeb foliations which are not quasi-regular
is the following (Theorem~\ref{thm:irregular-main}).

\begin{theorem}
\label{thm:irregular-main-intro} Let $(M^{5},\xi ,\eta ,\Phi ,g,f)$ be a
complete noncompact gradient shrinking Sasaki--Ricci soliton whose Reeb
foliation is not quasi-regular. The following are equivalent.

\begin{enumerate}
\item The Reeb foliation is globally developable.

\item $M$ admits a developable quasi-regular soliton approximation
(Definition~\ref{def:qr-approx}).

\item $R^{T}\equiv 0$.

\item $M$ is isomorphic to $\mathbb{H}_{2}$ or to an irrational Gaussian
suspension $M_{A,\ell }$.

\item $\pi _{1}(M)$ is infinite, or the Reeb flow of $M$ has no closed
orbit.
\end{enumerate}
\end{theorem}

Condition (2) presupposes that a quasi-regular soliton approximation
exists; whether every soliton whose Reeb closure is a torus of dimension
at least two admits one is open (Problem~\ref{prob:open}(2)).

\begin{remark}
\begin{enumerate}
\item If these conditions hold, then every developable quasi-regular soliton
approximation $(M_{i})$ of $M$ has, for all large $i$, Gaussian regular
covers
\begin{equation*}
\widetilde{M}_{i}\cong \mathbb{C}^{2}\times \mathbb{S}^{1}
\end{equation*}
with periods $L_{i}\rightarrow \infty $. In particular no FIK,
round-cylinder or BCCD branch of Theorem~\ref{thm:good} occurs in it.

\item If they fail, then the closure of the Reeb flow is a torus $T^{k}$,
$k\in \{2,3\}$, $R^{T}>0$, $\pi _{1}(M)$ is finite, the universal cover of $M$
is a finite cover of the same type whose Reeb foliation is again not
globally developable, and every quasi-regular soliton approximation of $M$
consists, for all large $i$, of solitons with nondevelopable leaf spaces. A
key ingredient is the fact that the Reeb flow of a soliton whose Reeb
closure is compact always has closed orbits (Corollary~\ref{cor:closed-orbit}).
\end{enumerate}
\end{remark}

In the nondevelopable case, we prove:

\begin{theorem}
\label{thm:nondev-structure-intro}Let $(M^{5},\xi ,\eta ,\Phi ,g,f)$ be a
complete noncompact gradient shrinking Sasaki--Ricci soliton whose Reeb
foliation is neither quasi-regular nor globally developable. Then the
following hold.

\begin{enumerate}
\item $M$ is irregular, with $H\cong T^{k}$, $k\in\{2,3\}$, and
$R^T>0$ on $M$, where $H$ is the closure of the Reeb flow in
$\operatorname{Isom}(M,g)$. On the universal cover, $\widehat H$ denotes
the closure of the lifted Reeb flow in the isometry group of the lifted
metric.

\item $\pi _{1}(M)$ is finite. The universal covering $\widehat{M}\rightarrow
M$ is finite and $M=\widehat{M}/\Gamma $ for a finite group $\Gamma \cong \pi
_{1}(M)$ of automorphisms of the lifted soliton which acts
freely and commutes with $\widehat{H}$.

\item If $(M_{i})$ is a quasi-regular soliton approximation of $M$, then for
all large $i$ the leaf space of $M_{i}$ is not developable, $R^{T}>0$ on
$M_{i}$, $\pi _{1}(M_{i})$ is finite, and the universal cover of $M_{i}$ is
quasi-regular, not regular, with a nondevelopable leaf space. Moreover the
generic Reeb periods of $M_{i}$ tend to infinity.
\end{enumerate}
\end{theorem}

This case is not empty: the solitons $\mathcal{S}_{a}$ with $a_{1}/a_{2}\notin
\mathbb{Q}$ of Example~\ref{ex:counterexample} in
Section~\ref{sec:nondev-models} belong to it.

\begin{theorem}
\label{thm:S3-family-intro} Let $\mathcal{S}_{a}$, $a=(a_{1},a_{2})$,
$a_{1},a_{2}>0$, $a_{1}+a_{2}=\lambda /2$, be the gradient shrinking
Sasaki--Ricci solitons on $\mathbb{C}\times \mathbb{S}^{3}$ of
Example~\ref{ex:counterexample} (Proposition~\ref{prop:S3-family}).

\begin{enumerate}
\item $\mathcal{S}_{a}$ is complete, has bounded curvature, and depends
real-analytically on $a$ in $C^{\infty }$ on compact subsets of
\begin{equation*}
\mathbb{C}\times \mathbb{S}^{3}.
\end{equation*}

\item If $a_{1}/a_{2}=m_{1}/m_{2}$ in lowest terms, then $\mathcal{S}_{a}$
is quasi-regular, with generic Reeb period $L=4\pi (m_{1}+m_{2})/\lambda $
and exceptional periods $L/m_{1}$ and $L/m_{2}$. Its leaf space is
\begin{equation*}
\mathbb{C}\times \mathbb{CP}_{[m_{1},m_{2}]}^{1},
\end{equation*}
which is not developable unless $m_{1}=m_{2}=1$; in that case
$\mathcal{S}_{a}$ is the regular cylinder lift with $d=1$.

\item If $a_{1}/a_{2}\notin \mathbb{Q}$, then $\mathcal{S}_{a}$ is irregular
with $H=T^{2}$. If $a^{(i)}\to a$ are positive pairs satisfying
$a_1^{(i)}+a_2^{(i)}=\lambda/2$ and
$a_1^{(i)}/a_2^{(i)}\in\mathbb Q$, then the solitons
$\mathcal S_{a^{(i)}}$ form a quasi-regular soliton approximation of
$\mathcal S_a$ on the fixed manifold $\mathbb C\times\mathbb S^3$;
for large $i$ it is nondevelopable, and its isotropy orders and generic
periods tend to infinity. The limit has positive nonconstant $R^{T}$.

\item $\mathcal{S}_{a}$ and $\mathcal{S}_{a^{\prime }}$ are isomorphic if
and only if $a^{\prime }=a$ or $a^{\prime }=(a_{2},a_{1})$. Hence, as
$a_{1}/a_{2}$ ranges over the irrational numbers larger than one, the
solitons $\mathcal{S}_{a}$ form a one-parameter family of pairwise
nonisomorphic, simply connected, irregular shrinkers with $R^{T}>0$.
\end{enumerate}
\end{theorem}

\begin{remark}
Theorem~\ref{thm:irregular-main} is sharp in the sense that, without
developability of the approximating solitons, a limit of quasi-regular
shrinkers which is not quasi-regular need not be transversely flat, and its
transverse scalar curvature may be nonconstant.
\end{remark}

\textbf{Global splitting: }When the transverse scalar curvature is bounded,
the local splitting theorem of \cite{PartI} shows that every developed
local limit at transverse infinity splits off a flat complex line.
We show that, up to finite quotients and the transversely flat models, the
solitons $\mathcal{S}_{a}$ are the only shrinkers which split globally
(Definition~\ref{def:split}), and that global splitting is characterized by
nonnegative transverse holomorphic bisectional curvature, in analogy with
the theorem of Naber \cite{Nab} and its extension by Munteanu--Wang
\cite{MWpos}.

\begin{theorem}
\label{thm:split-class-intro}Let $(M^{5},\xi ,\eta ,\Phi ,g,f)$ be a
connected complete gradient shrinking Sasaki--Ricci soliton with constant
$\lambda $, normalized by \eqref{eq:pre-normal}, which is globally split.

\begin{enumerate}
\item If $R^{T}\equiv 0$, then $M$ is isomorphic to $\mathbb{H}_{2}$ or to a
Gaussian suspension $M_{A,\ell }$; all of these are globally split.

\item Otherwise $R^{T}>0$, and
\begin{equation*}
M\cong \mathcal{S}_{a}/\Gamma
\end{equation*}
for a unique $a$ with $a_{1}\geq a_{2}>0$, $a_{1}+a_{2}=\lambda /2$, and a
finite group $\Gamma $ of automorphisms of $\mathcal{S}_{a}$ acting freely,
unique up to conjugation. If $a_{1}>a_{2}$, then $\Gamma $ is cyclic of some
order $p\geq 1$, generated by
\begin{equation*}
(z,w_{1},w_{2})\longmapsto (e^{2\pi ij/p}z,e^{2\pi i/p}w_{1},e^{2\pi
iq/p}w_{2}),\qquad \gcd (q,p)=1,\ j\in \mathbb{Z},
\end{equation*}
in the coordinates $w_{1}=\sqrt{1-t}\,e^{i\theta _{1}}$,
$w_2=\sqrt t\,e^{i\theta_2}$ of $\mathbb S^3$. Here $q$ and $j$ are
read modulo $p$. The lens space $L(p,q)$ is the quotient of $\mathbb S^3$
by the displayed action on $(w_1,w_2)$; this $q$ is an action weight,
not the Euler number of Theorem~\ref{thm:intro-topology}. The projection
$[z,w]\mapsto[w]$ makes $M$ a complex line bundle over $L(p,q)$ with a
flat unitary connection. Flatness here concerns that connection, not
the Sasakian metric on $M$. If $a_1=a_2$, then $M$ is a finite
quotient of the regular cylinder lift with $d=1$.
\end{enumerate}
\end{theorem}

\begin{theorem}
\label{thm:nonneg-op-intro}Let $(M^{5},\xi ,\eta ,\Phi ,g,f)$ be a connected
complete noncompact gradient shrinking Sasaki--Ricci soliton with
nonnegative transverse holomorphic bisectional curvature. Then $M$ is
globally split.
\end{theorem}

The proof of Theorem~\ref{thm:nonneg-op-intro} uses the compactness
criterion of Proposition~\ref{prop:pos-ricci}: a complete gradient shrinking
Sasaki--Ricci soliton with nonnegative transverse holomorphic bisectional
curvature and positive definite transverse Ricci tensor is compact. It is
the Sasakian counterpart of the compactness theorem of Munteanu--Wang
\cite{MWpos} for shrinkers with nonnegative sectional curvature and
positive Ricci curvature, and of its K\"{a}hler analogue due to Ni \cite{Ni} and Wu--Zhang \cite{WZh} (compare
Chang--Han--Wu \cite{CHW}), and it is used to force a zero eigenvalue of the
transverse Ricci tensor. We prove it in Section~\ref{sec:nonneg-bisec} by
the least-eigenvalue method of \cite{MWpos,WZh}. Since the Ricci curvature
of $g$ itself need not be nonnegative, the volume input is obtained from a
comparison for tubes about compact Reeb-closure orbits, proved in
Lemma~\ref{lem:transverse-CY}. This gives infinite total volume in the
noncompact
case. A least-eigenvalue barrier and a negative-gradient-flow argument force
the transverse scalar curvature to tend to infinity, contradicting the
integrated trace equation. The subsequent strong minimum principle gives a
zero Ricci eigenvalue everywhere. These arguments do not require a global
leaf space.

\begin{corollary}
\label{cor:nonneg-nondev-intro}Let $(M^{5},\xi ,\eta ,\Phi ,g,f)$ be a
complete noncompact gradient shrinking Sasaki--Ricci soliton whose Reeb
foliation is neither quasi-regular nor globally developable. The following
are equivalent:

\begin{enumerate}
\item[(a)] The transverse curvature operator of $M$ is nonnegative;

\item[(a$^{\prime }$)] The transverse holomorphic bisectional curvature of $M
$ is nonnegative;

\item[(b)] $M$ is globally split;

\item[(c)]
\begin{equation*}
M\cong \mathcal{S}_{a}/\Gamma
\end{equation*}
for a unique $a$ with $a_{1}>a_{2}$ and $a_{1}/a_{2}\notin \mathbb{Q}$, and
a cyclic group $\Gamma \cong \mathbb{Z}/p\mathbb{Z}$ as in
Theorem~\ref{thm:split-class}(2).
\end{enumerate}
\end{corollary}

\begin{corollary}
A connected complete five-dimensional gradient shrinking Sasaki--Ricci
soliton with constant $\lambda>0$, normalized by \eqref{eq:pre-normal},
whose Reeb foliation is not globally developable and which is globally
split is isomorphic to $\mathcal{S}_{a}/\Gamma $ for a unique $a$ with $a_{1}>a_{2}>0$,
$a_{1}+a_{2}=\lambda /2$, and a finite group $\Gamma $ of automorphisms of
$\mathcal{S}_{a}$ acting freely. Here $\mathcal{S}_{a}$ is quasi-regular if
$a_{1}/a_{2}\in \mathbb{Q}$ and irregular if $a_{1}/a_{2}\notin \mathbb{Q}$.
\end{corollary}

\begin{problem}
\label{prob:open}

\begin{enumerate}
\item Classify the nondevelopable quasi-regular five-dimensional shrinkers,
an orbifold version of \cite{LW,CDS,CCD,BCCD}. By
Corollary~\ref{cor:nondev-qr}, it suffices to treat simply connected total
spaces and
their finite free quotients.

\item Does every soliton whose Reeb closure is a torus of dimension at least
two (case $(T)$ of Proposition~\ref{prop:second-structure}) admit a
quasi-regular soliton approximation?

\item In the nondevelopable case of Theorem~\ref{thm:nondev-structure}, is
the transverse scalar curvature always bounded? (Compare \cite{LW,LZ} in the
K\"{a}hler case.)
\end{enumerate}
\end{problem}

\begin{remark}

\begin{enumerate}
\item The nondevelopable quasi-regular shrinkers are classified under global
splitting or nonnegativity of the transverse curvature, and for all solitons
exhibited in this paper quasi-regular soliton approximations exist and $R^{T}
$ is bounded. In general the three questions remain open. In contrast, in the
steady case there is no nondevelopable quasi-regular Sasaki--Ricci soliton
with nonnegative transverse holomorphic bisectional curvature; see \cite{CTZ}.

\item The questions of Problem~\ref{prob:open} concern the solitons of
Theorem~\ref{thm:nondev-structure} and the nondevelopable quasi-regular
solitons of Corollary~\ref{cor:nondev-qr}. So far the only known nonflat
examples in these classes are the globally split solitons
$\mathcal{S}_{a}/\Gamma $, which are \emph{cylindrical} at transverse
infinity: by \eqref{eq:ctr-RT}, $R^{T}$ is bounded below by the positive
constant $2\lambda c/(e^{c}-1)$, where $c=c(a)\geq 0$ is the constant of
\eqref{eq:Hc} (the bound reads $2\lambda$ for $c=0$).

\item In a forthcoming paper we shall construct, by a transverse Calabi
ansatz, a second family $\mathcal{F}_{b}$, again on $\mathbb{C}\times
\mathbb{S}^{3}$, which is not globally split and is \emph{conical}: $R^{T}>0$
decays like $1/f$. Its quasi-regular members are circle lifts of orbifold
versions of the FIK shrinker with a branch divisor along the exceptional
curve, and its irregular members are nondevelopable. These solitons, and the
partial answers to the three questions which they provide, are not used in
the present paper.
\end{enumerate}
\end{remark}

\bigskip

\textbf{Organization.} Section~\ref{sec:prelim} recalls Sasakian geometry,
Sasaki--Ricci solitons and the structure theorems for Reeb foliations.
Section~\ref{sec:local} recalls the results of \cite{PartI} used here.
Section~\ref{sec:structure} treats the regular and quasi-regular cases and
the proper $\mathbb{R}$-action. Section~\ref{sec:irregular} contains the
transversely flat models, the closed-orbit property, the quasi-regular
soliton approximations and the structure theorems for Reeb foliations which
are not quasi-regular. Section~\ref{sec:split} contains the global splitting
theorems. The computations for the explicit models and the proofs of
Section~\ref{sec:split} are collected in the appendices.

\section{Preliminaries}

\label{sec:prelim}

Sasakian geometry is the odd-dimensional analogue of K\"{a}hler geometry
(see \cite{BG}). In this section we recall the facts needed below and keep
the conventions of \cite{PartI}.

\subsection{Sasakian Manifolds and Sasakian Structures}

We recall the definition.

\begin{definition}
A Riemannian manifold $(M^{2n+1},g)$ is called a \emph{Sasakian manifold} if
its metric cone
\begin{equation*}
(C(M),\overline{g},J)=(M\times \mathbb{R}^{+},\,r^{2}g+dr^{2},\,J)
\end{equation*}
is K\"{a}hler. We identify $M$ with $\{r=1\}\subset C(M)$.
\end{definition}

The \emph{Reeb vector field} and the \emph{contact form} are
\begin{equation*}
\xi =J\Bigl(r\frac{\partial }{\partial r}\Bigr),\qquad \eta (Y)=g(\xi ,Y).
\end{equation*}
Then $\eta (\xi )=1$ and $\iota _{\xi }d\eta =0$, and $\xi $ is a Killing
field of unit length. The following conditions on a Riemannian manifold
$(M^{2n+1},g)$ are equivalent (see \cite{BG}):

\begin{enumerate}
\item there is a Killing vector field $\xi $ of unit length on $M$ such that
the tensor field $\Phi $ of type $(1,1)$ defined by $\Phi (Y)=\nabla _{Y}\xi $
satisfies
\begin{equation*}
(\nabla _{X}\Phi )(Y)=g(\xi ,Y)X-g(X,Y)\xi
\end{equation*}
for all vector fields $X$ and $Y$ on $M$;

\item there is a Killing vector field $\xi $ of unit length on $M$ such that
the Riemann curvature tensor satisfies
\begin{equation*}
R(X,\xi )Y=g(\xi ,Y)X-g(X,Y)\xi
\end{equation*}
for all vector fields $X$ and $Y$ on $M$;

\item the metric cone $(C(M),\overline{g})$ is K\"{a}hler.
\end{enumerate}

Moreover, there is a natural splitting
\begin{equation*}
TC(M)=L_{r\partial _{r}}\oplus L_{\xi }\oplus D,
\end{equation*}
where $L_{r\partial _{r}}$ and $L_{\xi }$ are the line bundles spanned by
$r\partial _{r}$ and $\xi $, and $D=\ker \eta $ is the contact distribution.
The complex structure $J$ of the cone and the tensor $\Phi $ are related by
\begin{equation*}
JY=\Phi (Y)-\eta (Y)\,r\frac{\partial }{\partial r},\qquad \Phi ^{2}=-I+\eta
\otimes \xi ,
\end{equation*}
for vector fields $Y$ on $M$, and
\begin{equation*}
g(\Phi X,\Phi Y)=g(X,Y)-\eta (X)\eta (Y).
\end{equation*}
Thus
\begin{equation*}
g=g^{T}+\eta \otimes \eta ,
\end{equation*}
where $g^{T}(X,Y)=g(\Phi X,\Phi Y)$ is the transverse metric, that is, the
restriction of $g$ to $D$ extended by zero in the direction of $\xi $. We call
$(M,\xi ,\eta ,\Phi ,g)$ a \emph{Sasakian structure}.

\subsection{Transverse Geometry and Transverse K\"{a}hler Structure}

\label{sec:transverse}

Let $\{U_{\alpha }\}_{\alpha \in A}$ be an open covering of the Sasakian
manifold $(M,\xi ,\eta ,\Phi ,g)$, and let
\begin{equation*}
\pi _{\alpha }:U_{\alpha }\rightarrow V_{\alpha }\subset \mathbb{C}^{n}
\end{equation*}
be foliation-chart submersions whose fibers are the connected plaques of
the Reeb foliation $\mathcal F_\xi$ in $U_\alpha$. After restricting to
sufficiently small components of an overlap, there are biholomorphic
transverse transition maps $\gamma_{\alpha\beta}$ with
$\pi_\alpha=\gamma_{\alpha\beta}\circ\pi_\beta$. These are local
transition maps: no inverse of a submersion $\pi_\beta$ and no global leaf
space are used.
For every $p\in U_{\alpha }$ the differential
\begin{equation*}
d\pi _{\alpha }:D_{p}\rightarrow T_{\pi _{\alpha }(p)}V_{\alpha }
\end{equation*}
is an isomorphism, where $D=\ker \eta \subset TM$. Since $\xi $ generates
isometries, the restriction of the Sasakian metric $g$ to $D$ induces a
well-defined Hermitian metric $g_{\alpha }^{T}$ on $V_{\alpha }$, and this
metric is K\"{a}hler. More precisely, let $z^{1},\ldots ,z^{n}$ be holomorphic
coordinates on $V_{\alpha }$, pulled back to $U_{\alpha }$ and denoted by the
same letters, and let $x$ be a coordinate along the leaves with $\xi =\partial
/\partial x$. Then $\{x,z^{1},\ldots ,z^{n}\}$ are foliated coordinates on
$U_{\alpha }$, and $D\otimes \mathbb{C}$ is spanned by the vector fields
\begin{equation*}
Z_{j}=\frac{\partial }{\partial z^{j}}-\eta \Bigl(\frac{\partial }{\partial
z^{j}}\Bigr)\frac{\partial }{\partial x},\qquad j=1,\ldots ,n,
\end{equation*}
and their conjugates. Since $\iota _{\xi }d\eta =0$,
\begin{equation*}
d\eta (Z_{j},\overline{Z_{k}})=d\eta \Bigl(\frac{\partial }{\partial
z^{j}},\frac{\partial }{\partial \overline{z}^{k}}\Bigr).
\end{equation*}
Hence the K\"{a}hler form $\omega _{\alpha }^{T}$ of $g_{\alpha }^{T}$ on
$V_{\alpha }$, which corresponds to the restriction of $\frac{1}{2}d\eta $ to
a slice $\{x=\mathrm{const}\}$ in $U_{\alpha }$, is closed. The collection of
K\"{a}hler metrics $\{g_{\alpha }^{T}\}$ on $\{V_{\alpha }\}$ is called a
\emph{transverse K\"{a}hler metric}; we write $g^{T}$, $\omega
^{T}=\frac{1}{2}d\eta $ and $J=\Phi |_{D}$ for the transverse metric, its
K\"{a}hler form and the transverse complex structure. For example, on the
Heisenberg group $\mathbb{H}_{n}=\mathbb{C}^{n}\times \mathbb{R}$ with
coordinates $(z,t)$, the contact form
\begin{equation*}
\eta
=dt+\frac{\sqrt{-1}}{2}\sum_{j}(z^{j}d\overline{z}^{j}-\overline{z}^{j}dz^{j})=dt+\sum_{j}(x_{j}dy_{j}-y_{j}dx_{j}),\qquad
z^{j}=x_{j}+\sqrt{-1}y_{j},
\end{equation*}
the Reeb field $\xi =\partial /\partial t$ and the frame $Z_{j}=\partial
/\partial z^{j}+\frac{\sqrt{-1}}{2}\overline{z}^{j}\partial /\partial t$ give
the flat model: a transverse K\"{a}hler potential is $\frac{1}{2}|z|^{2}$, the
transverse metric is Euclidean, and $d\eta =2\omega _{\mathrm{E}}$ with
$\omega _{\mathrm{E}}=\sum_{j}dx_{j}\wedge dy_{j}$. This is the Heisenberg
model \eqref{eq:heisenberg} used below.

Let $\widetilde{X},\widetilde{Y},\widetilde{W},\widetilde{Z}\in \Gamma (D)$
correspond under $d\pi _{\alpha }$ to $X,Y,W,Z\in \Gamma (TV_{\alpha })$. The
Levi-Civita connection $\nabla ^{T}$ of the transverse K\"{a}hler metric
$g^{T}$, its curvature and its Ricci tensor are given by
\begin{equation*}
\begin{array}{rcl}
\nabla _{X}^{T}Y & := & d\pi _{\alpha }(\nabla _{\widetilde{X}}\widetilde{Y}),
\\
\widetilde{\nabla _{X}^{T}Y} & = & \nabla _{\widetilde{X}}\widetilde{Y}+g(\Phi
\widetilde{X},\widetilde{Y})\xi , \\
\mathrm{Rm}^{T}(X,Y,Z,W) & = &
\mathrm{Rm}_{g}(\widetilde{X},\widetilde{Y},\widetilde{Z},\widetilde{W})+g(\Phi
\widetilde{Y},\widetilde{W})g(\Phi \widetilde{X},\widetilde{Z}) \\
& & -g(\Phi \widetilde{X},\widetilde{W})g(\Phi
\widetilde{Y},\widetilde{Z})+2g(\Phi \widetilde{X},\widetilde{Y})g(\Phi
\widetilde{Z},\widetilde{W}), \\
\mathrm{Ric}^{T}(X,Z) & = &
\mathrm{Ric}_{g}(\widetilde{X},\widetilde{Z})+2g(\widetilde{X},\widetilde{Z}).
\end{array}
\end{equation*}
In the last line $\mathrm{Ric}_{g}$ is the Ricci tensor of $g$; the
contraction of $\mathrm{Rm}_{g}$ over $D$ alone is
$\mathrm{Ric}_{g}|_{D}-g^{T}$, by \eqref{eqn-SR3} below. Here and throughout
we use the conventions
\begin{equation}
\begin{array}{rcl}
R(X,Y)Z & := & \nabla _{X}\nabla _{Y}Z-\nabla _{Y}\nabla _{X}Z-\nabla
_{\lbrack X,Y]}Z, \\
R(X,Y,Z,W) & := & \langle R(Z,W)Y,X\rangle , \\
R_{\alpha \beta \gamma \delta } & := & R(e_{\alpha },e_{\beta },e_{\gamma
},e_{\delta })=\langle R(e_{\gamma },e_{\delta })e_{\beta },e_{\alpha }\rangle
, \\
\mathrm{Ric}(X,Y) & := & \sum_{\alpha }R(X,e_{\alpha },Y,e_{\alpha }),
\end{array}
\label{eqn-SR1}
\end{equation}
where $\{e_{\alpha }\}$ is a local orthonormal frame; thus $R_{\alpha \beta
\alpha \beta }$ is a sectional curvature. Let $\{e_{i}\}_{i=1}^{2n}$ be a
local orthonormal frame of $D=\ker \eta $. On a Sasakian manifold
\begin{equation}
R(X,\xi )Y=g(\xi ,Y)X-g(X,Y)\xi ,  \label{eqn-SR2}
\end{equation}
and therefore
\begin{equation}
\mathrm{Ric}^{T}=\mathrm{Ric}_{g}+2g^{T}\qquad \text{on }D.  \label{eq:pre-ric}
\end{equation}
A function $v$ is basic if $\xi v=0$. For basic functions we use the real
transverse operators
\[
\operatorname{Hess}^T v=(\nabla^T)^2v,\qquad
\Delta_Bv=\operatorname{tr}_{g^T}\operatorname{Hess}^Tv,\qquad
L_fv=\Delta_Bv-\langle\nabla^Tf,\nabla^Tv\rangle_{g^T}.
\]
The trace is over a real orthonormal frame of $D$ of size $2n$; in
particular $R^T=\operatorname{tr}_{g^T}\mathrm{Ric}^T$ is the real
transverse scalar curvature. On basic functions $\nabla^Tv=\nabla_gv$
and $\Delta_Bv=\Delta_gv$.

For $\tau >0$, the \emph{transverse $D$-homothetic deformation} (see \cite{BG})
\begin{equation*}
\eta _{\tau }=\tau \eta ,\qquad \xi _{\tau }=\tau ^{-1}\xi ,\qquad \Phi _{\tau
}=\Phi ,\qquad g_{\tau }=\tau g+\tau (\tau -1)\eta \otimes \eta =\tau
g^{T}+\tau ^{2}\eta \otimes \eta
\end{equation*}
is a Sasakian structure with $g_{\tau }^{T}=\tau g^{T}$. It preserves $\nabla
^{T}$, $\mathrm{Ric}^{T}$ and $\operatorname{Hess}^{T}f$ as $(0,2)$-tensors,
while $R^{T}$ and $|\nabla ^{T}f|^{2}$ are multiplied by $\tau ^{-1}$.
Consequently, if $(g,f)$ satisfies \eqref{1} with constant $\lambda $, then
$(g_{\tau },f)$ satisfies it with constant $\lambda/\tau$ (see
\cite[Remark~\Iref{rem:Dhom}{2.1}]{PartI}). In the quasi-regular case the
generic Reeb period becomes $L_\tau=\tau L$. Thus $\lambda L$ and the
integral circle-orbibundle class are unchanged by this deformation;
in particular the integers $4\pi/(\lambda L)$ and $8\pi/(\lambda L)$,
when defined by the lifting conditions, are scale invariant.

\subsection{Sasaki--Ricci Solitons}

\begin{definition}
A Sasakian manifold $(M,\xi ,\eta ,\Phi ,g)$ is \emph{$\eta $-Einstein} if
there is a constant $A$ such that
\begin{equation*}
\mathrm{Ric}_{g}=Ag+(2n-A)\eta \otimes \eta .
\end{equation*}
\end{definition}

Since $g=g^{T}+\eta \otimes \eta $, the $\eta $-Einstein condition is
equivalent to
\begin{equation*}
\mathrm{Ric}^{T}=(A+2)g^{T}\quad \text{on }D=\ker \eta ,\qquad \text{that
is,}\qquad \rho ^{T}=(A+2)\omega ^{T},
\end{equation*}
where $\rho ^{T}=\mathrm{Ric}^{T}(\Phi \,\cdot ,\cdot )$ is the transverse
Ricci form; such a structure is also called \emph{transverse
K\"{a}hler--Einstein}. For $A=2n$ one has
\begin{equation*}
\mathrm{Ric}_{g}=2ng,
\end{equation*}
and the structure is \emph{Sasaki--Einstein}. In summary, $(M,g)$ is
Sasaki--Einstein if and only if the K\"{a}hler cone $(C(M),\overline{g})$ is
Ricci-flat, and if and only if the transverse K\"{a}hler structure of the Reeb
foliation $\mathcal{F}_{\xi }$ is K\"{a}hler--Einstein with
\begin{equation*}
\mathrm{Ric}^{T}=2(n+1)g^{T}.
\end{equation*}
We record the following curvature identities on a Sasakian manifold. By
\eqref{eqn-SR2},
\begin{equation}
R(e_{i},\xi ,e_{j},\xi )=g(e_{i},e_{j})=\delta _{ij}.  \label{eqn-SR3}
\end{equation}
Therefore
\begin{equation}
\mathrm{Ric}_{g}(\xi ,\xi )=\sum_{i=1}^{2n}R(e_{i},\xi ,e_{i},\xi )=2n.
\label{eqn-SR4}
\end{equation}
Consequently
\begin{equation*}
R^{T}=R_{g}+2n
\end{equation*}
and
\begin{equation}
\sum_{i=1}^{2n}\mathrm{Ric}_{g}(e_{i},e_{i})=R_{g}-2n=R^{T}-4n.  \label{eqn-SR6}
\end{equation}

\begin{definition}[{\cite{FOW,CLL}}]
Let $\mathcal X$ be a Reeb-invariant complex vector field whose
projection to the complexified normal bundle of the Reeb foliation is
of type $(1,0)$ and transversely holomorphic. It is \emph{Hamiltonian
holomorphic} with basic, possibly complex-valued, potential $f$ if
\begin{equation}
f=\sqrt{-1}\,\eta(\mathcal X),\qquad
\iota_{\mathcal X}\omega^T=\sqrt{-1}\,\overline\partial_B f.
\label{pr1}
\end{equation}
Here $\omega^T=\tfrac12d\eta=g^T(\Phi\,\cdot,\cdot)$,
$\overline\partial_B$ is the transverse $(0,1)$ differential on basic
functions, and $\iota$ denotes contraction. Put
$X=\mathcal X+\overline{\mathcal X}=2\operatorname{Re}\mathcal X$.
The Sasakian structure is a \emph{Sasaki--Ricci soliton} if the real
field $X$ satisfies
\[
\mathrm{Ric}^T+\tfrac12\mathcal L_Xg^T=(A+2)g^T
\]
for a real constant $A$. The soliton is expanding, steady or shrinking
according as $A<-2$, $A=-2$ or $A>-2$, respectively.

In the gradient case $f$ is real. The Hamiltonian identities give
\[
\mathcal X=(\nabla^Tf)^{1,0}-\sqrt{-1}f\xi,\qquad
(\nabla^Tf)^{1,0}=\tfrac12(\nabla^Tf-\sqrt{-1}J\nabla^Tf),\qquad
X=\nabla^Tf.
\]
Consequently $\tfrac12\mathcal L_Xg^T=\operatorname{Hess}^Tf$ and the
full real equation is
$\mathrm{Ric}^T+\operatorname{Hess}^Tf=(A+2)g^T$.
The Lie derivative always uses the real field $X$; $\mathcal X$ denotes
its complex Hamiltonian lift. Its $(1,1)$ components read
\begin{equation*}
R_{j\overline{k}}^T+f_{j\overline{k}}=(A+2)g_{j\overline{k}}^T,\qquad
\text{equivalently}\qquad
R_{j\overline{k}}+f_{j\overline{k}}=Ag_{j\overline{k}}^T,
\end{equation*}
where $R_{j\overline{k}}$ are the components of
$\mathrm{Ric}_g|_D$. The full real equation also includes its pure-type
components; the displayed $(1,1)$ equation alone is not the definition.
In the shrinking case we write $\lambda=A+2>0$.
At any scale, the transverse trace and divergence of the soliton equation give
(see \cite[Lemma~\Iref{lem:identities}{3.1}]{PartI})
\begin{equation}
R^{T}+\Delta _{B}f=2n\lambda ,\quad \nabla ^{T}R^{T}=2\mathrm{Ric}^{T}(\nabla
^{T}f,\cdot ),\quad R^{T}+|\nabla ^{T}f|^{2}-2\lambda f=\mathrm{const}.
\label{eq:pre-identities}
\end{equation}
We normalize the additive constant in $f$ so that
\begin{equation}
R^{T}+|\nabla ^{T}f|^{2}=2\lambda f,  \label{eq:pre-normal}
\end{equation}
and call the soliton \emph{normalized}; \eqref{eq:pre-normal} is the
normalized Hamilton identity.
\end{definition}

\begin{example}
Unlike compact gradient K\"{a}hler--Ricci shrinker surfaces, which are
classified (see \cite{Cao} and the references therein), compact
five-dimensional gradient shrinking Sasaki--Ricci solitons are not classified.
Examples include:

\begin{enumerate}
\item the regular Sasaki--Einstein metrics on $\mathbb{S}^{5}$,
$\mathbb{S}^{2}\times \mathbb{S}^{3}$ and $\#k(\mathbb{S}^{2}\times
\mathbb{S}^{3})$, $3\leq k\leq 8$, obtained as circle bundles over the
K\"{a}hler--Einstein del Pezzo surfaces \cite{BG};

\item the toric Sasaki--Einstein metrics on $\#k(\mathbb{S}^{2}\times
\mathbb{S}^{3})$ for every $k\geq 1$, most of which are quasi-regular or
irregular \cite{FOW}; see \cite{BG} and the references therein;

\item the regular lifts to $\mathbb{S}^{2}\times \mathbb{S}^{3}$ and to
$\#2(\mathbb{S}^{2}\times \mathbb{S}^{3})$ of the Koiso--Cao shrinker on
$\mathbb{CP}^{2}\#\overline{\mathbb{CP}}^{2}$ \cite{Koi,Cao96} (see also
\cite{Cao}) and of the Wang--Zhu shrinker on
$\mathbb{CP}^{2}\#2\overline{\mathbb{CP}}^{2}$ \cite{WZ}, respectively.
\end{enumerate}
\end{example}

\subsection{Types of Reeb Foliations and Structure Theorems}

\label{sec:types}

A Sasakian structure is \emph{quasi-regular} if all leaves of the Reeb
foliation $\mathcal{F}_{\xi }$ are compact, so that $\xi $ generates a locally
free circle action. It is \emph{regular} if this action is free, that is, the
Reeb flow is a free effective circle action with a common minimal period
$L>0$. We call the structure \emph{irregular} if some Reeb leaf is not closed.
On a noncompact manifold, there is also the possibility of a free proper
$\mathbb{R}$-action which has a smooth global quotient and is treated
separately. In that case every leaf is a closed but noncompact line, so the
structure is neither quasi-regular nor irregular; the phrase ``not
quasi-regular'' always includes this case. In the quasi-regular case,
the generic period $L$ is the least positive time for which the Reeb
flow is the identity on $M$. It is the minimal orbit period on the
principal stratum; an orbit with stabilizer of order $m$ has minimal
period $L/m$.

We recall the \emph{first structure theorem} for complete Sasakian manifolds.

\begin{proposition}[{\cite{BG,Kol2}}]
In the quasi-regular case the leaf space $\mathcal{Z}=M/\mathbb{S}^{1}$ is a
K\"{a}hler orbifold and $\pi :M\rightarrow \mathcal{Z}$ is a
Seifert circle orbibundle. In the regular case $Z=M/\mathbb{S}^{1}$ is a
smooth K\"{a}hler manifold and $\pi $ is a principal circle bundle. In both
cases
\begin{equation}
\frac{1}{2}d\eta =\pi ^{\ast }\omega _{Z},\qquad g=\pi ^{\ast }h_{Z}+\eta
\otimes \eta ,\qquad f=\pi ^{\ast }f_{Z},  \label{eq:pre-BW}
\end{equation}
and the gradient Sasaki--Ricci soliton descends to the gradient
K\"{a}hler--Ricci shrinker equation
\begin{equation*}
\mathrm{Ric}_{Z}+\operatorname{Hess}_{Z}f_{Z}=\lambda h_{Z}
\end{equation*}
on the leaf space with the same shrinking constant $\lambda>0$. After a
constant transverse $D$-homothety one may normalize $\lambda=2n+2$. Indeed,
under $h_Z\mapsto \tau h_Z$, the tensors $\mathrm{Ric}_Z$ and
$\operatorname{Hess}_Z
f_Z$ are unchanged as $(0,2)$ tensors, so the shrinking constant changes to
$\lambda/\tau$. Circle bundles over a smooth base are classified by their
integral Euler class \cite{Ste,Hat}, and the orbifold version is described
in \cite{Kol2,Kol1}. Without any regularity assumption, either the closure
of the Reeb flow in the isometry group of $(M,g)$ is a compact torus, or $\xi
$ generates a closed, free and proper $\mathbb{R}$-action, in which case
$Z=M/\mathbb{R}$ is a smooth complete K\"{a}hler manifold; this case
includes the Heisenberg model.
\end{proposition}

Next we recall the \emph{second structure theorem} for complete Sasakian
manifolds.

Let $(M^{2n+1},\xi ,\eta ,\Phi ,g)$ be a connected complete Sasakian
manifold, and let $H$ be the closure of the Reeb flow $\{\varphi
_{s}\}_{s\in \mathbb{R}}$ in $\operatorname{Isom}(M,g)$. By Myers--Steenrod
\cite{MyS}, $H$ is a closed connected abelian Lie subgroup. When $H$ is
compact,
its orbits are the closures of the Reeb leaves. We regard its Lie algebra
$\mathfrak{h}$ as a space of Killing fields on $M$ and write
\begin{equation}
\mu :M\longrightarrow \mathfrak{h}^{\ast },\qquad \langle \mu ,\zeta \rangle
=\mu _{\zeta }:=\eta (\zeta ),  \label{eq:contact-moment}
\end{equation}
for the contact moment map. An element $\zeta \in \mathfrak{h}$ is
\emph{rational} if the one-parameter group it generates is closed in $H$.

\begin{proposition}[{\cite{Mol,BG,MyS,Wad}}]
\label{prop:second-structure}
Exactly one of the following holds.

\begin{enumerate}
\item $(Q)$\quad $H\cong \mathbb{S}^{1}$: The Reeb flow is periodic and the
structure
is quasi-regular. Conversely, if every Reeb leaf is compact, then it holds;
in particular the minimal periods of the leaves are then bounded (Wadsley
\cite{Wad}).

\item $(P)$\quad $H\cong \mathbb{R}:$ The Reeb flow is a free proper
$\mathbb{R}$-action. Every leaf is a closed noncompact line, so the structure
is neither quasi-regular nor irregular.

\item $(T)$\quad $H\cong T^{k}$ with $2\leq k\leq n+1$: The structure is
irregular. The union of the principal $H$-orbits is open, dense and
$H$-invariant, and the Reeb leaf through each of its points $x$ is dense in
the
$k$-dimensional torus $Hx$ and is not closed.

\begin{enumerate}
\item Every $\zeta \in \mathfrak{h}$ is a Killing field which commutes with
$\xi $ and preserves $\eta $, $\Phi $ and every Reeb-invariant continuous
tensor field. In particular it preserves the potential of a gradient
Sasaki--Ricci soliton.

\item Each $\mu _{\zeta }$ is a basic $H$-invariant function, $\iota _{\zeta
}d\eta =-d\mu _{\zeta }$, and the $H$-orbits are isotropic for $d\eta $.

\item Let $\mathfrak{h}^{+}=\{\zeta \in \mathfrak{h}:\inf_{M}\mu _{\zeta
}>0\}$, a convex cone containing $\xi $. For $\zeta \in \mathfrak{h}^{+}$
the type-I deformation
\begin{equation}
\eta _{\zeta }=\frac{\eta }{\mu _{\zeta }},\qquad \Phi _{\zeta }=\Phi -\Phi
\zeta \otimes \eta _{\zeta },\qquad g_{\zeta }=\tfrac{1}{2}d\eta _{\zeta
}(\cdot ,\Phi _{\zeta }\cdot )+\eta _{\zeta }\otimes \eta _{\zeta }
\label{eq:typeI}
\end{equation}
is a Sasakian structure with Reeb field $\zeta $ and with the same CR
structure $(D,\Phi |_{D})$. It is quasi-regular if and only if $\zeta $ is
rational. If $v=X+a\xi $ with $X\in D$, and $\zeta _{D}$ denotes the
$D$-component of $\zeta $, then
\begin{equation}
g_{\zeta }(v,v)=\mu _{\zeta }^{-1}\bigl|X-a\mu _{\zeta }^{-1}\zeta
_{D}\bigr|_{g}^{2}+a^{2}\mu _{\zeta }^{-2},\qquad 2|\zeta _{D}|=|\nabla
^{T}\mu
_{\zeta }|.  \label{eq:typeI-metric}
\end{equation}

\item If $\sup_M\|\mu\|<\infty$ for a fixed norm on
$\mathfrak h^*$, then $\mathfrak h^+$ is open. There are rational
$\zeta_i\in\mathfrak h^+$ with $\zeta_i\to\xi$, and the structures
\eqref{eq:typeI} converge to $(\xi,\eta,\Phi,g)$ in $C^\infty$ on
compact subsets. If, in addition,
$\sup_M|\nabla^T\mu_{\zeta_a}|<\infty$ for each member $\zeta_a$ of
one fixed basis of $\mathfrak h$, then there is a common constant $C\geq1$ with
$C^{-1}g\leq g_{\zeta_i}\leq Cg$ for every $i$. In particular these approximating
metrics are complete.
\end{enumerate}
\end{enumerate}
\end{proposition}

\section{Local Results from Part I}

\label{sec:local}

In this section we recall, without proofs, the results of \cite{PartI} that
are used below. They hold without any regularity assumption on the Reeb
foliation. We keep the normalizations of \cite{PartI}.

\subsection{Reeb Foliation and Quadratic Growth}

In \cite[Sections~\Iref{sec:identities}{3} and~\Iref{sec:curvature}{4}]{PartI}
the soliton equation is used in the normalization
\begin{equation}
\mathrm{Ric}^{T}+\operatorname{Hess}^{T}f=\lambda h,\qquad \lambda
=2n+2,\qquad \xi
f=0,\qquad h=g|_{D},\quad D=\ker \eta .  \label{eq:soliton}
\end{equation}
It is a full real symmetric-tensor equation.

\begin{lemma}[{\cite[Lemma~\Iref{lem:reeb}{3.2}]{PartI}}]
\label{lem:reeb} Let $M$ be a
connected and complete Sasakian manifold.

\begin{enumerate}
\item Either the closure of the Reeb one-parameter group in
$\operatorname{Isom}(M,g) $ is a compact torus $H$, or

\item the Reeb group is a closed, free and proper $\mathbb{R}$-action.
\end{enumerate}

In the first case put $\mathcal{X}=M$, choose a compact orbit $P=Hp$ and let
$r=d(P,\cdot )$; this function is basic. In the second case put
$\mathcal{X}=Z=M/\mathbb{R}$, which is a smooth complete K\"{a}hler manifold,
and let $r$
be the distance on $Z$ from a fixed point. In both cases $r$ is proper on
$\mathcal{X}$.
\end{lemma}

\begin{lemma}[{\cite[Lemma~\Iref{lem:nonneg}{3.3}]{PartI}}]
\label{lem:nonneg}For $n=2$, the transverse scalar curvature of any complete noncompact gradient
shrinking Sasaki--Ricci soliton $(M^{5},\xi ,\eta ,\Phi ,g,f)$ is
nonnegative.
Either $R^{T}>0$ everywhere, or $R^{T}\equiv 0$. Normalize the additive
constant in $f$ by
\begin{equation}
R^{T}+|\nabla ^{T}f|^{2}=12f.  \label{eq:canonical}
\end{equation}
Then
\begin{equation*}
f\geq 0,
\end{equation*}
and
\begin{equation}
L_{f}f=24-12f,\qquad |\nabla ^{T}\sqrt{f}|\leq \sqrt{3}.  \label{eq:Lf-f}
\end{equation}
The second inequality holds in the Lipschitz sense also at zeros of $f$. (In
\cite{PartI} the lemma is proved in every dimension $2n+1$ and at every
$\lambda >0$, with $L_{f}f=2n\lambda -2\lambda f$ and $|\nabla
^{T}\sqrt{f}|\leq \sqrt{\lambda /2}$; it is used for $n=1$ in
Appendix~\ref{app:B}.)
\end{lemma}

\begin{lemma}[{\cite[Lemma~\Iref{lem:growth}{3.4}]{PartI}}]
\label{lem:growth}Let $(M^{5},\xi ,\eta ,\Phi ,g,f)$ be a connected complete
noncompact
gradient shrinking Sasaki--Ricci soliton satisfying \eqref{eq:soliton}, with
$f$ normalized by \eqref{eq:canonical}, and let $\mathcal{X}$ and $r$ be as
in Lemma~\ref{lem:reeb}. There are constants $C$ and $a$ such that
\begin{equation}
\frac{(ar-C)_{+}^{2}}{2\lambda }\leq f\leq \frac{(\lambda r+C)^{2}}{2\lambda
}\qquad \text{on }\mathcal{X}.  \label{eq:potential-growth}
\end{equation}
In particular, $f$ is proper on $\mathcal{X}$. In the compact-torus case it
is proper on $M$. In the proper $\mathbb{R}$ case, its sublevels need not be
compact on $M$. Here $a=4$ for $\mathcal{X}=M$ and $a=\lambda =6$ for
$\mathcal{X}=Z.$
\end{lemma}

\subsection{Local Splitting}

\begin{theorem}[{\cite[Proposition~\Iref{prop:flat}{3.2} and
Theorem~\Iref{thm:curvature}{4.1}]{PartI}}]
\label{thm:intro-curvature} Let $(M^{5},\xi ,\eta ,\Phi ,g,f)$ be a connected,
complete, noncompact gradient
shrinking Sasaki--Ricci soliton with $\lambda =6$, with $f$ normalized by
\eqref{eq:canonical}. If the transverse scalar curvature $R^{T}$ is bounded
from above, then the transverse curvature
tensor $\mathrm{Rm}^{T}$ and the curvature tensor $\mathrm{Rm}_{g}$ of $g$
are bounded. More precisely,

\begin{enumerate}
\item either $R^{T}\equiv 0$ and $\mathrm{Rm}^{T}\equiv 0$, or

\item $R^{T}>0$ on $M$, and there are constants $C,c>0$ such that
\begin{equation*}
R^{T}f\geq c,\qquad |\mathrm{Rm}^{T}|+|\mathrm{Ric}^{T}|\leq CR^{T},\qquad |\mathrm{Rm}_{g}|\leq
C(1+R^{T})=C(5+R_{g})
\end{equation*}
on $M$.
\end{enumerate}
\end{theorem}

\begin{proposition}[{\cite[Proposition~\Iref{prop:derivatives}{5.2}]{PartI}}]
\label{prop:derivatives}Let $(M^{5},\xi ,\eta ,\Phi ,g,f)$ be a complete
gradient shrinking Sasaki--Ricci soliton, normalized by a $D$-homothety so
that $\lambda =\frac{1}{2}$. If the transverse scalar curvature $R^{T}$ is
bounded from above, then $\mathrm{Rm}^{T}$ is bounded by
Theorem~\ref{thm:intro-curvature}, and the associated ancient flow $g^{T}(t)$ of
\cite[Proposition~\Iref{lem:selfsim}{5.1}]{PartI} satisfies
\begin{equation*}
|(\nabla ^{T})^{k}\mathrm{Rm}^{T}(g^{T}(t))|\leq C_{k}|t|^{-1-k/2}.
\end{equation*}
\end{proposition}

\begin{proposition}[{\cite[Proposition~\Iref{prop:unwrap}{5.3}]{PartI}}]
\label{prop:unwrap} Let $(M,\xi ,\eta ,\Phi ,g,f)$ be a complete gradient
shrinking Sasaki--Ricci
soliton with bounded transverse curvature. Then $\sup_{M}|\nabla
^{k}\mathrm{Rm}_{g}|<\infty $ for every $k\geq 0$, and every sequence of
basepoints
admits developed local limits in the sense of
\cite[Definition~\Iref{def:developed}{5.2}]{PartI}: along a subsequence, the
structures
pulled back by the exponential maps to the tangent balls of a fixed radius
converge smoothly on the balls of half that radius.
\end{proposition}

\begin{theorem}[{\cite[Theorem~\Iref{thm:splitting}{6.8}]{PartI}; compare
\cite{Nab}}]
\label{thm:splitting} Let $p_{i}\in M$ with $f(p_{i})\rightarrow \infty $ on a
complete five-dimensional gradient shrinking Sasaki--Ricci soliton $M$ with
bounded transverse curvature (for instance, with $R^{T}$ bounded above, by
Theorem~\ref{thm:intro-curvature}). Every developed local transverse limit
at these basepoints has orthonormal parallel fields $V$, $JV$. Consequently
it splits locally as
\begin{equation}
(U,h_{\mathrm{E}})\times (\Sigma ,h_{\Sigma }),\qquad U\subset \mathbb{C},
\label{eq:local-split}
\end{equation}
with its product K\"{a}hler structure. Suppose a nonflat product limit
satisfies the local soliton equation. That is, on a neighborhood containing
its product chart, there is a smooth function $F$ with
\begin{equation}
\mathrm{Ric}(h_{\lim })+\operatorname{Hess}F=\lambda h_{\lim }.  \label{LS}
\end{equation}
Then
\begin{equation*}
F(z,y)=\frac{\lambda }{2}|z|^{2}+\ell (z)+\psi (y),\qquad
\operatorname{Hess}_{h_{\Sigma }}\psi =(\lambda -K)h_{\Sigma },
\end{equation*}
where $\ell $ is affine and $K$ is the Gaussian curvature of $h_{\Sigma }$.
In particular,
\begin{equation}
dK=K\,d\psi .  \label{eq:dK}
\end{equation}
If its surface factor, together with its potential, extends to a complete
smooth connected two-dimensional gradient shrinker with the same constant
$\lambda $, then its local metric has
\begin{equation*}
K\equiv \lambda .
\end{equation*}
It is locally the round sphere factor of $\mathbb{C}\times \mathbb{CP}^{1}$.
More generally, when a complete smooth global K\"{a}hler limit is available,
the parallel fields give a global splitting
\begin{equation*}
(\mathbb{C},h_{\mathrm{E}})\times (N,k)
\end{equation*}
of its simply connected cover.
\end{theorem}

\subsection{Transverse Structure at Infinity}

\label{sec:partI-end}

Let $(M^{5},\xi ,\eta ,\Phi ,g,f)$ be connected, complete and noncompact,
and satisfy the full real transverse equation
\begin{equation*}
Q+\operatorname{Hess}^{T}f=6h,\qquad \xi f=0,\qquad h=g|_{D},\quad D=\ker \eta .
\end{equation*}
Write
\begin{equation*}
T=\mathrm{Rm}^{T},\quad Q=\mathrm{Ric}^{T},\quad R^{T}=R_{g}+4,\qquad \omega
^{T}=\tfrac{1}{2}\,d\eta .
\end{equation*}
Choose the compatible sign for the transverse complex structure $J$ so $\omega
^{T}(U,V)=h(JU,V)$. Our curvature components are $T_{ijkl}=\left\langle
R^{T}(e_{k},e_{l})e_{j},e_{i}\right\rangle $. All lower-case indices
in this subsection run from one to four.

We work on the exhaustion space $\mathcal{X}$ of Lemma~\ref{lem:reeb}. When
the Reeb closure is compact, infinity is measured on $M$ by distance to a
fixed compact orbit. When the action is proper $\mathbb{R}$, it is measured
on the complete quotient $Z$. The normalized potential is proper on
$\mathcal{X}$ (Lemma~\ref{lem:growth}). Assume $R^{T}$ is bounded above. Either $R^{T}=T=0$ or
$R^{T}>0$ everywhere. In the nonflat case,
\begin{equation}
c(1+r)^{2}\leq f\leq C(1+r)^{2}\text{.}  \label{eq:input-growth}
\end{equation}

We retain the distance function $r$ from Lemma~\ref{lem:reeb}: it is
distance to a fixed compact Reeb-closure orbit when $\mathcal X=M$,
and distance to a fixed point when $\mathcal X=Z$. In the first case
it differs from distance to a chosen point of that orbit by a bounded
amount. Constants may depend on a fixed compact sublevel and the soliton,
not only on $\sup R^T$.

Fix $F>\max\{0,\sup R^T/6\}$. Then
$|\nabla^Tf|^2=12f-R^T\geq6f>0$ on $\{f\geq F\}$. On $\{f\geq F\}$ define the $J$-invariant real
two-plane
\begin{equation*}
E=\operatorname{span}\{\nabla ^{T}f,J\nabla ^{T}f\}^{\perp }\subset D,
\end{equation*}
and let $P_{E}$ be the orthogonal projection of $\Lambda ^{2}D$ onto its
one-dimensional subspace $\Lambda ^{2}E$. The transverse curvature operator
$\mathcal{R}^{T}$ is normalized so that a surface of Gaussian curvature $K$
has curvature operator $K$ on its oriented area form.

\begin{theorem}[{\cite[Theorem~\Iref{thm:general}{6.9}]{PartI}}]
\label{thm:general}Let $(M^{5},\xi ,\eta ,\Phi ,g,f)$ be as above, with
$R^{T}$ bounded from above. In the nonflat case one has, on $\{f\geq F\}$,
\begin{equation}
\Bigl\Vert\mathcal{R}^{T}-\frac{R^{T}}{2}P_{E}\Bigr\Vert_{\mathrm{op}}\leq
\frac{CR^{T}}{\sqrt{f}}\leq C\sqrt{\frac{R^{T}}{f}}\leq \frac{C}{1+r}.
\label{eq:rankbounded}
\end{equation}
If in addition $R^{T}\rightarrow 0$ uniformly as $f\rightarrow \infty $,
then
\begin{equation}
\frac{c}{f}\leq R^{T}\leq \frac{C}{f},\qquad |(\nabla ^{T})^{j}T|\leq
C_{j}f^{-1-j/2}\quad (j\geq 0),\qquad
\Bigl\Vert\mathcal{R}^{T}-\frac{R^{T}}{2}P_{E}\Bigr\Vert_{\mathrm{op}}\leq
Cf^{-2}.  \label{eq:decay-f}
\end{equation}
Equivalently,
\begin{equation*}
c(1+r)^{-2}\leq R^{T}\leq C(1+r)^{-2},
\end{equation*}
and
\begin{equation*}
|(\nabla ^{T})^{j}T|\leq C_{j}(1+r)^{-2-j},
\end{equation*}
and
\begin{equation}
\Bigl\Vert\mathcal{R}^{T}-\frac{R^{T}}{2}P_{E}\Bigr\Vert_{\mathrm{op}}\leq
C(1+r)^{-4}.
\end{equation}
In particular the transverse curvature operator is bounded below by
$-CR^{T}(1+r)^{-1}\mathrm{Id}$, hence by
$-C(1+r)^{-1}\mathrm{Id}$, in the bounded-scalar case and by
$-C(1+r)^{-4}\mathrm{Id}$ in the decaying branch. Every developed local
transverse limit along a
sequence $p_{i}$ with $f(p_{i})\rightarrow \infty $ consequently has
nonnegative curvature operator. Its surface factor in \eqref{eq:local-split}
has $K\geq 0$.
\end{theorem}

\section{Structure Theorems for Quasi-Regular Sasaki--Ricci Shrinkers}

\label{sec:structure}

In this section we derive the structure at infinity. The statements of
Section~\ref{sec:local}, like those of
\cite[Section~\Iref{sec:end}{6.5}]{PartI} recalled in
Section~\ref{sec:partI-end}, are local and hold for
every Reeb type. The global classification of
Sections~\ref{sec:cds-decay}--\ref{sec:quasi} uses a complete smooth quotient.
Example~\ref{ex:counterexample} in Section~\ref{sec:nondev-models} shows that
the
roundness of the limits at infinity cannot be obtained by a purely local
argument: it gives an explicit irregular obstruction.

\subsection{Complete Smooth Quotient and Global Splitting}

\label{sec:cds-decay}

\begin{theorem}
\label{thm:quotient} Suppose that the transverse geometry
of~\eqref{eq:soliton} is realized by a Riemannian submersion onto a connected
complete smooth K\"{a}hler surface $(Z,h_{Z},J_{Z})$, where
$Z=M/\mathbb{S}^{1}$ or $Z=M/\mathbb{R}$, with
\begin{equation*}
\pi :M\longrightarrow Z,\qquad h=\pi ^{\ast }h_{Z}\text{ on }D,\quad f=\pi
^{\ast }f_{Z},\quad d\eta =2\pi ^{\ast }\omega _{Z}.
\end{equation*}
Then
\begin{equation*}
\sup_{M}R^{T}<\infty .
\end{equation*}
In particular this holds for a free proper circle Reeb action and for the
closed proper $\mathbb{R}$-action in Lemma~\ref{lem:reeb}. It also holds
when a complete orbifold transverse quotient admits a complete smooth
manifold orbifold cover to which the metric and potential lift, that is,
when the quotient orbifold is developable as in
Definition~\ref{def:developable}.
\end{theorem}

\begin{proof}
The descended equation is
\begin{equation*}
\mathrm{Ric}_{h_{Z}}+\operatorname{Hess}_{h_{Z}}f_{Z}=6h_{Z}.
\end{equation*}
For $\widehat{h}_{Z}=12h_{Z}$, constant scaling leaves the connection and
the $(0,2)$-Ricci tensor unchanged, and hence
\begin{equation*}
\mathrm{Ric}_{\widehat{h}_{Z}}+\operatorname{Hess}_{\widehat{h}_{Z}}f_{Z}=\tfrac{1}{2}\widehat{h}_{Z}.
\end{equation*}
This is a complete smooth K\"{a}hler--Ricci shrinker surface. The scalar
curvature bound of Li--Wang~\cite[Theorem~1.2]{LW} gives
$\sup_{Z}R_{\widehat{h}_{Z}}<\infty $. Since $R^{T}=12\,\pi ^{\ast
}R_{\widehat{h}_{Z}}$, the assertion follows. In the cover case, apply the
same theorem upstairs
and descend the pointwise bound.
\end{proof}

\begin{remark}
Quasi-regularity gives an orbifold quotient. It does not supply a smooth
manifold orbifold cover. Thus quasi-regularity alone is not a justification
for invoking a theorem stated for smooth K\"{a}hler surfaces. Likewise,
local transverse K\"{a}hler charts do not form a complete global K\"{a}hler
surface. These points decide whether the theorem cited in the proof of
Theorem~\ref{thm:quotient} applies.
\end{remark}

Using the results of Conlon--Deruelle--Sun \cite{CDS}, we obtain the
following.

\begin{theorem}
\label{thm:cds-decay}Let $(M^{5},\xi ,\eta ,\Phi ,g,f)$ be a complete
noncompact regular gradient shrinking Sasaki--Ricci soliton with Reeb period
$L>0$. Assume that
\begin{equation}
R^{T}(p)\longrightarrow 0\qquad \text{uniformly as}\qquad
f(p)\longrightarrow \infty .  \label{eq:cds-global-decay}
\end{equation}
Then the complete smooth quotient $Z=M/\mathbb{S}^{1}$ is, up to biholomorphic
pullback, either the Gaussian K\"{a}hler shrinker on $\mathbb{C}^{2}$ or the
FIK shrinker on
\begin{equation*}
\operatorname{Bl}_{0}\mathbb{C}^{2}=\operatorname{Tot}\bigl(\mathcal{O}(-1)\rightarrow
\mathbb{CP}^{1}\bigr).
\end{equation*}
\end{theorem}

\begin{proof}
Because the fibers are compact, $f_{Z}$ is proper and
\eqref{eq:cds-global-decay} is exactly the assertion $R_{h_{Z}}\rightarrow 0$
at infinity on the complete smooth complex surface $Z$. The full real
equation descends, and its Hessian is $J_{Z}$-invariant. Thus $Z$ is an
ordinary gradient K\"{a}hler shrinker, not merely a solution of a $(1,1)$
equation. Set $h_{\mathrm{C}}=\lambda h_{Z}$ and $f_{\mathrm{C}}=f_{Z}-2$.
This is the normalization of \cite{CDS}, where the shrinking soliton
equation reads
\begin{equation*}
\mathrm{Ric}+\operatorname{Hess}f=g
\end{equation*}
and the potential is normalized by $|\nabla f|^{2}+R-2f=2n$ with $n=\dim
_{\mathbb{C}}=2$ \cite[(1.2) and Remark~2.24]{CDS}. Since
$R_{h_{\mathrm{C}}}=\lambda ^{-1}R_{h_{Z}}$ and $|\nabla
f|_{h_{\mathrm{C}}}^{2}=\lambda
^{-1}|\nabla f|_{h_{Z}}^{2}$, the identity $R_{h_{Z}}+|\nabla
f_{Z}|_{h_{Z}}^{2}=2\lambda f_{Z}$ of \eqref{eq:pre-normal} gives $|\nabla
f_{\mathrm{C}}|_{h_{\mathrm{C}}}^{2}+R_{h_{\mathrm{C}}}-2f_{\mathrm{C}}=4$.
\cite[Theorem~E(3)]{CDS} applies directly and gives the Gaussian/FIK
classification. Completeness, complex dimension two, the full gradient
equation, and global scalar decay verify that theorem's hypotheses.
\end{proof}

For comparison we recall the following result of Munteanu--Wang
\cite{MW2} on the ends of four-dimensional shrinkers.

\begin{proposition}[{\cite{MW2}}]
\label{P51}Let $(N^{4},g,f)$ be a complete four-dimensional
gradient shrinking Ricci soliton with bounded scalar curvature. Suppose that
$R\geq \delta >0$ throughout an end $E$. Then either $E$ is smoothly
asymptotic to a round cylinder $\mathbb{R}\times \mathbb{S}^{3}/\Gamma $, or
every escaping sequence on a fixed gradient trajectory in $E$ has pointed
smooth limit $\mathbb{R}^{2}\times \mathbb{S}^{2}$ or its $\mathbb{Z}/2\mathbb{Z}$
quotient. The limit in the latter alternative is fixed by the trajectory.
The round metrics have the soliton normalization.
\end{proposition}

Proposition~\ref{P51} is recorded for comparison only and is not used in
the proofs below: the statement along a single trajectory needed in
Theorem~\ref{T51} is obtained instead from the fixed-basepoint
tangent-flow argument in its proof.

Furthermore, using the result of Munteanu--Wang, Cifarelli--Conlon--Deruelle
proved the following:

\begin{proposition}[{\cite[Proposition~3.9]{CCD}}]
\label{P55} Let $(N,g,f,J)$ be a simply
connected complete noncompact gradient shrinking K\"{a}hler--Ricci soliton
of complex dimension two with bounded scalar curvature $R$. Assume that
$\nabla f$ has an escaping integral curve along which $R\nrightarrow0$
at infinity. Then $N$ is
biholomorphic to either $\mathbb{CP}^{1}\times \mathbb{C}$ or the blow-up of
$\mathbb{CP}^{1}\times \mathbb{C}$ at one fixed point of the standard torus
action.
\end{proposition}

Combining \cite{MW2}, \cite{CCD} and \cite{BCCD}, we obtain the following.

\begin{theorem}
\label{T51} Let $(M^{5},\xi ,\eta ,\Phi ,g,f)$ be a complete noncompact
regular gradient shrinking Sasaki--Ricci soliton, and suppose that an
escaping integral curve $\gamma $ of $\nabla ^{T}f$ has $R^{T}\nrightarrow 0$.
Then
\begin{equation}
Z\cong _{\mathrm{bihol}}\mathbb{C}\times \mathbb{CP}^{1}\quad \text{or}\quad
\operatorname{Bl}_{p}(\mathbb{C}\times \mathbb{CP}^{1}).  \label{eq:twobases}
\end{equation}
The metric is the product Gaussian-round shrinker in the first case and the
BCCD shrinker in the second, up to biholomorphism and at the fixed
normalization. On the projected trajectory the unique pointed smooth
K\"{a}hler limit at infinity is the round cylinder: for every sequence
$p_{i}\in
\gamma $ tending to infinity, with $q_{i}=\pi (p_{i})$,
\begin{equation}
(Z,h_{Z},J_{Z},q_{i})\longrightarrow (\mathbb{C}\times
\mathbb{CP}^{1},h_{\mathrm{E}}+h_{\lambda },J_{\mathrm{prod}},q_{\infty
}),\qquad K_{h_{\lambda
}}=6,\qquad R^{T}(p_{i})\longrightarrow 12.  \label{eq:quotientlimit}
\end{equation}
\end{theorem}

\begin{proof}
Theorem~\ref{thm:quotient} provides the curvature bound on the complete
smooth quotient, so $Z$ is a complete smooth gradient K\"{a}hler shrinker
with bounded scalar curvature. CCD \cite[Theorem~A]{CCD} gives its
biholomorphic classification \eqref{eq:twobases}. For simply connected $Z,$
this is Proposition~\ref{P55}, and the simple-connectivity assumption is
removed in \cite[Theorem~A(i)]{CCD}. BCCD \cite[Theorem~A]{BCCD} supplies
the existence and uniqueness on the blow-up. The complete bounded-scalar
classification in \cite[Theorem~B]{BCCD} includes the product metric on
$\mathbb{C}\times \mathbb{CP}^{1}$.

For the trajectory statement we argue on the complete smooth quotient $Z$,
whose curvature is bounded by \cite{LW}. After the constant scaling
$\widehat{h}_{Z}=12h_{Z}$ of the proof of Theorem~\ref{thm:quotient}, which
gives the normalization $\lambda =\frac{1}{2}$ of
\cite[Section~\Iref{sec:splitting}{6}]{PartI}, the associated ancient flow
$h(t)=(-t)\Psi
_{t}^{\ast }\widehat{h}_{Z}$ is transversely $(C,\kappa )$-controlled by
\cite[Corollary~\Iref{cor:controlled}{6.1}]{PartI}. All statements below are
invariant under this scaling after rescaling the constants. Let
$x_0=\pi(\gamma(0))$, choosing the origin of the forward trajectory
before all points in the sequence. Write
$\widehat\nabla f_Z=\nabla^{\widehat h_Z}f_Z=12^{-1}\nabla^{h_Z}f_Z$.
The map $\tau\mapsto\Psi_{-\tau}(x_0)$ parametrizes the same forward
gradient trajectory, with the parameter change $s=-\log\tau$ for the
field $\widehat\nabla f_Z$. Thus every escaping sequence $q_i$ on that
trajectory, after omitting finitely many points, is
$q_i=\Psi_{-\tau_i}(x_0)$ with $\tau_i\downarrow0$, and
\begin{equation*}
(Z,\widehat{h}_{Z},q_{i})\cong (Z,\tau _{i}^{-1}h(-\tau _{i}),x_{0}).
\end{equation*}
By \cite[Theorem~\Iref{thm:asymp}{6.1}]{PartI} (Naber's theorem \cite{Nab}
for the quotient flow), a subsequence converges smoothly to the time $-1$
slice $(Z_{\infty },h_{\infty },q_{\infty })$ of a complete normalized
gradient shrinking K\"{a}hler--Ricci soliton, and by
\cite[Theorems~\Iref{P54}{6.6} and~\Iref{thm:splitting}{6.8}]{PartI} its
universal cover
splits as
\begin{equation*}
(\mathbb{C},h_{\mathrm{E}})\times (N,k),
\end{equation*}
where $(N,k)$ with the restricted potential is a complete
 two-dimensional gradient shrinker with normalized constant $1/2$.
By the classification of such surfaces (see \cite{Ham,BM} and
\cite[Theorem~\Iref{thm:splitting}{6.8}]{PartI}) and orientability,
$(N,k)$ is the flat plane or the round sphere of Gaussian curvature
$1/2$. Hence the scalar curvature in the metric $\widehat h_Z$
converges along the subsequence to $0$ or $1$. Returning to the
original metric $h_Z$, for which $\lambda=6$, gives limiting values
$0$ or $12=2\lambda$. From this point onward all metrics and scalar
curvatures are expressed at this original scale. Thus every cluster
value of $R^T$ along $\gamma$ lies in $\{0,2\lambda\}$.
The set of cluster values of a bounded continuous function along a curve
tending to infinity is connected, so either $R^{T}\rightarrow 0$ or
$R^{T}\rightarrow 2\lambda $ along $\gamma $, and our hypothesis excludes the
first alternative. Hence every limit along $\pi (\gamma )$ is a free
quotient of $\mathbb{C}\times \mathbb{S}^{2}$ with the product of the flat
and the round metric. A holomorphic isometry of this product which preserves
the limiting potential $\frac{\lambda }{2}|z-z_{0}|^{2}+1$ fixes $z_{0}$ and
acts on $\mathbb{S}^{2}$ by a rotation. Since every rotation of $\mathbb{S}^{2}$ has a fixed point $y_{0}$, every
such isometry fixes $(z_{0},y_{0})$; hence the only such isometry that can
occur in a free quotient is the identity. Therefore every limit is
\begin{equation*}
(\mathbb{C}\times \mathbb{CP}^{1},h_{\mathrm{E}}+h_{\lambda },J_{\mathrm{prod}})
\end{equation*}
itself, independently of the sequence, which proves \eqref{eq:quotientlimit}.
\end{proof}

\begin{corollary}
\label{cor:four-models-body}Every complete noncompact regular
five-dimensional gradient shrinking Sasaki--Ricci soliton has, up to
biholomorphic pullback, exactly one of the smooth quotient metrics
\begin{equation}
\mathbb{C}^{2}\ \text{(Gaussian)},\qquad \operatorname{Bl}_{0}\mathbb{C}^{2}\
\text{(FIK)},\qquad \mathbb{C}\times \mathbb{CP}^{1}\ \text{(product)},\qquad
\operatorname{Bl}_{p}(\mathbb{C}\times \mathbb{CP}^{1})\ \text{(BCCD)}.
\end{equation}
\end{corollary}

\begin{proof}
By Theorem~\ref{thm:quotient} the quotient $Z$ has bounded scalar curvature.
If $R_{h_{Z}}\rightarrow 0$ along every integral curve of $\nabla f_{Z}$ which
tends to infinity, then $R_{h_{Z}}\rightarrow 0$ at infinity by
\cite[Lemma~2.7]{CCD}, and Theorem~\ref{thm:cds-decay} applies. Otherwise some
escaping integral curve of $\nabla ^{T}f$ has $R^{T}\nrightarrow 0$, and
Theorem~\ref{T51} applies. The four surfaces are pairwise not biholomorphic:
$\mathbb{C}^{2}$ contains no compact complex curve,
$\operatorname{Bl}_{p}(\mathbb{C}\times \mathbb{CP}^{1})$ has second Betti
number two, and the remaining two surfaces, which have second Betti number
one, are distinguished by the self-intersection of a generator of $H_{2}$,
which is $-1$ on $\operatorname{Bl}_{0}\mathbb{C}^{2}$ and $0$ on
$\mathbb{C}\times \mathbb{CP}^{1}$.
\end{proof}

\subsection{Quasi-Regular Structures and the Free Proper-Line Case}

\label{sec:quasi}

We begin with smooth circle orbibundles and their lifting criterion. We then
reduce the developable quasi-regular case to the regular one.

\medskip

We give the local argument needed in the present noncompact setting. Over a
sufficiently small orbifold chart $(U,G)$ centered at a point fixed by $G$,
a circle orbibundle has the form
\begin{equation}
(U\times \mathbb{S}^{1})/G,\qquad \gamma \cdot (z,u)=(\gamma z,\chi (\gamma )u),
\label{eq:localmodel}
\end{equation}
where $\chi :G\rightarrow \mathbb{S}^{1}$ is a character. Its restriction to
each
stabilizer is part of the bundle data.

\begin{lemma}
\label{lem:faithful} A circle orbibundle has a smooth total space with no
residual orbifold isotropy if and only if each local stabilizer acts
faithfully on the circle fiber. In particular, every local group must be
cyclic. If $G=\mathbb{Z}/m\mathbb{Z}$ and a generator acts on the fiber by
$e^{2\pi ib/m}$, the condition is
\begin{equation}
\gcd (b,m)=1.  \label{eq:coprime}
\end{equation}
\end{lemma}

\begin{proof}
An element fixing $(z,u)$ in \eqref{eq:localmodel} must fix $z$ and satisfy
$\chi (\gamma )=1$. Thus a faithful character on every stabilizer makes the
diagonal action free. Conversely, at the center of a chart an element of $\ker
\chi $ fixes the entire circle above that center. A finite subgroup of
$\mathbb{S}^{1}$ is cyclic, and the character of a cyclic group is faithful
precisely
under \eqref{eq:coprime}.
\end{proof}

\begin{lemma}
\label{lem:finite} Let $\mathcal{Z}$ be a complete noncompact gradient
shrinking K\"{a}hler orbifold whose potential is globally defined. If
$\mathcal{Z}$ is developable, its universal orbifold cover is a finite cover.
\end{lemma}

\begin{proof}
Lift the metric, complex structure, and potential to the smooth universal
orbifold cover $\widetilde{Z}$. An orbifold Riemannian covering of a
complete orbifold is complete, so this is an ordinary complete gradient
K\"{a}hler shrinker. Its potential $\widetilde{f}$ is proper by the standard
quadratic potential estimate, after a constant normalization and a constant
metric scaling; see \cite[Theorem~2.21(i)]{CDS} or
\cite[Theorem~1.1]{CaZ}. Each deck transformation preserves
$\widetilde{f}$ exactly, since it is a pullback. Its orbit through any point
lies in a compact level set. Properness of the deck action makes the set of
group elements carrying that point into this compact set finite. Hence $\Gamma
$ is finite. Properness is used only on the smooth cover $\widetilde{Z}$ and
developability of $\mathcal{Z}$ is a hypothesis of the lemma. In particular,
this construction does not use the global noncollapsing or splitting
results of \cite{PartI}.
\end{proof}

\begin{theorem}
\label{thm:good} Let $(M^{5},\xi ,\eta ,\Phi ,g,f)$ be a connected complete
noncompact quasi-regular gradient shrinking Sasaki--Ricci soliton with
constant $\lambda >0$ and generic Reeb period $L$, and suppose that its leaf
space $\mathcal{Z}$ is developable. There is a finite ordinary covering
$\widetilde{M}\rightarrow M$ preserving the pulled-back Sasakian structure
and potential, such that the Reeb flow on $\widetilde{M}$ is regular with
period $L$. Its smooth K\"{a}hler quotient $\widetilde{Z}$ is one of the
four smooth models of Corollary~\ref{cor:four-models}. More precisely,

\begin{enumerate}
\item If $R^{T}\rightarrow 0$ uniformly at transverse infinity,
$\widetilde{Z}$ is Gaussian or FIK, and the regular bundle formulas of
Theorem~\ref{thm:decay-lifts} apply to $\widetilde{M}$.

\item If an escaping gradient trajectory has $R^{T}\nrightarrow 0$,
$\widetilde{Z}$ is the round K\"{a}hler cylinder or the BCCD surface, and
Theorems~\ref{T51} and \ref{thm:regular} apply to $\widetilde{M}$ along a
lift of that trajectory.
\end{enumerate}

The original manifold is a finite free Sasakian quotient
\begin{equation}
M=\widetilde{M}/\Gamma ,  \label{eq:finitequotient}
\end{equation}
where $\Gamma =\pi _{1}^{\mathrm{orb}}(\mathcal{Z})$ acts on $\widetilde{M}$
preserving $\eta ,\xi ,\Phi ,g,\widetilde{f}$ and commuting with the circle
action.
\end{theorem}

\begin{proof}
Lemma~\ref{lem:finite} gives the finite smooth K\"{a}hler cover
$\widetilde{Z}$. Pull back $M\rightarrow \mathcal{Z}$ in the \emph{orbifold}
sense. The
result is an ordinary principal circle bundle $\widetilde{M}\rightarrow
\widetilde{Z}$, not the naive fiber product of coarse spaces. In the local
model \eqref{eq:localmodel}, its map to $M$ is the map from $U\times
\mathbb{S}^{1}$
to the free finite quotient. It is therefore an ordinary
local-diffeomorphism covering of smooth total spaces. The pullback is
connected because both base and fiber are connected. Its circle generator
has the same generic normalization \eqref{eq:normalization}, so its regular
Reeb period is $L$.

The lifted base satisfies the complete smooth K\"{a}hler gradient equation;
Li--Wang \cite{LW} bound its curvature without any a priori assumption. The
smooth classification used in Corollary~\ref{cor:four-models}, with its
separation of the sources CDS, CCD and BCCD, now applies on $\widetilde{Z}$.
The full bounded-scalar list is also recorded in \cite[Theorem~B]{BCCD}.
Because the covers are finite, escape and global decay are preserved and
gradient trajectories lift. Applying the regular lifting calculation on the
smooth cover proves the two refinements. Deck transformations of
$\widetilde{M}$ preserve the pulled-back tensors and act freely; this gives
\eqref{eq:finitequotient}.
\end{proof}

To relate this construction to setting (S) of \cite{PartI}, write
$q:\widetilde M\to M$ and $\pi:\widetilde M\to\widetilde Z$ for the
covering and quotient maps. After a constant normalization to
$\lambda=\tfrac12$, the canonical quotient flow is
\[
\widetilde h_Z(t)=(-t)\Psi_t^*\widetilde h_Z,\qquad t<0,
\]
where $\Psi_{-1}=\mathrm{id}$ and $\Psi_t$ is generated by
$(-t)^{-1}\nabla\widetilde f_Z$. The linear growth bound for the
soliton gradient makes these biholomorphisms globally defined; hence
every time slice is complete. The lifted soliton vector field is
$q$-related to the original field, projects under $\pi$ to
$\nabla\widetilde f_Z$, commutes with the Reeb flow, and is invariant
under the deck action. Thus the associated transverse flow is realized
on this same smooth quotient for all $t<0$, with $q$, $\pi$, and the
complex structure fixed. This verifies setting (S) for the canonical
flow; the quotient volumes and noncollapsing constants refer to
$\widetilde Z$.

\subsection{Integral Liftings and Topology}

We will construct the regular or quasi-regular gradient Sasaki--Ricci
shrinkers $M$ from the quotient space $Z$ as follows.

\noindent (I) \emph{Regular lifts.}

\begin{proposition}
\label{prop:lifting} Let $(Z,h_{Z},J_{Z},f_{Z})$ be a complete K\"{a}hler
surface satisfying
\begin{equation*}
\mathrm{Ric}_{Z}+\operatorname{Hess}_{Z}f_{Z}=\lambda h_{Z}.
\end{equation*}
A regular Sasakian circle lift with Reeb period $L>0$ and $d\eta =2\pi
^{\ast }\omega _{Z}$ exists if and only if
\begin{equation}
e_{\mathbb{R}}=\frac{2}{L}[\omega _{Z}]=\frac{4\pi }{\lambda
L}\,c_{1}(Z)_{\mathbb{R}}  \label{eq:euler}
\end{equation}
is the real image of an integral Euler class $e\in H^{2}(Z;\mathbb{Z})$.
After such a class is specified, the lift has
\begin{equation}
g=\pi ^{\ast }h_{Z}+\eta ^{2},\qquad f=\pi ^{\ast }f_{Z}
\label{eq:liftmetric}
\end{equation}
and is complete and satisfies the soliton equation. On a simply
connected base with fixed principal circle bundle, curvature and period,
the connection is unique up to gauge, meaning a circle-equivariant
bundle automorphism covering the identity on $Z$.
\end{proposition}

\begin{proof}
Normalize the circle action to have parameter period $2\pi $. Its
fundamental vector is $(L/2\pi )\xi $, and its connection is $\theta =(2\pi
/L)\eta $. Write $d\theta=\pi^*F_\theta$, where
$F_\theta=(4\pi/L)\omega_Z$ is the curvature form on the base. Then
\begin{equation*}
e_{\mathbb R}=[F_\theta/(2\pi)]=\frac{2}{L}[\omega_Z].
\end{equation*}
With $d^{c}f=-df\circ J$, the full real soliton equation implies $\rho
_{Z}+\tfrac{1}{2}dd^{c}f_{Z}=\lambda \omega _{Z}$. Since $[\rho _{Z}]=2\pi
c_{1}(Z)_{\mathbb{R}}$, this proves the second equality in \eqref{eq:euler}.

Conversely take the principal circle bundle of an integral class with that
real image. A connection with curvature representing it can be adjusted by a
basic one-form to have curvature $(4\pi /L)\omega _{Z}$. Multiplying the
connection by $L/(2\pi )$ gives $\eta $ with the required period and
curvature. The standard K\"{a}hler connection construction~\cite{BG}, which
is globally the local contactization of
\cite[Proposition~\Iref{prop:contactization}{2.1}]{PartI}, gives a Sasakian
structure. Its
transverse metric and Hessian are pulled back from the base, so the soliton
equation follows. The bundle projection is proper because its fibers
are compact and the bundle is locally trivial. A closed metric ball
in the total space projects into a compact ball of the complete base
and is therefore contained in a compact inverse image. It is compact,
which proves completeness by Hopf--Rinow. Finally two such
connections differ by a closed basic one-form; on a simply connected base it
is exact and is removed by a gauge transformation. Circle bundles over the
smooth base are classified by their integral Euler class \cite{Ste,Hat}.
\end{proof}

\begin{theorem}
\label{thm:regular}In Theorem~\ref{T51} the following statements hold.

\begin{enumerate}
\item For $Z=\mathbb{C}\times \mathbb{CP}^{1}$,
$d:=e(\mathbb{CP}^{1})=\frac{8\pi }{\lambda L}\in \mathbb{Z}_{>0}$, and
\begin{equation}
\qquad M\cong \mathbb{C}\times L(d,1)\qquad  \label{eq:cylinder-topology}
\end{equation}
with $\pi _{1}(M)=H^{2}(M;\mathbb{Z})=\mathbb{Z}/d\mathbb{Z}.$ Here $L(d,1)$
is the lens space defined before Theorem~\ref{thm:intro-topology}, and
every positive integer $d$ occurs.

\item For $Z=\operatorname{Bl}_{p}(\mathbb{C}\times \mathbb{CP}^{1})$, let $F$
be a
generic sphere fiber and $E$ the exceptional curve. Then, for $q:=\frac{4\pi
}{\lambda L}\in \mathbb{Z}_{>0}$,
\begin{equation}
e(F)=2q,\quad e(E)=q,  \label{eq:BCCD-euler}
\end{equation}
and
\begin{equation}
\pi ^{-1}(E)\cong L(q,1),  \label{eq:BCCD-topology}
\end{equation}
with
\begin{equation*}
\pi _{1}(M)=\mathbb{Z}/q\mathbb{Z},\qquad H^{2}(M;\mathbb{Z})=\mathbb{Z}\oplus
\mathbb{Z}/q\mathbb{Z}.
\end{equation*}
Every positive integer $q$ is realized by a regular lift of the BCCD
shrinker.
\end{enumerate}
\end{theorem}

\begin{proof}
For the cylinder, $c_{1}(Z)(\mathbb{CP}^{1})=2$, so \eqref{eq:euler} gives
\eqref{eq:cylinder-topology}. The base retracts onto its sphere factor.
The circle bundle is therefore pulled back from
$p_d:L(d,1)\to\mathbb{CP}^{1}$, whose Euler class is $d u$ in our
curvature convention. Choosing $L=8\pi /(\lambda
d)$ realizes any $d>0$, and $\pi
_{1}(M)=H^{2}(M;\mathbb{Z})=\mathbb{Z}/d\mathbb{Z}$ are the invariants of
$L(d,1)$.

For the blow-up, $H_{2}(Z;\mathbb{Z})=\mathbb{Z}[F]\oplus \mathbb{Z}[E]$, and
\begin{equation*}
c_{1}(Z)(F)=2,\qquad c_{1}(Z)(E)=1.
\end{equation*}
The second identity follows from $E^{2}=-1$ and adjunction. Thus
\eqref{eq:euler} forces \eqref{eq:BCCD-euler}; conversely $q\,c_{1}(Z)$ is
an integral class and gives a lift for every $q>0$. The base is simply
connected. The homotopy sequence of its circle bundle gives
\begin{equation*}
\pi _{2}(Z)\xrightarrow{\partial}\mathbb{Z}\longrightarrow \pi
_{1}(M)\longrightarrow 0,\qquad \operatorname{im}\partial=q\mathbb{Z}.
\end{equation*}
This yields $\pi _{1}(M)=\mathbb{Z}/q\mathbb{Z}$. The Gysin
sequence~\cite{Hat} gives
\begin{equation*}
H^{2}(M;\mathbb{Z})=H^{2}(Z;\mathbb{Z})/\langle e\rangle
=\mathbb{Z}^{2}/\langle (2q,q)\rangle =\mathbb{Z}\oplus
\mathbb{Z}/q\mathbb{Z}.
\end{equation*}
Restriction of the bundle to $E$ gives the indicated lens space. In
particular $b_{2}(M;\mathbb{Q})=1$, so these total spaces cannot be
$\mathbb{C}\times \mathbb{S}^{3}$ or its lens-space finite quotient, both of
which
have zero rational second cohomology.
\end{proof}

\begin{remark}
\label{rem:cylmodel} The total space is not determined by the base's
biholomorphic type alone. The Reeb period $L$ and integral Euler class $d$
are essential. Cylindrical lifts occur on
\begin{equation*}
\mathbb{C}\times L(d,1)
\end{equation*}
for every positive integer $d$. BCCD lifts have different cohomology. For
example, let $h_{\lambda }$ be round of Gaussian curvature $\lambda $, and
let $\eta _{d}$ be a connection on $L(d,1)\rightarrow \mathbb{CP}^{1}$ with
curvature $2\omega _{\lambda }$ and Reeb period $L=8\pi /(\lambda d)$. On
$\mathbb{C}\times L(d,1),$ set $\eta =\eta _{d}+x\,dy-y\,dx$ with
\begin{equation}
g=dx^{2}+dy^{2}+\pi ^{\ast }h_{\lambda }+\eta ^{2},\qquad f=1+\frac{\lambda
}{2}(x^{2}+y^{2}).  \label{eq:cylmodel}
\end{equation}
The transverse Ricci tensor is zero on $\mathbb{C}$ and $\lambda h_{\lambda
} $ on the sphere, while the transverse Hessian is $\lambda h_{\mathrm{E}}$
on $\mathbb{C}$ and zero on the sphere. Thus the full soliton equation and
the Hamilton identity \eqref{eq:pre-normal} hold, and
\begin{equation}
R^{T}=2\lambda ,\qquad R_{g}=2\lambda -4.  \label{eq:cylscalar}
\end{equation}
For the normalization $\lambda =6$, one has $R_{g}\equiv 8$. These metrics
are complete, their potentials are proper, and every $d>0$ is allowed. They
are connection metrics, not ordinary product Riemannian metrics.
\end{remark}

\medskip

\noindent We now justify the last assertion of Theorem~\ref{thm:intro-cyl}.
Let $M$ be regular with Reeb period $L$ and quotient $Z$ as in
Theorem~\ref{T51}, at an arbitrary scale $\lambda $, and let $F_{i}$ be
embeddings,
defined on an exhaustion of $\mathbb{C}\times \mathbb{CP}^{1}$, which
realize the convergence \eqref{eq:quotientlimit} along $\pi (\gamma )$. Fix
$R>0$ and put $B_{R}=\{|z|<R\}\times \mathbb{CP}^{1}$, which retracts onto
$\{0\}\times \mathbb{CP}^{1}$. For large $i$ the pulled-back circle bundle
$F_{i}^{\ast }M\rightarrow B_{R}$ has the integral Euler number
\begin{equation*}
\frac{2}{L}\int_{\{0\}\times \mathbb{CP}^{1}}F_{i}^{\ast }\omega
_{Z}\longrightarrow \frac{2}{L}\cdot \frac{4\pi }{\lambda }=\frac{8\pi
}{\lambda L}
\end{equation*}
by Proposition~\ref{prop:lifting}, since the round sphere of Gaussian
curvature $\lambda $ has area $4\pi /\lambda $. Hence it equals $d=8\pi
/(\lambda L)$ for large $i$, and $F_{i}^{\ast }M$ is isomorphic over $B_{R}$
to the restriction of the cylindrical lift $\mathbb{C}\times L(d,1)$ of
\eqref{eq:cylmodel}. Under such isomorphisms the pulled-back contact forms
and the form $\eta $ of \eqref{eq:cylmodel} have the same period
normalization. They differ by $\pi ^{\ast }\alpha _{i}$ for one-forms $\alpha
_{i}$ on $B_{R}$ with $d\alpha _{i}=2(F_{i}^{\ast }\omega _{Z}-\omega
_{\infty })\rightarrow 0$ smoothly, where $\omega _{\infty }$ is the
K\"{a}hler form of the limit. Since $B_{R}$ is simply connected, a gauge
transformation changes $\alpha_i$ by an exact form. We impose the
Coulomb condition $d^*\alpha_i=0$ with the boundary condition
$\alpha_i(\nu)=0$ on $\partial B_R$, using the fixed limit metric and
its outward unit normal $\nu$. The resulting harmonic ambiguity is
zero because $H^1(B_R;\mathbb R)=0$. The estimates for this Hodge
boundary problem and $d\alpha_i\to0$ give
$\alpha_i\to0$ smoothly on $B_{R/2}$. As the Sasakian structure is
determined by $(h_{Z},J_{Z})$ and the connection through \eqref{eq:pre-BW},
the total spaces converge, after these gauge transformations, smoothly on $\pi
^{-1}(B_{R/2})$ to the cylindrical lift \eqref{eq:cylmodel} with this $d$. In
the BCCD case $d=2q=e(F)$, in accordance with \eqref{eq:BCCD-euler}.

\begin{theorem}
\label{thm:decay-lifts}In Theorem~\ref{thm:cds-decay}, fix the original
scale $\lambda $ and the Reeb period $L$.

\begin{enumerate}
\item If the quotient is Gaussian, every $L>0$ is allowed and the regular
lift has underlying manifold
\begin{equation*}
\mathbb{C}^{2}\times \mathbb{S}^{1}.
\end{equation*}
Its integral Euler class vanishes, its fundamental group is $\mathbb{Z}$,
and $H^{2}(M;\mathbb{Z})=0$.

\item If the quotient is the FIK surface, let $E\cong \mathbb{CP}^{1}$ be
its exceptional curve. Then, for $q:=\frac{4\pi }{\lambda L}\in
\mathbb{Z}_{>0}$,
\begin{equation}
\qquad e(E)=q,\qquad \pi _{1}(M)=\mathbb{Z}/q\mathbb{Z},\qquad
H^{2}(M;\mathbb{Z})=\mathbb{Z}/q\mathbb{Z}.  \label{eq:fik-lift-invariants}
\end{equation}
Every positive integer $q$ occurs. Let
$p_q:L(q,1)\to\mathbb{CP}^{1}$ be the principal circle bundle with
$e(p_q)=q u$ defined before Theorem~\ref{thm:intro-topology}. Then the
total space is diffeomorphic to
\begin{equation}
M\cong \operatorname{Tot}\bigl(p_{q}^{\ast }\mathcal{O}(-1)\longrightarrow
L(q,1)\bigr),\qquad \pi ^{-1}(E)=L(q,1).  \label{eq:fik-totalspace}
\end{equation}
For $q=1$, the pulled-back complex line bundle is trivial and
\begin{equation*}
M\cong \mathbb{C}\times \mathbb{S}^{3}.
\end{equation*}
For $q>1,$ its first Chern class is $-1$ modulo $q$ and the line bundle is
nontrivial.
\end{enumerate}
\end{theorem}

\begin{proof}
These are bundle calculations using Proposition~\ref{prop:lifting}, not
additional claims of the CDS classification theorem. On $\mathbb{C}^{2}$,
$H^{2}(\mathbb{C}^{2};\mathbb{Z})=0$ and $\omega _{Z}$ is exact. The lifting
criterion therefore gives the trivial principal circle bundle for every $L$.
Its connection metric is the central-period quotient of the Heisenberg model
\begin{equation}
\eta _{0}=dt+\sum_{a=1}^{2}(x_{a}\,dy_{a}-y_{a}\,dx_{a}),\quad
g_{0}=|dz|^{2}+\eta _{0}^{2},\quad f_{0}=\frac{\lambda }{2}|z|^{2}.
\label{eq:heisenberg}
\end{equation}
This complete transverse Gaussian has
$\mathrm{Rm}^{T}=\mathrm{Ric}^{T}=R^{T}=0$ and $R_{g}=-4$.

On $Z=\operatorname{Tot}(\mathcal{O}(-1)\rightarrow \mathbb{CP}^{1})$, the zero
section is a deformation retract, so $H^{2}(Z;\mathbb{Z})=\mathbb{Z}$ and $\pi
_{1}(Z)=0$. Adjunction gives $c_{1}(Z)\cdot E=1$. Formula~\eqref{eq:euler} now
forces \eqref{eq:fik-lift-invariants}'s integrality
condition. Conversely, $qc_{1}(Z)$ is integral for every positive integer $q$,
so Proposition~\ref{prop:lifting} constructs a complete regular FIK lift
with $L=4\pi /(\lambda q)$.

For clarity, the sign and normalization of $p_q$ can be checked directly.
Let $H:\mathbb S^3\to\mathbb{CP}^{1}$ be the Hopf projection and let
\[
\alpha=\left.\sum_{j=0}^{1}(x_j\,dy_j-y_j\,dx_j)\right|_{\mathbb S^3},
\qquad z_j=x_j+iy_j.
\]
Then $\alpha(iz)=1$ and $d\alpha=2H^*\omega_{\mathrm{FS}}$, with
$\int_{\mathbb{CP}^{1}}\omega_{\mathrm{FS}}=\pi$. If
$r_q:\mathbb S^3\to L(q,1)$ is the quotient map, the connection
$\theta_q$ for the effective circle action on $L(q,1)$ is determined by
$r_q^*\theta_q=q\alpha$. Hence
\[
d\theta_q=p_q^*(2q\omega_{\mathrm{FS}}),\qquad
\frac{1}{2\pi}\int_{\mathbb{CP}^{1}}2q\omega_{\mathrm{FS}}=q,
\]
which proves $e(p_q)=q u$. This is the positive-curvature convention
for $e$ used throughout this paper. The same lens space is the unit-circle
bundle of $\mathcal O(-q)$, whose usual complex first Chern class is
$-q u$; these two sign conventions should therefore be distinguished.

Let $\rho:Z\to\mathbb{CP}^{1}$ be the projection of
$Z=\operatorname{Tot}(\mathcal O(-1))$. Since $c_1(Z)=\rho^*u$,
classification by the Euler class identifies $M\to Z$ with the pullback
circle bundle $\rho^*(p_q)$. Its total space is the fiber product
\[
\{(v,\ell)\in Z\times L(q,1):\rho(v)=p_q(\ell)\}.
\]
Projecting to $\ell$ instead of $v$ identifies this same space with the
complex line bundle $p_q^*\mathcal O(-1)\to L(q,1)$: its fiber at
$\ell$ is the tautological line over $p_q(\ell)$. This proves
\eqref{eq:fik-totalspace}. The fibration homotopy sequence has connecting
homomorphism
\begin{equation*}
\pi _{2}(Z)=\mathbb{Z}\rightarrow \pi _{1}(\mathbb{S}^{1})=\mathbb{Z}
\end{equation*}
with image $q\mathbb{Z}$, proving the fundamental group formula. The
Gysin sequence gives
\begin{equation*}
H^{2}(M;\mathbb{Z})=H^{2}(Z;\mathbb{Z})/\langle e\rangle
=\mathbb{Z}/q\mathbb{Z}.
\end{equation*}
Likewise $H^{2}(L(q,1);\mathbb{Z})=\mathbb{Z}/q\mathbb{Z}$, and $p_{q}^{\ast
}:H^{2}(\mathbb{CP}^{1};\mathbb{Z})\rightarrow H^{2}(L(q,1);\mathbb{Z})$ is
reduction modulo $q$. Therefore $c_{1}(p_{q}^{\ast }\mathcal{O}(-1))=-1\bmod
q$. For $q=1,$ the base is $\mathbb{S}^{3}$ with zero second cohomology, so
the complex line bundle is trivial. All these uses of the Gysin and homotopy
sequences are integral.
\end{proof}

\bigskip

\noindent (II) \emph{Orbifold lifts.}

Here $H^2_{\mathrm{orb}}(\mathcal Z;\mathbb Z)$ denotes the ordinary
integral cohomology of the orbifold classifying space, which classifies
principal circle orbibundles. It is not the cohomology of the coarse
space or Chen--Ruan cohomology. The subscript $\mathbb R$ denotes its
image in orbifold de Rham cohomology. For each local stabilizer $G_x$,
the restriction of the bundle gives a character
$\chi_x:G_x\to\mathbb S^1$ on the fiber, as in
\eqref{eq:localmodel}; this integral isotropy data is part of the
bundle and is not determined by its real curvature class alone.

Normalize the effective circle parameter to have period $2\pi $. Its
generator and connection are
\begin{equation}
U=\frac{L}{2\pi }\xi ,\qquad \theta =\frac{2\pi }{L}\eta ,\qquad \theta
(U)=1.  \label{eq:normalization}
\end{equation}
If $d\theta =\pi ^{\ast }F_{\theta }$, the Euler convention in this paper is
\begin{equation*}
e_{\mathbb{R}}=[F_{\theta }/(2\pi )].
\end{equation*}
This is the same curvature convention used to define $e(p_q)=q u$
before Theorem~\ref{thm:intro-topology}.

\begin{theorem}
\label{thm:orb-lift} Let $(\mathcal{Z},h_{\mathcal{Z}},J,f_{\mathcal{Z}})$
be a connected complete effective K\"{a}hler orbifold satisfying
\begin{equation*}
\mathrm{Ric}_{\mathcal{Z}}+\operatorname{Hess}_{\mathcal{Z}}f_{\mathcal{Z}}=\lambda
h_{\mathcal{Z}},\qquad \lambda >0.
\end{equation*}
Fix $L>0$. A smooth quasi-regular Sasakian circle lift with generic Reeb
period $L$, metric and curvature as in \eqref{eq:pre-BW}, exists if and only
if there is a class $e\in H_{\mathrm{orb}}^{2}(\mathcal{Z};\mathbb{Z})$ such
that
\begin{equation}
e_{\mathbb{R}}=\frac{2}{L}[\omega _{\mathcal{Z}}]=\frac{4\pi }{\lambda
L}\,c_{1}^{\mathrm{orb}}(T\mathcal{Z})_{\mathbb{R}}  \label{eq:euler-orb}
\end{equation}
and its local characters are faithful. Every such lift is complete and
satisfies the soliton equation with constant $\lambda $ and $f=\pi ^{\ast
}f_{\mathcal{Z}}$.
\end{theorem}

\begin{proof}
The curvature equality and \eqref{eq:normalization} give $F_{\theta }=(4\pi
/L)\omega _{\mathcal{Z}}$ and hence the first equality of
\eqref{eq:euler-orb}. With $d^{c}f=-df\circ J$, the real soliton equation
implies
\begin{equation*}
\rho _{\mathcal{Z}}+\tfrac{1}{2}dd^{c}f_{\mathcal{Z}}=\lambda \omega
_{\mathcal{Z}}.
\end{equation*}
Since $[\rho _{\mathcal{Z}}]=2\pi
c_{1}^{\mathrm{orb}}(T\mathcal{Z})_{\mathbb{R}}$, it gives the second
equality. Smoothness forces the local
condition of Lemma~\ref{lem:faithful}.

Conversely choose the circle orbibundle represented by $e$ and an invariant
connection $\theta _{0}$. Its curvature and $(4\pi /L)\omega _{\mathcal{Z}}$
have the same orbifold de Rham class, so their difference is $da$ for a
smooth orbifold one-form $a$. Replacing $\theta _{0}$ by $\theta _{0}+\pi
^{\ast }a$ prescribes exactly that curvature. Set $\eta =(L/2\pi )\theta $
and $g=\pi ^{\ast }h_{\mathcal{Z}}+\eta ^{2}$. The local diagonal actions
are free by faithfulness. Thus these tensors define a smooth positive
definite metric and a contact form on the total space.
As in the proof of Proposition~\ref{prop:lifting}, the transverse structure
is pulled back from $\mathcal{Z}$, so the local contactization of
\cite[Proposition~\Iref{prop:contactization}{2.1}]{PartI} shows that the
structure
is Sasakian, with $\Phi $ the horizontal lift of $J$, and the soliton
equation with constant $\lambda $ holds for $f=\pi ^{\ast }f_{\mathcal{Z}}$.
The generic Reeb period is $L$ by \eqref{eq:normalization}. Completeness
follows as in Proposition~\ref{prop:lifting}: the compact-fiber
orbibundle projection is proper (as seen in its finite local charts),
and a closed metric ball lies in the compact inverse image of a
closed ball of the complete orbifold base.
\end{proof}

\begin{remark}
\label{rem:weighted-example} Regularity in Theorem~\ref{T51} and
developability in Theorem~\ref{thm:good} cannot be dropped. Let $m_{1},m_{2}$
be coprime positive integers with $(m_{1},m_{2})\neq (1,1)$, and let
$\mathbb{CP}_{[m_{1},m_{2}]}^{1}$ be the weighted projective line, the
quotient of $\mathbb{S}^{3}\subset \mathbb{C}^{2}$ by the circle action
\begin{equation*}
e^{i\theta }\cdot (w_{1},w_{2})=(e^{im_{1}\theta }w_{1},e^{im_{2}\theta
}w_{2}).
\end{equation*}
It is a bad orbifold, so it carries no metric of constant curvature. By Wu
\cite{Wu} it carries a gradient shrinking Ricci soliton, which after scaling
and normalization satisfies
\begin{equation*}
\mathrm{Ric}_{h_{m}}+\operatorname{Hess}_{h_{m}}f_{m}=\lambda h_{m},\qquad
R_{h_{m}}+|\nabla f_{m}|^{2}=2\lambda f_{m};
\end{equation*}
in real dimension two this is the full real K\"{a}hler equation, and $f_{m}$
is not constant. Put
\begin{equation*}
\mathcal{Z}=\mathbb{C}\times \mathbb{CP}_{[m_{1},m_{2}]}^{1},\qquad
h_{\mathcal{Z}}=|dz|^{2}+h_{m},\qquad f_{\mathcal{Z}}=\tfrac{\lambda
}{2}|z|^{2}+f_{m}.
\end{equation*}
This is a complete effective K\"{a}hler orbifold shrinker. Let $e$ be the
pullback of the orbifold Euler class of the weighted Hopf fibration
\begin{equation*}
\mathbb{S}^{3}\rightarrow \mathbb{CP}_{[m_{1},m_{2}]}^{1}.
\end{equation*}
Its local characters are faithful because $\mathbb{S}^{3}$ is smooth.
Orbifold integration over
$\mathbb{CP}_{[m_1,m_2]}^1$ gives $1/(m_1m_2)$ for $e$ and
$1/m_1+1/m_2$ for $c_1^{\mathrm{orb}}$. These pairings use the
orbifold fundamental class, so they can be rational even though $e$
is an integral orbifold class. Thus
$c_{1}^{\mathrm{orb}}(T\mathcal{Z})_{\mathbb{R}}=(m_{1}+m_{2})e_{\mathbb{R}}$,
and \eqref{eq:euler-orb} holds with $L=4\pi (m_{1}+m_{2})/\lambda $.
Theorem~\ref{thm:orb-lift} therefore gives a complete gradient shrinking
Sasaki--Ricci
soliton on
\begin{equation*}
M=\mathbb{C}\times \mathbb{S}^{3}
\end{equation*}
with $R^{T}=\pi ^{\ast }R_{h_{m}}$ bounded. It is quasi-regular but not
regular. The orbits over the two orbifold points have minimal periods
$L/m_{1}$ and $L/m_{2}$, and one of $m_{1},m_{2}$ is at least $2$. Its leaf
space is not developable, and it is not one of the four models of
Corollary~\ref{cor:four-models}.
\end{remark}

\begin{corollary}
\label{cor:Rcase} Under \eqref{eq:soliton}, if the Reeb action is closed,
free and proper with group $\mathbb{R}$, then its smooth quotient is the
Gaussian shrinker on $\mathbb{C}^{2}$ and the total Sasakian metric is the
Heisenberg model \eqref{eq:heisenberg}.
\end{corollary}

\begin{proof}
Since the Reeb action is free and proper, $\pi :M\rightarrow Z=M/\mathbb{R}$
is a principal $\mathbb{R}$-bundle. The quotient is a complete smooth
K\"{a}hler manifold and its Riemannian metric is the submersion metric, and
horizontal lifts of finite-length quotient paths exist on the complete total
space. The basic function $f$ descends and satisfies
\begin{equation*}
\mathrm{Ric}_{h_{Z}}+\operatorname{Hess}_{h_{Z}}f_{Z}=\lambda h_{Z}.
\end{equation*}
The structure group $\mathbb{R}$ is contractible, so there is a smooth
section $s:Z\rightarrow M$. The Sasakian curvature identity gives
\begin{equation*}
2\omega _{Z}=d(s^{\ast }\eta ).
\end{equation*}
Consequently $\omega _{Z}$ is exact.

\textbf{$Z$ is the Gaussian $\mathbb{C}^{2}$.} No compact complex curve
$C\subset Z$ can exist by Stokes' theorem. In particular $Z$ is noncompact.
By Theorem~\ref{thm:quotient} the quotient has bounded scalar curvature. The
complete smooth K\"{a}hler shrinking-surface classification cited in
\cite{CDS,CCD,BCCD,LW} then restricts $Z$ to the Gaussian, FIK,
product-cylinder,
or BCCD model. Each of the latter three contains a compact holomorphic
rational curve (the exceptional divisor or a $\mathbb{CP}^{1}$ factor), so $Z
$ is the Gaussian $\mathbb{C}^{2}$.

\textbf{Identification with the Heisenberg model.} Choose the global
$\mathbb{R}$-fiber coordinate $t$ determined by the section. Write
\begin{equation*}
\eta =dt+\pi ^{\ast }\alpha
\end{equation*}
and then $d\alpha =2\omega _{\mathrm{E}}$. Since
$H^{1}(\mathbb{C}^{2};\mathbb{R})=0$, $\alpha -\alpha _{0}=du$ for the
standard Heisenberg
primitive $\alpha _{0}=\sum_{j}(x_{j}dy_{j}-y_{j}dx_{j})$. Replacing $t$ by
$t+u(z)$ identifies the contact form, horizontal complex structure and metric
with \eqref{eq:heisenberg}. Finally, the quotient soliton equation implies
$\operatorname{Hess}f_{Z}=\lambda h_{\mathrm{E}}$, so
\begin{equation*}
f_{Z}(z)=\frac{\lambda }{2}|z-z_{\ast }|^{2}+b.
\end{equation*}
The Hamilton identity \eqref{eq:pre-normal}, which here reads $|\nabla
f_{Z}|^{2}=2\lambda f_{Z}$, forces $b=0$. A Heisenberg translation moves
$z_{\ast }$ to the origin. Thus the normalized potential is exactly
\begin{equation*}
f_{0}=\frac{\lambda }{2}|z|^{2}.
\end{equation*}
\end{proof}

\section{Structure Theorems for Irregular Sasaki--Ricci Shrinkers}

\label{sec:irregular}

For an arbitrary Reeb type, we have the following local structure.

\begin{theorem}
\label{thm:irregular-local} Let $M^{5}$ be a complete gradient shrinking
Sasaki--Ricci soliton with bounded transverse scalar curvature, of arbitrary
Reeb type. Let $\gamma $ be an escaping gradient trajectory with
$R^{T}\nrightarrow 0$. Then there are $p_{i}\in \gamma $ and $\delta >0$ such
that $R^{T}(p_{i})\geq \delta $, and a subsequence has a developed local
transverse limit
\begin{equation*}
(U,h_{\mathrm{E}})\times (V,h_{V}),\qquad K_{V}\geq 0,\qquad K_{V}(q_{\infty
})\geq \delta /2.
\end{equation*}
Here $U$ is an open subset of the flat complex line, $V$ is a transverse
surface chart, and $q_\infty\in V$ is the surface component of the limiting
basepoint. If the local soliton equation holds for this limit, the surface is a
positive-curvature local shrinking soliton \eqref{LS} satisfying
\begin{equation*}
\operatorname{Hess}\psi =(\lambda -K_{V})h_{V},\qquad dK_{V}=K_{V}\,d\psi .
\end{equation*}
Roundness follows under the complete smooth extension hypothesis of Theorem
\ref{thm:splitting}.
\end{theorem}

\begin{proof}
Choose the positive-scalar subsequence from nondecay.
Theorem~\ref{thm:intro-curvature}, Proposition~\ref{prop:unwrap},
Theorem~\ref{thm:splitting} and Theorem~\ref{thm:general} give the local
product with $K_{V}\geq 0$ and
$2K_{V}(q_{\infty })=\lim_i R^T(p_i)$. Under
\eqref{LS},
Theorem~\ref{thm:splitting} gives the equations. The identity
$dK_{V}=K_{V}\,d\psi $ propagates positivity on a connected smaller chart.
The additional roundness statement is the last assertion of
Theorem~\ref{thm:splitting}. Its necessity is demonstrated in
Example~\ref{ex:counterexample}.
\end{proof}

However, globally the situation is different and roundness may fail as in
Example \ref{ex:counterexample}.

\subsection{Globally Developable Five-Dimensional Models}

\label{sec:nondev-models}

The following are the globally developable (Definition~\ref{def:globdev})
quasi-regular and irregular shrinkers.

\begin{example}[Gaussian suspension]
\label{ex:suspension} The Heisenberg model $\mathbb{H}_{2}$ has Reeb field
$\xi _{0}=\partial _{t}$, structure tensor $\Phi _{0}$,
and for $\langle u,v\rangle :={\textstyle\sum_{a}}\overline{u}_{a}v_{a}$
\begin{equation*}
\eta _{0}=dt+\operatorname{Im}\langle z,dz\rangle ,\quad g_{0}=|dz|^{2}+\eta
_{0}^{2},\quad f_{0}=\tfrac{\lambda }{2}|z|^{2}.
\end{equation*}
For $A\in U(2)$ and $\ell >0,$ let
\begin{equation}
\sigma _{A,\ell }(z,t)=(Az,t+\ell )\quad \text{and}\quad M_{A,\ell
}=\mathbb{H}_{2}/\langle \sigma _{A,\ell }\rangle  \label{eq:suspension}
\end{equation}
with the structure induced by $(\eta _{0},\Phi _{0},g_{0},f_{0})$. Thus
$M_{I,\ell }$ is the regular Gaussian lift with period $\ell $.

\begin{enumerate}
\item if $A\neq I$ has finite order, $M_{A,\ell }$ is a quasi-regular finite
quotient of a regular Gaussian lift;

\item if $A$ has infinite order, $M_{A,\ell }$ is irregular
(Proposition~\ref{prop:transflat}).
\end{enumerate}
\end{example}

\begin{proposition}
\label{prop:transflat}Let $(M^{5},\xi ,\eta ,\Phi ,g,f)$ be a connected
complete gradient shrinking Sasaki--Ricci soliton with constant $\lambda >0$,
normalized by \eqref{eq:pre-normal}, and assume $R^{T}\equiv 0$.

\begin{enumerate}
\item $\mathrm{Rm}^{T}\equiv 0$, and the universal cover of $M$, with the
lifted structure and potential, is isomorphic to the Heisenberg model
$(\mathbb{H}_{2},\eta _{0},g_{0},f_{0})$ of \eqref{eq:heisenberg}.

\item $M$ is isomorphic either to $\mathbb{H}_{2}$ or to a Gaussian
suspension $M_{A,\ell }$, $A\in U(2)$, $\ell >0$. Moreover $M_{A,\ell }\cong
M_{A^{\prime },\ell ^{\prime }}$ if and only if $\ell =\ell ^{\prime }$ and
$A^{\prime }$ is conjugate to $A$ in $U(2)$. The Reeb leaves of
$\mathbb{H}_{2}$ are lines, $M_{A,\ell }$ is quasi-regular if and only if $A$
has
finite order, and otherwise it is in case $(T)$.

\item Let $A$ have finite order $m$. The central leaf $\{z=0\}$ of $M_{A,\ell
}$ has period $\ell $. The leaf through $[z,t]$ has period $k\ell
, $ where $k$ is the least positive integer with $A^{k}z=z$, the generic
period is $m\ell $, and the leaf space is the developable orbifold
\begin{equation*}
\lbrack \mathbb{C}^{2}/\langle A\rangle ].
\end{equation*}
In particular $M_{A,\ell }$ is regular if and only if $A=I$.
\end{enumerate}
\end{proposition}

\begin{proof}
\emph{Step 1.} By \cite[Proposition~\Iref{prop:flat}{3.2}]{PartI},
$R^{T}\equiv 0$ implies
\begin{equation*}
\mathrm{Ric}^{T}\equiv 0\quad \text{and}\quad \mathrm{Rm}^{T}\equiv 0.
\end{equation*}
Hence $\operatorname{Hess}^{T}f=\lambda h$. A Sasakian manifold whose transverse
metric is flat has constant $\Phi $-holomorphic sectional curvature $-3$,
and by Tanno's classification \cite{Tan} a complete simply connected
Sasakian manifold of constant $\Phi $-holomorphic sectional curvature $-3$
is isomorphic to the Heisenberg group with its standard Sasakian structure.
With the normalization $d\eta =2\omega ^{T}$, this is $(\mathbb{H}_{2},\eta
_{0},g_{0})$ of \eqref{eq:heisenberg}.

\emph{Step 2.} On $\mathbb{H}_{2}$ the lifted potential $\widehat{f}$ is
basic with $\operatorname{Hess}^{T}\widehat{f}=\lambda h_{\mathrm{E}}$, so
\begin{equation*}
\widehat{f}=\frac{\lambda }{2}|z-z_{\ast }|^{2}+b.
\end{equation*}
Since $R^{T}=0$, the identity \eqref{eq:pre-normal} reads $|\nabla
^{T}\widehat{f}|^{2}=2\lambda \widehat{f}$ and forces $b=0$. A left
translation
of the Heisenberg group, which is an automorphism of $(\eta _{0},g_{0})$,
moves $z_{\ast }$ to the origin. This proves (1). Steps~1 and~2 hold
verbatim in every dimension $2n+1$. For $n=1,$ they identify a simply
connected transversely flat soliton with $\mathbb{H}_{1}$.

\emph{Step 3.} An automorphism $\phi $ of $(\mathbb{H}_{2},\eta _{0},\Phi
_{0},g_{0},f_{0})$ commutes with the Reeb flow and induces a holomorphic
isometry $z\mapsto Az+w$ of the leaf space $\mathbb{C}^{2}$ preserving
$\frac{\lambda }{2}|z|^{2}$. Hence $w=0$ and $A\in U(2)$. Writing $\phi
(z,t)=(Az,t+s(z))$, the condition $\phi ^{\ast }\eta _{0}=\eta _{0}$ gives
$ds=0$, because $\operatorname{Im}\langle Az,A\,dz\rangle =\operatorname{Im}\langle
z,dz\rangle $. Thus
\begin{equation*}
\operatorname{Aut}(\mathbb{H}_{2},f_{0})=\{(z,t)\mapsto (Az,t+s):A\in U(2),\ s\in
\mathbb{R}\}.
\end{equation*}
The deck group $\Gamma \cong \pi _{1}(M)$ is a discrete subgroup of this
group acting freely. An element $(A,s)$ with $s=0$ fixes the points $(0,t)$,
so $(A,s)\mapsto s$ is injective on $\Gamma $. Its image is a discrete
subgroup of $\mathbb{R}$, since $U(2)$ is compact and $\Gamma $ acts
properly discontinuously. Hence either $\Gamma $ is trivial and $M\cong
\mathbb{H}_{2}$, or $\Gamma =\langle \sigma _{A,\ell }\rangle $ for a unique
generator with $\ell >0$, and $M\cong M_{A,\ell }$. An isomorphism $M_{A,\ell
}\cong M_{A^{\prime },\ell ^{\prime }}$ lifts to an element $(B,s_{0})$ of
$\operatorname{Aut}(\mathbb{H}_{2},f_{0})$ conjugating the deck
groups and their positive generators, so $\ell ^{\prime }=\ell $ and
$A^{\prime }=BAB^{-1}$. The converse is clear.

\emph{Step 4.} The Reeb flow of $M_{A,\ell }$ is $[z,t]\mapsto \lbrack
z,t+s] $, and $[z,t+s]=[z,t]$ if and only if $s=k\ell $ and $A^{k}z=z$ for
some $k\in \mathbb{Z}$. Hence the leaf through $[z,t]$ is closed if and only
if $A^{k}z=z$ for some $k\geq 1$, with minimal period $k\ell $ for the least
such $k$. If $A$ has finite order $m$, all leaves are closed and the central
leaf has period $\ell $. The leaves through the complement of the finitely
many proper subspaces $\ker (A^{k}-I)$, $0<k<m$, have period $m\ell $, and
the leaf space is $[\mathbb{C}^{2}/\langle A\rangle ]$ which is a
developable orbifold. The structure is regular if and only if $m=1$. If $A$
has infinite order, the leaves through points $z$ with $A^{k}z\neq z$ for
all $k\neq 0$ are not closed, and the closure of the Reeb flow is the
compact abelian group generated by the Reeb flow and the closure of $\langle
A\rangle $ in $U(2)$ which is a torus of dimension at least two. Thus $M$ is
in case $(T)$.
\end{proof}

\begin{theorem}
\label{thm:globdev}Let $(M^{5},\xi ,\eta ,\Phi ,g,f)$ be a complete
noncompact gradient shrinking Sasaki--Ricci soliton whose Reeb foliation is
globally developable. Then exactly one of the following holds.

\begin{enumerate}
\item $M$ is quasi-regular and its orbifold leaf space is developable, so
Theorem~\ref{thm:good} applies.

\item $M$ is in the proper-line case and is isomorphic to $\mathbb{H}_{2}$.

\item $M$ is in case $(T)$ and is isomorphic to an irrational Gaussian
suspension $M_{A,\ell }$.
\end{enumerate}

\noindent In cases (2) and (3),
\begin{equation*}
R^{T}\equiv 0.
\end{equation*}
\end{theorem}

\begin{proof}
Let $\mathcal{D}:\widehat{M}\rightarrow \widehat{Z}$ be as in Definition
\ref{def:globdev}, and let $\Gamma =\pi _{1}(M)$. Since $\mathcal{D}$ is a
Riemannian submersion from a complete manifold, it is a locally trivial
fibration by Hermann's theorem \cite{Her}. Thus either all lifted leaves are
circles or all are lines.

(i) Suppose they are circles. Proposition~\ref{prop:second-structure},
applied to the complete Sasakian manifold $\widehat{M}$, shows that the
lifted Reeb flow is periodic, $\widehat{\varphi }_{L}=\mathrm{Id}$. Deck
transformations commute with it, so $\varphi _{L}=\mathrm{Id}$ on $M$,
and $M $ is quasi-regular. Upstairs the circle action is actually free: a
fiber
foliation of a submersion has trivial leaf holonomy, whereas an exceptional
orbit in an effective locally free Reeb circle action has nontrivial
transverse holonomy. Hence all upstairs periods equal $L$. The homotopy
sequence of $\mathcal{D}$ shows that $\widehat{Z}$ is simply connected.
Since $\mathcal{D}$ is proper (a locally trivial fibration with circle
fibers) and $\Gamma $-equivariant, $\Gamma $ acts properly on $\widehat{Z}$
with finite stabilizers, and the map $\widehat{Z}\rightarrow
\mathcal{Z}=M/\mathbb{S}^{1}$ induced by the foliated covering
$\widehat{M}\rightarrow M$
identifies $\mathcal{Z}$ with the orbifold quotient $[\widehat{Z}/\Gamma ]$,
up to the kernel of the action. Hence $\mathcal{Z}$ is developable, and (1)
holds.

(ii) Suppose they are lines. Then the lifted Reeb flow acts freely, and
properly because its orbits are the fibers of the locally trivial fibration
$\mathcal{D}$. Corollary~\ref{cor:Rcase}, applied to $\widehat{M}$, gives
$\widehat{M}\cong \mathbb{H}_{2}$, so $R^{T}\equiv 0$ on $M$. By Proposition
\ref{prop:transflat}, $M$ is $\mathbb{H}_{2}$ or a suspension $M_{A,\ell }$.
The latter is quasi-regular with developable leaf space if $A$ has finite
order, which is (1), and is in case (T) if $A$ has infinite order, which is
(3). The three cases are distinguished by
Proposition~\ref{prop:second-structure}.
\end{proof}

\begin{corollary}
\label{cor:globdev-class}A complete noncompact five-dimensional gradient
shrinking Sasaki--Ricci soliton with globally developable Reeb foliation is
isomorphic to exactly one of the following:

\begin{enumerate}
\item a finite free quotient $\widetilde{M}/\Gamma $, as in
\eqref{eq:finitequotient}, of a regular lift $\widetilde{M}$ of one of the
four models \eqref{eq:fourmodels};

\item the Heisenberg model $\mathbb{H}_{2}$; or

\item an irrational Gaussian suspension $M_{A,\ell }$, determined by $\ell $
and the conjugacy class of $A$.
\end{enumerate}

\noindent The FIK, round-cylinder and BCCD models occur only in $(1)$, with
generic Reeb period at most $8\pi /\lambda $.
\end{corollary}

\begin{proof}
Combine Theorems~\ref{thm:globdev} and \ref{thm:good},
Proposition~\ref{prop:transflat}, and Lemma~\ref{lem:period-bound} applied to
$\widetilde{M}$, whose Reeb period is the generic period of $M$ by the proof
of Theorem~\ref{thm:good}.
\end{proof}

\begin{theorem}
\label{thm:pi1}Let $(M^{5},\xi ,\eta ,\Phi ,g,f)$ be a connected complete
noncompact gradient shrinking Sasaki--Ricci soliton. Let
$p:\widehat{M}\rightarrow M$ be its universal covering, with the lifted
structure and
potential $\widehat{f}=f\circ p$ and deck group $\Gamma \cong \pi _{1}(M)$,
and let $\widehat{H}$ be the closure of the lifted Reeb flow in
$\operatorname{Isom}(\widehat{M},\widehat{g})$. Exactly one of the following holds.

\begin{enumerate}
\item $\widehat{H}$ is not compact. Then the lifted Reeb flow is a free
proper $\mathbb{R}$-action, $\widehat{M}\cong \mathbb{H}_{2}$, and $M$ is
isomorphic to $\mathbb{H}_{2}$ or to a Gaussian suspension $M_{A,\ell }$;
in particular $R^{T}\equiv 0$.

\item $\widehat{H}$ is a torus. Then $M$ is in case $(Q)$ or $(T)$,
$\widehat{f}$ is proper on $\widehat{M}$, $\Gamma $ is finite, and $\Gamma $
commutes with $\widehat{H}$. The covering $\widehat{H}\rightarrow H$ induced
by $p$ is surjective with finite kernel, so $\dim \widehat{H}=\dim H$, and
$\widehat{M}$ is quasi-regular if and only if $M$ is.
\end{enumerate}

\noindent Consequently $\pi _{1}(M)$ is infinite if and only if $M$ is
isomorphic to a Gaussian suspension $M_{A,\ell }$ for some $A\in U(2)$ and
$\ell >0$, and then $\pi _{1}(M)\cong \mathbb{Z}$. In particular, if
$R^{T}\not\equiv 0$, then $\pi _{1}(M)$ is finite.
\end{theorem}

\begin{proof}
By Lemma~\ref{lem:reeb}, applied to the connected complete Sasakian manifold
$\widehat{M}$, either $\widehat{H}$ is a compact torus or the lifted Reeb
flow is a closed, free and proper $\mathbb{R}$-action. These are the two
cases and they are exclusive.

(i) The lifted soliton is complete, noncompact and normalized. Corollary
\ref{cor:Rcase}, applied to $\widehat{M}$, shows that it is transversely flat.
Proposition~\ref{prop:transflat}, applied to the simply connected
$\widehat{M}$, gives $\widehat{M}\cong \mathbb{H}_{2}$. Since $p$ is a local
isomorphism, $R^{T}\equiv 0$ on $M$, and Proposition~\ref{prop:transflat}
gives the assertion on $M$.

(ii) By Lemma~\ref{lem:growth}, in the compact-torus case for $\widehat{M}$,
the function $\widehat{f}$ is proper. It is invariant under $\Gamma $, which
acts properly discontinuously. For $x\in \widehat{M},$ the orbit $\Gamma x$
lies in the compact set $\{\widehat{f}\leq \widehat{f}(x)\}$. So $\Gamma $
is finite as in the proof of Lemma~\ref{lem:finite}. Deck transformations
commute with the lifted Reeb flow, hence with its closure $\widehat{H}$.
Thus every $h\in \widehat{H}$ descends to an isometry $\rho (h)$ of $M$ with
$p\circ h=\rho (h)\circ p$, and $\rho :\widehat{H}\rightarrow
\operatorname{Isom}(M,g)$ is a continuous homomorphism, since compact subsets
of $M$ are images of compact subsets of $\widehat{M}$. Its image is compact, contains the Reeb
flow of $M$, and is contained in the closure $H$ of that flow. Then $\rho
(\widehat{H})=H$. Its kernel $\widehat{H}\cap \Gamma $ is finite. Hence $H$
is compact, $M$ is in case $(Q)$ or $(T)$, and $\dim \widehat{H}=\dim H$. If
$\varphi _{L}=\mathrm{Id}$ on $M$, then $\widehat{\varphi }_{L}\in
\Gamma $
has finite order $m$ and $\widehat{\varphi }_{mL}=\mathrm{Id}$.
Conversely $\widehat{\varphi }_{L}=\mathrm{Id}$ implies $\varphi
_{L}\circ p=p\circ
\widehat{\varphi }_{L}=p$.

(iii) For the last assertions: if $\pi _{1}(M)$ is infinite, (2) fails, so
(1) holds. The Heisenberg model $\mathbb{H}_{2}$ is simply connected, so
$M\cong M_{A,\ell }$, whose fundamental group is the infinite cyclic group
generated by $\sigma _{A,\ell }$. Conversely $\pi _{1}(M_{A,\ell })\cong
\mathbb{Z}$. Finally (1) implies $R^{T}\equiv 0$.
\end{proof}

As a consequence, for quasi-regular shrinkers with nondevelopable leaf space
we have the following.

\begin{corollary}
\label{cor:nondev-qr}Let $M$ be a complete noncompact quasi-regular
five-dimensional gradient shrinking Sasaki--Ricci soliton whose orbifold
leaf space is not developable. Then

\begin{enumerate}
\item $R^{T}>0$, $\pi _{1}(M)$ is finite, and $M=\widehat{M}/\Gamma $, where
the universal cover $\widehat{M}$ is quasi-regular with a nondevelopable
leaf space.

\item $\Gamma \cong \pi _{1}(M)$ is a finite group of automorphisms of the
lifted soliton which acts freely and commutes with the circle action. In
particular $\widehat{M}$ is not regular.
\end{enumerate}
\end{corollary}

\begin{proof}
By Proposition~\ref{prop:transflat}, a transversely flat quasi-regular
soliton is a suspension $M_{A,\ell }$ with $A$ of finite order, whose leaf
space $[\mathbb{C}^{2}/\langle A\rangle ]$ is developable. Hence
$R^{T}\not\equiv 0$, and $R^{T}>0$ by Lemma~\ref{lem:nonneg}.
Theorem~\ref{thm:pi1}(2) gives the finiteness of $\Gamma $, its commutation
with the
lifted circle action, and the quasi-regularity of $\widehat{M}$.

Suppose the leaf space of $\widehat{M}$ were developable.
Theorem~\ref{thm:good}, applied to the simply connected $\widehat{M}$, would
show that $\widehat{M}$ itself is regular, with a smooth quotient
$\widehat{Z}$ which
is simply connected by the homotopy sequence. The finite group $\Gamma $
commutes with the free circle action, so it acts on $\widehat{Z}$ by
holomorphic isometries. As in the proof of Theorem~\ref{thm:globdev}, the
leaf space of $M$ is the orbifold quotient $[\widehat{Z}/\Gamma ]$, up to
the kernel of the action. This quotient is developable, a contradiction.
\end{proof}

The following example gives, for $a_{1}/a_{2}\notin \mathbb{Q}$, irregular
shrinkers whose Reeb foliation is not globally developable
(Proposition~\ref{prop:S3-family}).

\begin{example}
\label{ex:S3}\label{ex:counterexample} Represent $\mathbb{S}^{3}$ by
$(t,\theta _{1},\theta _{2})$, $0\leq t\leq 1$, with both angles of period
$2\pi $, collapsing the $\theta _{2}$ circle at
$t=0$ and the $\theta _{1}$ circle at $t=1$. Choose $a=(a_{1},a_{2})$ with
$a_{1},a_{2}>0$, set $\lambda =2(a_{1}+a_{2})$, and let $c=c(a)$ be the
constant of \eqref{eq:Hc}; thus $c=0$ if and only if $a_{1}=a_{2}$, and $c>0$
if and only if $a_{1}>a_{2}$. There is a smooth Sasakian structure $S_{a}$
on $\mathbb{S}^{3}$ with Reeb field
\begin{equation*}
\xi _{S}=a_{1}\partial _{\theta _{1}}+a_{2}\partial _{\theta _{2}}
\end{equation*}
and a basic potential $\psi $ satisfying
\begin{equation*}
\mathrm{Ric}_{S}^{T}+\operatorname{Hess}_{S}^{T}\psi =\lambda h_{S}^{T}.
\end{equation*}
Its transverse Gaussian curvature $K$ is positive, and nonconstant unless
$a_{1}=a_{2}$; compare Wang--Zhang \cite{WZ2}, who proved the existence and
uniqueness of gradient Sasaki--Ricci solitons on Sasakian three-spheres. Set
\begin{equation*}
M=\mathbb{C}\times \mathbb{S}^{3},
\end{equation*}
write $z=x+iy$, and define
\begin{gather}
\alpha =x\,dy-y\,dx,\qquad \eta =\eta _{S}+\alpha ,\qquad \xi =\xi _{S},
\label{eq:ctr-eta} \\
g=dx^{2}+dy^{2}+h_{S}^{T}+(\eta _{S}+\alpha )^{2},\qquad f=\frac{\lambda
}{2}(x^{2}+y^{2})+ct+\frac{c}{e^{c}-1},  \label{eq:ctr-g}
\end{gather}
where $c/(e^{c}-1)$ is read as $1$ for $c=0$. Then $\mathcal{S}_{a}:=(M,\xi
,\eta ,\Phi ,g,f)$, with $\Phi $ induced by the product
transverse complex structure, is a complete noncompact gradient shrinking
Sasaki--Ricci soliton with shrinking constant $\lambda $, proper basic
potential, and Hamilton identity $R^{T}+|\nabla ^{T}f|^{2}=2\lambda f$. It
is irregular if and only if $a_{1}/a_{2}\notin \mathbb{Q}$
(Proposition~\ref{prop:S3-family}). Its transverse scalar curvature is
\begin{equation}
R^{T}(z,t)=2K(t)=\frac{2\lambda ce^{ct}}{e^{c}-1}>0,  \label{eq:ctr-RT}
\end{equation}
read as $2\lambda $ for $c=0$; it is bounded, and nonconstant unless
$a_{1}=a_{2}$. Explicit formulas are given in
\eqref{eq:V}--\eqref{eq:S3metric} in Appendix~\ref{app:A}.

Moreover, let $a_{1}\neq a_{2}$, let $q_{\ast }\in \mathbb{S}^{3}$ be a
point of the endpoint orbit on which $K$ attains its maximum, and put
$p_{r}=(r,q_{\ast })$ for $r\geq 0$. Then $r\mapsto p_{r}$ is, after
reparametrization, an escaping integral curve of $\nabla ^{T}f$ along which
$R^{T}\equiv 2K(q_{\ast })>2\lambda $. The pointed manifolds $(M,p_{r})$ are
all isomorphic as pointed Sasakian manifolds, by isomorphisms which change $f
$ by affine functions of $z$. Hence every developed local limit along this
curve is a local product of a flat plane and a surface of nonconstant
positive curvature which satisfies the local soliton equation \eqref{LS},
and no neighborhood of its basepoint is transversely isometric to the round
K\"{a}hler cylinder (Appendix~\ref{app:A}). Thus the roundness assertion of
Theorem~\ref{thm:splitting} needs its extension hypothesis, and the
regularity assumption in Theorem~\ref{T51} cannot be dropped; compare
Theorem~\ref{thm:irregular-local} and Remark~\ref{rem:weighted-example}.
\end{example}

\subsection{Quasi-Regular Soliton Approximation}

\label{sec:qr-approx}

In this subsection we treat Reeb foliations which are not
quasi-regular by the
strategy used for compact transverse Fano Sasakian five-manifolds in
\cite[Section~4]{CHLWHT}. The compact second structure theorem \cite{Ru,BG}
has
the conditional complete-manifold extension in
Proposition~\ref{prop:second-structure}. A soliton approximation, when
independently
constructed, can be controlled by a compactness theorem, and the structure
theorem of the quasi-regular case (Theorem~\ref{thm:good} and
Corollary~\ref{cor:four-models}) is passed to the limit.

Two features of the complete noncompact setting require changes. First,
when a quasi-regular sequence converges to a limit that is not
quasi-regular, its generic Reeb periods tend to infinity
(Lemma~\ref{lem:period-blowup}). In the explicit examples, the orders of
orbifold isotropy groups or of the finite quotient groups also become
unbounded
(Propositions~\ref{prop:suspension-approx} and \ref{prop:S3-family}). Since
the noncollapsing theory of
\cite[Section~\Iref{sec:noncollapsing}{6.2}]{PartI} is formulated on a
complete smooth quotient (setting (S) there),
compactness is formulated here for the five-dimensional total spaces, with
an injectivity-radius hypothesis on the Sasakian metric.

Second, the deformation of the Reeb field supplied by the second structure
theorem does not preserve the soliton equation and no Sasaki--Ricci flow
limit theorem is available here to restore it. The approximation is
therefore formulated at the level of solitons (Definition~\ref{def:qr-approx}).

\begin{definition}
\label{def:pointed-conv}We consider pointed complete Sasakian manifolds with
basic functions. The sequence $\mathcal{M}_{i}$ \emph{converges} to
$\mathcal{M}_{\infty }$ as
\begin{equation*}
\mathcal{M}_{i}=(M_{i},\xi _{i},\eta _{i},\Phi
_{i},g_{i},f_{i},p_{i})\rightarrow \mathcal{M}_{\infty }=(M_{\infty },\xi
_{\infty },\eta _{\infty },\Phi _{\infty },g_{\infty },f_{\infty },p_{\infty
})
\end{equation*}
if there are an exhaustion of $M_{\infty }$ by relatively compact open sets
$U_{i}\ni p_{\infty }$ and smooth embeddings
\begin{equation*}
F_{i}:U_{i}\rightarrow M_{i}
\end{equation*}
with $F_{i}(p_{\infty })=p_{i}$ such that
\begin{equation*}
F_{i}^{\ast }g_{i}\rightarrow g_{\infty },\qquad F_{i}^{\ast }\xi
_{i}\rightarrow \xi _{\infty },\qquad F_{i}^{\ast }f_{i}\rightarrow f_{\infty
}\quad \text{in}\quad
C^{\infty }
\end{equation*}
on every fixed compact subset of $M_{\infty }$ once it lies in $U_i$.
For vector fields, the pullback means
$F_i^*\xi_i=(dF_i)^{-1}(\xi_i\circ F_i)$. The functions $f_i$ here have
their stated additive normalization; no further constants are subtracted
in this definition. Then also
\begin{equation*}
F_{i}^{\ast }\eta _{i}\rightarrow \eta _{\infty }\quad \text{and}\quad
F_{i}^{\ast
}\Phi _{i}\rightarrow \Phi _{\infty },
\end{equation*}
since $\eta =g(\xi ,\cdot )$ and $\Phi =\nabla \xi $.
\end{definition}

\begin{definition}
\label{def:qr-approx}Let $(M^{5},\xi ,\eta ,\Phi ,g,f)$ be a complete
noncompact gradient shrinking Sasaki--Ricci soliton with constant $\lambda>0$,
normalized by \eqref{eq:pre-normal}.
A \emph{quasi-regular soliton approximation} of $M$ consists of
quasi-regular complete noncompact five-dimensional gradient shrinking
Sasaki--Ricci solitons $(M_{i},\xi _{i},\eta _{i},\Phi _{i},g_{i},f_{i})$
with the same constant $\lambda $, normalized by \eqref{eq:pre-normal}, and
of points $p_{i}\in M_{i}$ and $p\in M$, such that $(M_{i},\xi _{i},\eta
_{i},\Phi _{i},g_{i},f_{i},p_{i})$ converges to $(M,\xi ,\eta ,\Phi ,g,f,p)$
in the sense of Definition~\ref{def:pointed-conv}. It is \emph{developable}
if every orbifold leaf space $\mathcal{Z}_{i}=M_{i}/\mathbb{S}^{1}$ is
developable in
the sense of Definition~\ref{def:developable}.
\end{definition}

\begin{remark}
\label{rem:compact-fano}For compact transverse Fano Sasakian five-manifolds
the method of \cite[Section~4]{CHLWHT} starts from the quasi-regular
structures of Proposition~\ref{prop:second-structure}(3)(d), which exist
because $M$ is compact. The Sasaki--Ricci flow, its stability under
perturbation of the Reeb field \cite{He}, and a compactness theorem for
quasi-regular transverse Fano Sasaki--Ricci solitons, based on the
transverse $\alpha $-invariant, log canonical thresholds and Birkar's
boundedness theorem \cite{Bir} and modeled on Guo--Phong--Song--Sturm
\cite{GPSS}, then produce quasi-regular solitons converging to the limit of
the
irregular flow. In the complete noncompact setting three points differ:
quasi-regular approximations of the Reeb field need a bounded contact moment
map (Proposition~\ref{prop:second-structure}(3)(d)); the leaf spaces of the
approximations degenerate, since their generic Reeb periods tend to infinity
(Lemma~\ref{lem:period-blowup}); and, since no flow argument is available,
the approximating structures have to be solitons themselves
(Definition~\ref{def:qr-approx}).
\end{remark}

We now prove a compactness theorem for five-dimensional shrinkers.

\begin{theorem}
\label{thm:compactness-total}Let $(M_{i}^{5},\xi _{i},\eta _{i},\Phi
_{i},g_{i},f_{i})$ be connected complete gradient shrinking Sasaki--Ricci
solitons with the same constant $\lambda >0$, normalized by
\eqref{eq:pre-normal}, and let $p_{i}\in M_{i}$. Assume that there are
constants $\Lambda ,i_{0},C_{0}$ with $f_{i}(p_{i})\leq C_{0}$ and
\begin{equation}
\sup_{M_{i}}|\mathrm{Rm}_{g_{i}}^{T}|\leq \Lambda ,\qquad
\operatorname{inj}_{g_{i}}(p_{i})\geq i_{0}>0.  \label{eq:compactness-hyp}
\end{equation}
Then a subsequence converges, in the sense of
Definition~\ref{def:pointed-conv}, to a complete gradient shrinking
Sasaki--Ricci soliton $(M_{\infty },\xi _{\infty },\eta _{\infty },\Phi
_{\infty },g_{\infty
},f_{\infty },p_{\infty })$ with constant $\lambda $, normalized by
\eqref{eq:pre-normal}, with $\sup |\mathrm{Rm}_{g_{\infty }}^{T}|\leq
\Lambda $. If the $M_{i}$ are noncompact, so is $M_{\infty }$.
\end{theorem}

\begin{proof}
After a fixed $D$-homothety we may assume $\lambda =\frac{1}{2}$. The
associated flows
\begin{equation*}
g_{i}^{T}(t)=(-t)\Psi _{i,t}^{\ast }g_{i}^{T},
\end{equation*}
$t\in \lbrack -2,-1]$, have transverse curvature bounded in terms of $\Lambda
$, so the local Shi-type estimates \cite{Shi,Co1} bound
all transverse covariant derivatives of $\mathrm{Rm}_{g_{i}}^{T}$ at $t=-1$
in terms of $\Lambda $. By the Sasakian identities the full curvature of
$g_{i}$ and all its covariant derivatives are then uniformly bounded (cf.\
\cite[Proposition~\Iref{prop:unwrap}{5.3}]{PartI}). With the
injectivity-radius bound, the Cheeger--Gromov compactness theorem gives a
subsequence converging to a complete pointed Riemannian manifold $(M_{\infty
},g_{\infty },p_{\infty })$. The unit Killing fields $\xi _{i}$ have
uniformly bounded covariant derivatives and converge to a unit Killing field
$\xi _{\infty }$. Then $\eta _{i}$ and $\Phi _{i}=\nabla \xi _{i}$ converge
(Definition~\ref{def:pointed-conv}), and the Sasakian identities pass to the
limit. Since
\begin{equation*}
|\nabla ^{T}f_{i}|^{2}\leq 2\lambda f_{i},\quad f_{i}(p_{i})\leq C_{0}
\end{equation*}
and
\begin{equation*}
\operatorname{Hess}^{T}f_{i}=\lambda h_{i}-\mathrm{Ric}_{g_{i}}^{T},
\end{equation*}
the potentials are locally uniformly bounded with all derivatives and
converge to a basic function $f_{\infty }$ satisfying the soliton equation
and \eqref{eq:pre-normal}. If $M_{\infty }$ were compact, then for large $i$
the embedding $F_{i}$ would be defined on all of $M_{\infty }$, and its
image would be open and closed in the connected manifold $M_{i}$, so $M_{i}$
would be compact.
\end{proof}

\begin{remark}
\label{rem:compactness-scope}

\begin{enumerate}
\item The injectivity-radius bound concerns the five-dimensional metric
$g_{i}$, not a leaf space. Under the curvature bound it may be replaced by a
lower bound $\operatorname{Vol}_{g_{i}}B_{g_{i}}(p_{i},1)\geq v_{0}>0$, by
Cheeger--Gromov--Taylor \cite{CGT}. For a regular structure, it fails if the
Reeb periods tend to zero, while Reeb periods tending to infinity are
allowed. If a quasi-regular sequence has a limit which is not quasi-regular,
its generic periods necessarily diverge by Lemma~\ref{lem:period-blowup}(2).

\item The Reeb type is not preserved under convergence, in either direction.
The regular lifts of the Gaussian shrinker converge to the Heisenberg model
as their period tends to infinity, quasi-regular Gaussian suspensions
converge to irregular ones (Proposition~\ref{prop:suspension-approx}), and
irregular solitons on $\mathbb{C}\times \mathbb{S}^{3}$ converge to the
regular cylinder lift (Proposition~\ref{prop:S3-family}(4)).
\end{enumerate}
\end{remark}

\begin{lemma}
\label{lem:period-blowup}Let $\mathcal{M}_{i}\rightarrow \mathcal{M}_{\infty
}$ in the sense of Definition~\ref{def:pointed-conv}, with the notation
there, and assume that every $M_{i}$ is quasi-regular, with generic Reeb
period $L_{i}$.

\begin{enumerate}
\item If the Reeb leaf of $M_{\infty }$ through $q$ is not compact and
$q_{i}=F_{i}(q)$, then the minimal periods $\tau _{i}(q_{i})$ of the leaves
through $q_{i}$ tend to infinity. In particular $L_{i}\rightarrow \infty $.

\item If $\liminf_{i}L_{i}<\infty $, then $M_{\infty }$ is quasi-regular,
and $\liminf_{i}L_{i}$ is an integer multiple of its generic Reeb period.
\end{enumerate}
\end{lemma}

\begin{proof}
Write $X_{i}=F_{i}^{\ast }\xi _{i}$ on $U_{i}$, and denote the Reeb flows by
$\varphi _{t}^{i}$ and $\varphi _{t}^{\infty }$. These flows are defined for
all real times, since the Reeb fields have unit length on complete
manifolds. We first establish convergence of the pulled-back flows on
bounded time intervals. Fix a compact set $K\subset M_{\infty }$ and $T>0$.
The set
\begin{equation*}
K_{T}=\{\varphi _{t}^{\infty }(q):q\in K,\ |t|\leq T+1\}
\end{equation*}
is compact. Choose a relatively compact open neighborhood $U$ of $K_{T}$ and
$\delta >0$ such that
\begin{equation*}
d_{g_{\infty }}(K_{T},M_{\infty }\setminus U)>2\delta .
\end{equation*}
For all sufficiently large $i$, $\overline{U}\subset U_{i}$ and
$X_{i}\rightarrow \xi _{\infty }$ smoothly on $\overline{U}$. Uniform
continuous dependence for ordinary differential equations, in finitely many
coordinate charts covering $K_{T}$, shows that the local $X_{i}$-flows $\psi
_{t}^{i}$ with initial points in $K$ exist for $|t|\leq T$, remain in $U$,
and converge uniformly on $[-T,T]\times K$ to $\varphi _{t}^{\infty }$. By
uniqueness of solutions,
\begin{equation*}
F_{i}(\psi _{t}^{i}(q))=\varphi _{t}^{i}(F_{i}(q))\qquad (q\in K,\ |t|\leq
T).
\end{equation*}

We also need a positive lower bound for periods at each fixed $q\in
M_{\infty }$. Choose a flow-box coordinate $u$ near $q$ with $du(\xi
_{\infty })>0$. On a smaller relatively compact neighborhood, this
derivative is bounded below by a positive constant, and so is $du(X_{i})$
for all large $i$. Choose $\varepsilon _{q}>0$ so that the limit trajectory
from $q$ for $0\leq t\leq \varepsilon _{q}$ stays in the interior of that
neighborhood. The preceding bounded-time argument, or its local first-exit
proof, keeps the corresponding $X_{i}$-trajectories there for all large $i$.
The coordinate $u$ is strictly increasing on them. Since $F_{i}$ is
injective, no positive period of the Reeb orbit through $F_{i}(q)$ is
smaller than $\varepsilon _{q}$.

(1) If $\tau_i(F_i(q))$ were bounded along a subsequence, choose a further
subsequence with $\tau_i(F_i(q))\to\tau_*$, where $\tau_*\geq\varepsilon_q>0$.
For a time bound larger than these periods, the identity above and
injectivity of $F_i$ give $\psi^i_{\tau_i(F_i(q))}(q)=q$. Uniform
bounded-time convergence, together with continuity of the limit flow in
time, yields $\varphi^\infty_{\tau_*}(q)=q$. The limit leaf is therefore
compact, contrary to the hypothesis. Thus the minimal periods tend to
infinity. Each minimal period divides the generic period, so
$L_i\geq\tau_i(F_i(q))\to\infty$ as well.

(2) Choose a subsequence with $L_{i}\rightarrow L_{\ast
}:=\liminf_{i}L_{i}<\infty $. The period lower bound at one fixed point
gives $L_{\ast }>0$. For any fixed $q\in M_{\infty }$, apply the
bounded-time convergence with a time bound larger than $L_{\ast }+1$. Since
$\varphi _{L_{i}}^{i}=\mathrm{Id}$, it follows that $\varphi _{L_{\ast
}}^{\infty }(q)=q$. As $q$ was arbitrary, $\varphi _{L_{\ast }}^{\infty
}=\mathrm{Id}$ on all of $M_{\infty }$. The global period kernel is a nonzero
closed subgroup of $\mathbb{R}$. It is not all of $\mathbb{R}$ because $\xi
_{\infty }$ never vanishes. It is therefore $L_{0}\mathbb{Z}$ for some
$L_{0}>0$, and the induced $\mathbb{R}/L_{0}\mathbb{Z}$ action is effective
and locally free. Thus the limit is quasi-regular. Since the circle group is
abelian, its principal stabilizer fixes the dense principal stratum and
hence the whole manifold. The effectiveness makes that stabilizer trivial.
Consequently $L_{0}$ is the generic minimal Reeb period and $L_{\ast }$ is a
positive integer multiple of $L_{0}$.
\end{proof}

\begin{lemma}
\label{lem:period-bound}Let $(M^{5},\xi ,\eta ,\Phi ,g,f)$ be a complete
noncompact regular gradient shrinking Sasaki--Ricci soliton with Reeb period
$L$, quotient $Z$ and integral Euler class $e$. For every compact
holomorphic curve $C\subset Z$,
\begin{equation}
L\leq \frac{4\pi }{\lambda }\,c_{1}(Z)\cdot \lbrack C].
\label{eq:period-bound}
\end{equation}
In particular, if $Z$ is the FIK, the round-cylinder or the BCCD model, then
\begin{equation*}
L\leq 8\pi /\lambda .
\end{equation*}
Only the Gaussian quotient admits regular lifts of arbitrary period.
\end{lemma}

\begin{proof}
By Proposition~\ref{prop:lifting},
\begin{equation*}
\langle e,[C]\rangle =\frac{2}{L}\int_{C}\omega _{Z}=\frac{4\pi }{\lambda
L}\,c_{1}(Z)\cdot \lbrack C].
\end{equation*}
The middle term is positive, since $\omega _{Z}$ restricts to an area form
on the complex curve $C$, and the left-hand side is an integer. Hence it is
at least one, which is \eqref{eq:period-bound}. The FIK and BCCD surfaces
contain the exceptional curve, with $c_{1}(Z)\cdot \lbrack E]=1$, and the
round cylinder contains $\{0\}\times \mathbb{CP}^{1}$, with $c_{1}(Z)\cdot
\lbrack \mathbb{CP}^{1}]=2$ (see the proofs of Theorems~\ref{thm:regular}
and~\ref{thm:decay-lifts}). This gives $L\leq 8\pi /\lambda $. The Gaussian
quotient admits every $L>0$ by Theorem~\ref{thm:decay-lifts}(1).
\end{proof}

\begin{corollary}
\label{cor:nongaussian}Let $\mathcal{M}_{i}\rightarrow \mathcal{M}_{\infty }$
in the sense of Definition~\ref{def:pointed-conv}, where the $M_i$ are
quasi-regular complete noncompact gradient shrinking Sasaki--Ricci solitons
with the same constant $\lambda$ and with developable orbifold leaf spaces. Let
$\widetilde{M}_{i}\rightarrow M_{i}$ be the finite coverings of
Theorem~\ref{thm:good}, whose lifted Reeb foliations are regular,
with smooth quotients $\widetilde{Z}_{i}$.

\begin{enumerate}
\item If $\widetilde{Z}_{i}$ is the FIK, the round-cylinder or the BCCD
model for infinitely many $i$, then $M_{\infty }$ is quasi-regular with
generic Reeb period at most $8\pi /\lambda $.

\item If $M_{\infty }$ is not quasi-regular, then for all large $i$ the
quotient $\widetilde{Z}_{i}$ is the Gaussian shrinker, $\widetilde{M}_{i}$
is the regular Gaussian lift with period $L_{i}\rightarrow \infty $, and
$\mathrm{Rm}^{T}\equiv 0$ on $M_{i}$ and on $M_{\infty }$.
\end{enumerate}
\end{corollary}

\begin{proof}
By the proof of Theorem~\ref{thm:good}, the Reeb period of $\widetilde{M}_{i}
$ is the generic period $L_{i}$ of $M_{i}$. (1) If $\widetilde{Z}_{i}$
is the FIK, round-cylinder or BCCD model, then $L_{i}\leq 8\pi /\lambda $ by
Lemma~\ref{lem:period-bound}. If this happens for infinitely many $i$,
Lemma~\ref{lem:period-blowup}(2) shows that $M_{\infty }$ is quasi-regular
with
generic period at most $8\pi /\lambda $. (2) If $M_{\infty }$ is not
quasi-regular, then by (1) the quotient $\widetilde{Z}_{i}$ is the Gaussian
shrinker for all large $i$, so $\widetilde{M}_{i}$ is a regular Gaussian
lift and $M_{i}$, being covered by it, is transversely flat. Moreover
$L_{i}\rightarrow \infty $, since a finite $\liminf_{i}L_{i}$ would make
$M_{\infty }$ quasi-regular by Lemma~\ref{lem:period-blowup}(2). Finally,
transverse flatness passes to the smooth limit.
\end{proof}

The following proposition gives quasi-regular soliton approximations of
Gaussian suspensions.

\begin{proposition}
\label{prop:suspension-approx}

\begin{enumerate}
\item The regular Gaussian lifts $M_{I,L}$ converge to $\mathbb{H}_{2}$ as
$L\rightarrow \infty $, in the sense of Definition~\ref{def:pointed-conv}.

\item Let $X\in \mathfrak{u}(2)$, $A=\exp X$, $\ell >0$, and
$N=\mathbb{C}^{2}\times (\mathbb{R}/\ell \mathbb{Z})$. The map $(z,t)\mapsto
(e^{tX/\ell
}z,t)$ descends to a diffeomorphism $N\rightarrow M_{A,\ell }$, under which
the structure of $M_{A,\ell }$ becomes
\begin{equation}
\begin{gathered} \eta _{X}=\bigl(1+\ell ^{-1}\operatorname{Im}\langle
z,Xz\rangle \bigr)dt+\operatorname{Im}\langle z,dz\rangle ,\qquad \xi
_{X}=\partial _{t}-\ell ^{-1}V_{X}, \\
g_{X}=\bigl|dz+\ell
^{-1}Xz\,dt\bigr|^{2}+\eta _{X}^{2},\qquad f=\tfrac{\lambda }{2}|z|^{2},
\end{gathered}  \label{eq:suspension-model}
\end{equation}
where $V_{X}$ is the linear vector field $z\mapsto Xz$. These tensors do not
depend on $t$ and depend polynomially on $X$. Hence $X_{i}\rightarrow X$
implies convergence of $(N,\xi _{X_{i}},\eta _{X_{i}},g_{X_{i}},f)$ to $(N,\xi
_{X},\eta _{X},g_{X},f)$ on the fixed manifold $N$.

\item If $A$ has infinite order, choose $U\in U(2)$ with $U^{\ast
}XU=\operatorname{diag}(2\pi i\theta _{1},2\pi i\theta _{2})$ and rational
pairs $r^{(i)}\rightarrow (\theta _{1},\theta _{2})$, and put
$X_{i}=U\operatorname{diag}(2\pi ir_{1}^{(i)},2\pi ir_{2}^{(i)})U^{\ast }$.
Then $A_{i}=\exp X_{i}$ has
finite order $m_{i}\rightarrow \infty $, the quasi-regular suspensions
$M_{A_{i},\ell }$ form a developable quasi-regular soliton approximation of
$M_{A,\ell }$, their generic periods $m_{i}\ell $ tend to infinity, and the
central leaf keeps the period $\ell $. The Reeb fields $\xi _{X_{i}}$ and $\xi
_{X}$ lie in the Lie algebra of the three-torus generated by $\partial
_{t}$ and $U\operatorname{diag}(e^{i\phi _{1}},e^{i\phi _{2}})U^{\ast }$,
which acts
on $N$ by automorphisms of each of these structures. The fields $\xi
_{X_{i}} $ are rational and converge to $\xi _{X}$.
\end{enumerate}
\end{proposition}

\begin{proof}
(1) Left translations of the Heisenberg group are automorphisms of $(\eta
_{0},g_{0})$ and commute with the central translations $(z,t)\mapsto (z,t+s)$,
so the latter move every point by the same distance $\delta (s)$, and $\delta
(s)\rightarrow \infty $ as $|s|\rightarrow \infty $. Hence for every $R$ the
quotient map $\mathbb{H}_{2}\rightarrow M_{I,L}$ is injective on
$B_{g_{0}}((0,0),R)$ once $\delta (s)>2R$ for all $|s|\geq L$, and it is an
isomorphism of the structures there.

(2) With $\Psi (z,t)=(e^{tX/\ell }z,t),$ one has $\sigma _{A,\ell }\circ
\Psi =\Psi \circ \sigma _{I,\ell }$. Since $e^{X}$ commutes with $e^{tX/\ell
}$, $\Psi $ descends to a diffeomorphism $N\rightarrow M_{A,\ell }$. For
$w=e^{tX/\ell }z,$ one has $dw=e^{tX/\ell }(dz+\ell ^{-1}Xz\,dt)$ and
unitarity gives
\begin{equation*}
\operatorname{Im}\langle w,dw\rangle =\operatorname{Im}\langle z,dz\rangle
+\ell ^{-1}\operatorname{Im}\langle z,Xz\rangle \,dt,
\end{equation*}
and
\begin{equation*}
|dw|^{2}=|dz+\ell ^{-1}Xz\,dt|^{2},\qquad |w|=|z|.
\end{equation*}
The Reeb field $\partial _{t}$ of $\mathbb{H}_{2}$ pulls back to $\partial
_{t}-\ell ^{-1}V_{X}$, and $\eta _{X}(\xi _{X})=1$. This is
\eqref{eq:suspension-model}.

(3) The approximation follows from (2). The order of $A_{i}$ is the least
common denominator of $r_{1}^{(i)},r_{2}^{(i)}$, which tends to infinity
because one of $\theta _{1},\theta _{2}$ is irrational.
Proposition~\ref{prop:transflat}(3) gives the periods and the developability.
The maps $(z,t)\mapsto (U\operatorname{diag}(e^{i\phi _{1}},e^{i\phi
_{2}})U^{\ast }z,t+s)$
commute with every $X^{\prime }\in
U\operatorname{diag}(i\mathbb{R},i\mathbb{R})U^{\ast }$ and preserve
\eqref{eq:suspension-model} for $X^{\prime }$.
Finally, $\exp (s\xi _{X_{i}})(z,t)=(e^{-sX_{i}/\ell }z,t+s)$ is the
identity of $N$ exactly when $s=m\ell $ with $A_{i}^{m}=I$. Hence $\xi
_{X_{i}}$ is rational.
\end{proof}

\subsection{Irregular Structure Theorem: Developable Quasi-Regular
Approximation}

We can now relate the irregular case to the structure theorem of the
quasi-regular case. The main irregular structure theorem is the following.

\begin{theorem}
\label{thm:irregular-main}Let $(M^{5},\xi ,\eta ,\Phi ,g,f)$ be a complete
noncompact gradient shrinking Sasaki--Ricci soliton whose Reeb foliation is
not quasi-regular. The following are equivalent.

\begin{enumerate}
\item[(1)] The Reeb foliation is globally developable.

\item[(2)] $M$ admits a developable quasi-regular soliton approximation
(Definition~\ref{def:qr-approx}).

\item[(3)] $R^{T}\equiv 0$.

\item[(4)] $M$ is isomorphic to $\mathbb{H}_{2}$ or to an irrational
Gaussian suspension $M_{A,\ell }$.

\item[(5)] $\pi _{1}(M)$ is infinite, or the Reeb flow of $M$ has no closed
orbit.
\end{enumerate}

\noindent If they hold, then every developable quasi-regular soliton
approximation $(M_i)$ of $M$ has, for all large $i$, finite covers
$\widetilde M_i\to M_i$ whose lifted Reeb structures are regular Gaussian
lifts. Their underlying smooth manifolds satisfy
\begin{equation*}
\widetilde{M}_{i}\cong \mathbb{C}^{2}\times \mathbb{S}^{1}.
\end{equation*}
Their Reeb periods $L_i$, equal to the generic periods of $M_i$, tend to
infinity. In particular no FIK,
round-cylinder or BCCD branch of Theorem~\ref{thm:good} occurs in it. If
they fail, then every quasi-regular soliton approximation of $M$ consists,
for all large $i$, of solitons with nondevelopable leaf spaces. This
alternative occurs, for instance, for the solitons $\mathcal{S}_{a}$,
$a_{1}/a_{2}\notin \mathbb{Q}$, of Example~\ref{ex:counterexample}
(Proposition~\ref{prop:S3-family}(3)). The structure of $M$ in this
alternative is described in Theorem~\ref{thm:nondev-structure}.
\end{theorem}

Before the proof of the main theorem, we show that, unless the Reeb flow is
a proper $\mathbb{R}$-action, it always has closed orbits. For a basic
function $u,$ we write
\begin{equation}
\zeta _{u}=u\,\xi +\tfrac{1}{2}\Phi \nabla ^{T}u  \label{eq:contact-field}
\end{equation}
for the contact vector field with Hamiltonian $u$. It satisfies $\eta (\zeta
_{u})=u$ and $\mathcal{L}_{\zeta _{u}}\eta =0$ (cf.~\cite{BG}).

\begin{proposition}[{\cite[Lemma~\Iref{lem:growth}{3.4}]{PartI}, \cite{BG,MyS}}]
\label{prop:soliton-torus} Let $(M^{2n+1},\xi ,\eta ,\Phi ,g,f)$ be a connected
complete gradient shrinking Sasaki--Ricci soliton whose Reeb closure $H$ is
compact.

\begin{enumerate}
\item $\zeta _{f}$ is a complete Killing field and an infinitesimal
automorphism of $(\xi ,\eta ,\Phi ,g)$. It commutes with $\xi $ and with $H$,
and $\zeta _{f}f=0$.

\item The closure $H_{f}$ of the group generated by the flows of $\xi $ and
$\zeta _{f}$ in $\operatorname{Isom}(M,g)$ is a compact torus of dimension at most
$n+1$, which contains $H$, acts effectively by automorphisms of the Sasakian
structure and preserves $f$. We call $H_{f}$ the \emph{soliton torus}.
\end{enumerate}
\end{proposition}

\begin{proof}
(1) Since $\xi $ preserves $\Phi $ and $f$, $[\xi ,\zeta _{f}]=0$. In the
convention $d\eta =2h(\Phi \cdot ,\cdot ),$ one has
\begin{equation*}
d\eta (\Phi Y,\cdot )=-2h(Y,\cdot )
\end{equation*}
on $D$ and then
\begin{equation*}
\mathcal{L}_{\zeta _{f}}\eta =df+\iota _{\zeta _{f}}d\eta =0.
\end{equation*}
The transverse part $\frac{1}{2}J\nabla ^{T}f$ of $\zeta _{f}$ is
transversely holomorphic, since $\nabla ^{T}f$ is (the Hessian of $f$ is
$J$-invariant by the full real equation \eqref{eq:soliton}), and it preserves
$\omega ^{T}=\frac{1}{2}d\eta .$ Hence it preserves $J$, $\omega ^{T}$ and
$g^{T}$, and $\zeta _{f}$ preserves $\xi ,\eta ,\Phi $ and $g=g^{T}+\eta
^{2}$. Killing fields on complete manifolds are complete. Finally
\begin{equation*}
\zeta _{f}f=\frac{1}{2}g(\Phi \nabla ^{T}f,\nabla ^{T}f)=0,
\end{equation*}
and $\zeta _{f}$ commutes with the closure $H$ of the Reeb flow by
continuity.

(2) By \cite[Lemma~\Iref{lem:growth}{3.4}]{PartI} the potential is proper.
The isometries of $(M,g)$ preserving $f$ form a closed subgroup which maps
the compact set $\{f\leq f(x)\}$ into itself. Since $\operatorname{Isom}(M,g)$
acts properly \cite{MyS}, this subgroup is compact. It contains the flows of
$\xi $ and $\zeta _{f}$, which preserve $f$ by (1), so $H_{f}$ is a compact
connected abelian Lie group, i.e.\ a torus, and it acts by automorphisms as
these form a closed subgroup. For $X,Y$ in the Lie algebra of $H_{f},$ one
has $X(\eta (Y))=\eta ([X,Y])=0$, hence $d\eta (X,Y)=0$. At a point of a
principal orbit, where the infinitesimal action is injective, the
projections of these fields to $D$ therefore span an isotropic subspace of
the symplectic vector space $(D,d\eta )$ of dimension at most $n$ and at
most one direction is mapped into $\mathbb{R}\xi $. Hence $\dim H_{f}\leq
n+1 $.
\end{proof}

\begin{proposition}
\label{prop:min-set}Under the assumptions of
Proposition~\ref{prop:soliton-torus}, let $C$ be a connected component of the
critical set $\operatorname{Crit}(f)$. Then

\begin{enumerate}
\item $f$ is constant on $C$, say $f|_{C}=a$, and $C$ is a connected
component of the zero set of the Killing field $\zeta _{f}-a\xi $.
Consequently, $C$ is a closed, totally geodesic, $\Phi $-invariant
submanifold tangent to $\xi $.

\item Moreover, with the induced structure, $C$ is a Sasakian manifold of
odd dimension. The minimum set $f^{-1}(\min f)$ is nonempty, compact and a
finite union of such components.
\end{enumerate}
\end{proposition}

\begin{proof}
By Sard's theorem the set of critical values of $f$ has measure zero, so $f$
is constant on the connected set $C$. For a constant $a$, $\zeta _{f}-a\xi
=(f-a)\xi +\frac{1}{2}\Phi \nabla ^{T}f$ vanishes exactly at the critical
points of $f$ in $f^{-1}(a)$, because $\Phi $ is injective on $D$. Hence
$C\subset \mathrm{Zero}(\zeta _{f}-a\xi )\subset \operatorname{Crit}(f)$, and $C$
is a connected component of $\mathrm{Zero}(\zeta _{f}-a\xi )$. The zero set
of a Killing field is a disjoint union of closed totally geodesic
submanifolds, and at $x\in C$ one has $T_{x}C=\ker \nabla (\zeta _{f}-a\xi
)(x)$. Since $\zeta _{f}-a\xi $ is an infinitesimal automorphism commuting
with $\xi $, $C$ is invariant under the Reeb flow, and $\nabla (\zeta
_{f}-a\xi )(x)$ is skew, annihilates $\xi (x)$ and commutes with $\Phi _{x}$.
Thus $T_{x}C=\mathbb{R}\xi \oplus (T_{x}C\cap D_{x})$ with a $\Phi $-invariant
second summand.

A $\Phi $-invariant submanifold tangent to $\xi $ is Sasakian with the
induced structure, since its metric cone is a complex submanifold of the
K\"{a}hler cone of $M$. Finally $f$ is proper and bounded below, so
$f^{-1}(\min f)$ is nonempty and compact. It is a finite union of components
of $\mathrm{Zero}(\zeta _{f}-(\min f)\xi ).$
\end{proof}

\begin{corollary}
\label{cor:closed-orbit}Let $(M^{5},\xi ,\eta ,\Phi ,g,f)$ be a complete
noncompact gradient shrinking Sasaki--Ricci soliton. If its Reeb closure $H$
is compact, then every connected component of the minimum set of $f$
contains a closed Reeb orbit. Consequently, if the Reeb flow of $M$ has no
closed orbit, then $M$ is in case $(P)$ and $M\cong \mathbb{H}_{2}$.
\end{corollary}

\begin{proof}
Let $C$ be a component of $f^{-1}(\min f)$. By Proposition~\ref{prop:min-set},
it is a compact Sasakian manifold whose Reeb flow is the restriction of
that of $M$. Let $T_{C}$ be the closure of this flow in
$\operatorname{Isom}(C)$, which is a torus. We argue by induction on $\dim C+\dim T_{C}$.

If $\dim T_{C}=1$, the Reeb flow of $C$ factors through a circle action, and
all its orbits are closed.

If $\dim T_{C}\geq 2$, choose $\zeta $ in the Lie algebra of $T_{C}$ which
is not a multiple of $\xi $. The function $\eta (\zeta )$ attains its
maximum on the compact manifold $C$ at some $x$, and $d(\eta (\zeta
))=-\iota _{\zeta }d\eta $, together with the nondegeneracy of $d\eta $ on
$D\cap TC$, gives $\zeta (x)=\eta (\zeta )(x)\,\xi (x)$. Hence the nonzero
Killing field $\zeta -\eta (\zeta )(x)\xi $ of $C$ vanishes at $x$. As in
the proof of Proposition~\ref{prop:min-set}, the component $C^{\prime }$ of
its zero set through $x$ is a compact Sasakian submanifold of smaller
dimension. It is invariant under the connected group $T_{C}$, which commutes
with $\zeta $, and the closure of its Reeb flow is contained in the image of
$T_{C}$ in $\operatorname{Isom}(C^{\prime })$. This image is a torus of dimension
at most $\dim T_{C}-1$, since the one-parameter group generated by $\zeta
-\eta (\zeta )(x)\xi $ acts trivially on $C^{\prime }$. By induction
$C^{\prime }$ contains a closed Reeb orbit. For the last assertion, if the
Reeb flow has no closed orbit, then $H$ is not compact, so the Reeb flow is
a closed, free and proper $\mathbb{R}$-action by Lemma~\ref{lem:reeb}. Hence
Corollary~\ref{cor:Rcase} gives $M\cong \mathbb{H}_{2}$.
\end{proof}

We can now prove the main irregular structure theorem.

\begin{proof}[Proof of Theorem~\protect\ref{thm:irregular-main}]
(i) (1)$\Rightarrow $(4): By Theorem~\ref{thm:globdev}. (4)$\Rightarrow $(1):
The universal cover is $\mathbb{H}_{2}\rightarrow \mathbb{C}^{2}$.

(ii) (4)$\Rightarrow $(3): The Heisenberg model is transversely flat.
(3)$\Rightarrow $(4): by Proposition~\ref{prop:transflat}, since $M$ is not
quasi-regular.

(iii) \ (4)$\Rightarrow $(2): By Proposition~\ref{prop:suspension-approx}(1)
and (3). \ (2)$\Rightarrow $(3): By Corollary~\ref{cor:nongaussian}(2).

(iv) (4)$\Rightarrow $(5): $\pi _{1}(M_{A,\ell })\cong \mathbb{Z}$, and the
Reeb leaves of $\mathbb{H}_{2}$ are lines. (5)$\Rightarrow $(4): if $\pi
_{1}(M)$ is infinite, Theorem~\ref{thm:pi1} gives $M\cong M_{A,\ell }$, and
$A$ has infinite order by Proposition~\ref{prop:transflat}(2),(3), since $M$
is not quasi-regular. If the Reeb flow has no closed orbit, Corollary
\ref{cor:closed-orbit} gives $M\cong \mathbb{H}_{2}$.

(v) The statement about approximations is Corollary~\ref{cor:nongaussian}(2).
Finally, if (1)--(5) fail and a quasi-regular soliton approximation had
infinitely many developable members, these would form a developable
quasi-regular soliton approximation, and (2) would hold.
\end{proof}

\subsection{Nondevelopable Quasi-Regular Soliton Approximation}

We first describe the solitons for which (1)--(5) of
Theorem~\ref{thm:irregular-main} fail, and then show by explicit examples that
this case
is not empty.

\begin{theorem}
\label{thm:nondev-structure}Let $(M^{5},\xi ,\eta ,\Phi ,g,f)$ be a complete
noncompact gradient shrinking Sasaki--Ricci soliton whose Reeb foliation is
neither quasi-regular nor globally developable. Then the following hold.

\begin{enumerate}
\item $M$ is in case $(T)$, with $H\cong T^{k}$ and $k\in \{2,3\}$, and
\begin{equation*}
R^{T}>0
\end{equation*}
on $M$.

\item $\pi _{1}(M)$ is finite. The universal covering $\widehat{M}\rightarrow
M$ is finite, $\widehat{M}$ is in case $(T)$ with the same $k$,
its Reeb foliation is not globally developable, and $M=\widehat{M}/\Gamma $
for a finite group $\Gamma \cong \pi _{1}(M)$ of automorphisms of the lifted
soliton which acts freely and commutes with $\widehat{H}$.

\item If $(M_{i})$ is a quasi-regular soliton approximation of $M$, then for
all large $i,$ the leaf space of $M_{i}$ is nondevelopable, and $M_{i}$ has
the properties of Corollary~\ref{cor:nondev-qr}: $R^{T}>0$ on $M_{i}$, $\pi
_{1}(M_{i})$ is finite, and the universal cover of $M_{i}$ is quasi-regular
with a nondevelopable leaf space. Moreover the generic Reeb periods of $M_{i}
$ tend to infinity.
\end{enumerate}
\end{theorem}

\begin{proof}
(1) Condition (1) of Theorem~\ref{thm:irregular-main} fails, hence so do (3)
and (4). Thus $R^{T}\not\equiv 0$, and $R^{T}>0$ by Lemma~\ref{lem:nonneg}.
In case $(P)$, Corollary~\ref{cor:Rcase} would give
\begin{equation*}
M\cong \mathbb{H}_{2},
\end{equation*}
which is (4); hence $M$ is in case $(T)$, and $k\leq 3$ by
Proposition~\ref{prop:second-structure}.

(2) Since $R^{T}\not\equiv 0$, alternative (2) of Theorem~\ref{thm:pi1}
holds. It gives the finiteness of $\Gamma $, its commutation with
$\widehat{H}$, $\dim \widehat{H}=k$, and that $\widehat{M}$ is not
quasi-regular. Hence
$\widehat{M}$ is in case $(T)$. Its generic Reeb leaves are not closed, so
they are not the fibers of a submersion. Since $\widehat{M}$ is simply
connected, its Reeb foliation is not globally developable. Deck
transformations preserve the lifted structure and potential.

(3) By Theorem~\ref{thm:irregular-main}, for large $i$ the leaf spaces are
not developable, and Corollary~\ref{cor:nondev-qr} applies. Since the
principal Reeb leaves of $M$ are not compact, Lemma~\ref{lem:period-blowup}(1)
gives $L_{i}\rightarrow \infty $.
\end{proof}

Finally, the solitons of Example~\ref{ex:counterexample} provide
nondevelopable quasi-regular soliton approximations.

\begin{proposition}
\label{prop:S3-family}Fix $\lambda>0$ and let
$a=(a_1,a_2)$ range over the positive pairs satisfying
$a_1+a_2=\lambda/2$.

\begin{enumerate}
\item For every $a$, $\mathcal{S}_{a}$ is a complete gradient shrinking
Sasaki--Ricci soliton, normalized by \eqref{eq:pre-normal} with bounded
curvature. Moreover, $\mathcal{S}_{a}$ depends real-analytically on $a$ in
$C^{\infty }$ on compact subsets of
\begin{equation*}
\mathbb{C}\times \mathbb{S}^{3}.
\end{equation*}

\item If $a_{1}/a_{2}=m_{1}/m_{2}$ for positive integers with
$(m_{1},m_{2})=1$, then $\mathcal{S}_{a}$ is quasi-regular with generic Reeb
period $L=4\pi (m_{1}+m_{2})/\lambda $ and exceptional periods $L/m_{1}$ and
$L/m_{2}$ over $t=0$ and $t=1$. Its leaf space is
\begin{equation*}
\mathbb{C}\times \mathbb{CP}_{[m_{1},m_{2}]}^{1},
\end{equation*}
as in Remark~\ref{rem:weighted-example}. It is nondevelopable unless
$m_{1}=m_{2}=1$. For $m_{1}=m_{2}=1$, $\mathcal{S}_{a}$ is the regular
cylinder lift of Remark~\ref{rem:cylmodel} with $d=1$.

\item If $a_{1}/a_{2}\notin \mathbb{Q}$, then $\mathcal{S}_{a}$ is
irregular, in case $(T)$ with $H=T^{2}$ and bounded contact moment map
$\mu$ from \eqref{eq:contact-moment}. Choose positive pairs
$a^{(i)}\to a$ with
$a_1^{(i)}+a_2^{(i)}=\lambda/2$ and
$a_1^{(i)}/a_2^{(i)}\in\mathbb Q$. Then the solitons
$\mathcal{S}_{a^{(i)}}$ form a quasi-regular soliton approximation of $\mathcal{S}_{a}$ on the fixed
manifold $\mathbb{C}\times \mathbb{S}^{3}$. For large $i$ its members are
nondevelopable, and their isotropy orders $m_{1}^{(i)},m_{2}^{(i)}$ and
generic periods tend to infinity. Its Reeb fields are rational elements of
$\mathfrak{h}$ converging to $\xi _{a}$. The limit has positive nonconstant
$R^{T}$.

\item As $a\rightarrow (\lambda /4,\lambda /4)$, also along irrational
ratios, $\mathcal{S}_{a}$ converges on the fixed manifold to the regular
cylinder lift with $d=1$.

\item Let $a_{1}>a_{2}$. The critical set of $f$ is $\{0\}\times (\gamma
_{-}\cup \gamma _{+})$, where $\gamma _{-}=\{t=0\}$ and $\gamma _{+}=\{t=1\}$
are closed Reeb orbits of lengths $2\pi /a_{1}$ and $2\pi /a_{2}$, at the
levels $\psi _{-}=c/(e^{c}-1)$ and $\psi _{+}=ce^{c}/(e^{c}-1)$ of $f$.
$\{0\}\times \gamma _{-}$ is the minimum set. The soliton torus $H_{f}$ of
Proposition~\ref{prop:soliton-torus} is the three-torus generated by $\partial
_{\theta _{1}}$, $\partial _{\theta _{2}}$ and the rotations of $\mathbb{C}$.

\item $\max_{M}R^{T}/\min_{M}R^{T}=e^{|c|}$. Consequently
$\mathcal{S}_{a}\cong \mathcal{S}_{a^{\prime }}$ if and only if $a^{\prime
}=a$ or $a^{\prime }=(a_{2},a_{1})$.
\end{enumerate}
\end{proposition}

\begin{proof}
(1) By Example~\ref{ex:S3} and Appendix~\ref{app:A}, $\mathcal{S}_{a}$ is a
complete normalized soliton for every $a$; for $a_{1}=a_{2}$, that is $c=0$,
the construction uses $V=\frac{\lambda }{2}t(1-t)$ and gives $K\equiv
\lambda $. Since $\mathcal{H}^{\prime }(c)<0$ for $c\neq 0$ and
$\mathcal{H}(c)=\frac{1}{2}-\frac{c}{12}+O(c^{3})$, the implicit function
theorem shows
that $c(a)$ is real-analytic. The function $V$ of \eqref{eq:V} is
real-analytic in $(t,c)$, its singularity at $c=0$ being removable, and
$V(t)/(t(1-t))$ is real-analytic and positive on $[0,1]$. In the coordinates
$w_{1}=\sqrt{1-t}\,e^{i\theta _{1}}$, $w_{2}=\sqrt{t}\,e^{i\theta _{2}}$ of
$\mathbb{S}^{3}\subset \mathbb{C}^{2},$ one has $\eta
_{S}=a_{1}^{-1}\operatorname{Im}(\overline{w}_{1}dw_{1})+a_{2}^{-1}\operatorname{Im}(\overline{w}_{2}dw_{2})$.
The
endpoint expansions in Appendix~\ref{app:A} show that the coefficients of
$h_{S}^{T}$, $\xi _{S}$ and $f$ in these coordinates are real-analytic in
$(w,a)$. This proves the analytic dependence. The transverse curvature is
that of the product of the flat plane with the surface of curvature
\begin{equation*}
K=\lambda ce^{ct}/(e^{c}-1),
\end{equation*}
which is bounded.

(2) Write $a_{j}=\kappa m_{j}$ with $\kappa =\lambda /(2(m_{1}+m_{2}))$. The
Reeb flow $(\theta _{1},\theta _{2})\mapsto (\theta _{1}+a_{1}s,\theta
_{2}+a_{2}s)$ is periodic with generic period $2\pi /\kappa =L$. The
endpoint circles, on which one angle collapses, have periods $2\pi
/a_{1}=L/m_{1}$ and $2\pi /a_{2}=L/m_{2}$. The quotient of $\mathbb{S}^{3}$
by this circle action is the weighted projective line
$\mathbb{CP}_{[m_{1},m_{2}]}^{1}$ of Remark~\ref{rem:weighted-example}. It is
a bad
orbifold if $(m_{1},m_{2})\neq (1,1)$, so the leaf space is not developable.
If $m_{1}=m_{2}=1$, then $c=0$, the transverse metric is flat times round of
curvature $\lambda $, the potential is $1+\frac{\lambda }{2}|z|^{2}$, and
the Reeb flow is the free Hopf action of period $8\pi /\lambda $. This is
the regular lift with $d=1$, unique up to gauge by
Proposition~\ref{prop:lifting}.

(3) Such a sequence exists: choose positive rational numbers $r_i\to a_1/a_2$
and set
$a_1^{(i)}=\lambda r_i/(2(1+r_i))$ and
$a_2^{(i)}=\lambda/(2(1+r_i))$.
The approximation and convergence of the Reeb fields then follow from
(1) at the fixed value of $\lambda$. The orders $m_{j}^{(i)}$ tend to infinity because
$m_{1}^{(i)}/m_{2}^{(i)}\rightarrow a_{1}/a_{2}\notin \mathbb{Q}$; in
particular $(m_{1}^{(i)},m_{2}^{(i)})\neq (1,1)$ for large $i$, so these
leaf spaces are not developable by (2). The closure of the flow of $\xi _{a}$
is the torus $T^{2}$ acting on the angles, $\mu $ depends only on $t\in
\lbrack 0,1]$, and \eqref{eq:ctr-RT} gives $R^{T}=2K>0$, which is
nonconstant since $c\neq 0$.

(4) This is (1) at $a=(\lambda /4,\lambda /4)$ together with (2).

(5) On $\mathcal{S}_{a}=\mathbb{C}\times _{\mathrm{H}}S_{a}$, one has
$f=\frac{\lambda }{2}|z|^{2}+\psi $ with $\psi =ct+c/(e^{c}-1)$, and $\nabla
^{T}f$ is the sum of $\lambda $ times the horizontal lift of the position
vector of $\mathbb{C}$ and of $\nabla ^{T}\psi $, where $|\nabla ^{T}\psi
|^{2}=2c^{2}V(t)$.

Since $V>0$ on $(0,1)$ and $V(0)=V(1)=0$, the critical set is $\{z=0,\ t\in
\{0,1\}\}$. On $\{t=0\}$ the $\theta _{2}$-circle collapses, and on $\{t=1\}$ the
$\theta _{1}$-circle collapses. Hence the Reeb flow
\begin{equation*}
(\theta _{1},\theta _{2})\mapsto (\theta _{1}+a_{1}s,\theta _{2}+a_{2}s)
\end{equation*}
has periods $2\pi /a_{1}$ and $2\pi /a_{2}$ there and these are the lengths
since $|\xi |=1$. As $c>0$, the minimum of $\psi $ is attained at $t=0$.

To determine the closure defining $H_f$, put $b=c/(e^c-1)$ and let
$\partial_\varphi$ generate the rotations of $\mathbb C$, with period
$2\pi$. Formula~\eqref{eq:contact-field} gives
\[
\xi=a_1\partial_{\theta_1}+a_2\partial_{\theta_2},\qquad
\zeta_f=a_1b\partial_{\theta_1}+a_2(b+c)\partial_{\theta_2}
+\frac\lambda2\partial_\varphi.
\]
A character of the ambient three-torus, represented by
$(k_1,k_2,k_3)\in\mathbb Z^3$, vanishes on both flows precisely when
\[
k_1a_1+k_2a_2=0,\qquad k_2a_2c+k_3\lambda/2=0.
\]
If $a_1/a_2$ is irrational, these equations force all $k_j$ to vanish.
If $a_1/a_2$ is rational, put $r=a_2/(a_1+a_2)\in\mathbb Q$. Then
$c$ is irrational: otherwise \eqref{eq:Hc} would give
$e^c=1+c/(1-rc)\in\mathbb Q$ for the nonzero rational number $c$,
contradicting the transcendence of $e$: if $c=m/n$ with
nonzero integers $m,n$, then $e^m=(e^c)^n$ would be algebraic,
and hence so would $e$.
Here $1-rc=c/(e^c-1)>0$, so the denominator does not vanish.
Dividing the second character equation by $\lambda/2=a_1+a_2$ now gives
$k_2rc+k_3=0$, hence $k_2=k_3=0$ and then $k_1=0$.
There is therefore no nonzero integer character annihilating the two
flows, so their closure is the whole three-torus. This proves
$H_f=U(1)\times T^2$.

(6) By \eqref{eq:ctr-RT}, $R^{T}=2\lambda ce^{ct}/(e^{c}-1)$ with $t\in
\lbrack 0,1]$, so
\begin{equation*}
\max R^{T}/\min R^{T}=e^{|c|}
\end{equation*}
is an isomorphism invariant. It determines $|c|$ and hence, by \eqref{eq:Hc}
and $a_{1}+a_{2}=\lambda /2$, the unordered pair $\{a_{1},a_{2}\}$.
Conversely, exchanging $w_{1}$ and
$w_{2}$ (and $t$ with $1-t$) identifies $\mathcal{S}_{(a_{1},a_{2})}$ with
$\mathcal{S}_{(a_{2},a_{1})}$.
\end{proof}

\section{Naber-Type Global Splitting Theorem}

\label{sec:split}

The solitons $\mathcal{S}_{a}$ of Proposition~\ref{prop:S3-family} split
globally. That is, their transverse K\"{a}hler structure is the product of
the Gaussian plane and the transverse structure of a Sasakian three-sphere.
By Theorem~\ref{thm:splitting} every developed local limit at transverse
infinity splits off a flat complex line. It is then natural to ask which
solitons split globally, and whether they are always as in
Proposition~\ref{prop:S3-family}, possibly together with positivity of the
curvature of the
remaining factor. We show that, up to finite quotients and the transversely
flat models, the solitons $\mathcal{S}_{a}$ are the only globally split
shrinkers, and that they are not rigid unless $a_{1}=a_{2}$.

Finally, we show that nonnegativity of the transverse holomorphic
bisectional curvature forces global splitting (Theorem~\ref{thm:nonneg-op}).
Hence it characterizes these solitons among the nonflat ones.

\subsection{Globally Split Shrinkers}

\begin{definition}
\label{def:split}(1) Let $N^{3}$ be a Sasakian three-manifold with structure
tensors $\xi _{N},\eta _{N},\Phi _{N},g_{N}$ and transverse metric $h_{N}$,
and let $\psi _{N}$ be a basic function on $N$. The \emph{Heisenberg product}
\begin{equation*}
\mathbb{C}\times _{\mathrm{H}}N
\end{equation*}
is $\mathbb{C}\times N$ with $z=x+iy$, $\xi =\xi _{N},$ $f=\tfrac{\lambda
}{2}|z|^{2}+\psi _{N}$ and
\begin{gather}
\eta =\eta _{N}+x\,dy-y\,dx,\qquad  \notag \\
g=dx^{2}+dy^{2}+h_{N}+\eta ^{2},\qquad  \label{eq:H-product}
\end{gather}
and with the product transverse complex structure. As in
Example~\ref{ex:counterexample} (see Appendix~\ref{app:A}), it is Sasakian,
and it is a
normalized gradient shrinking Sasaki--Ricci soliton with constant $\lambda $
if and only if $(N,\psi _{N})$ is. It is complete if $N$ is complete.

(2) A five-dimensional gradient shrinking Sasaki--Ricci soliton is
\emph{split} if there is a globally defined smooth real section $V$ of
$D$ with $|V|=1$ and
\begin{equation*}
\nabla ^{T}V=0,
\end{equation*}
where the equality includes every direction in $TM$, with
$\nabla_\xi^T V=[\xi,V]$. It is \emph{globally split} if the lifted
soliton on its universal cover is split. Thus ``split'' requires the
parallel section on $M$ itself, whereas ``globally split'' allows a
nontrivial action of the deck group on the parallel complex line upstairs.
These terms concern the transverse connection and the Heisenberg product
\eqref{eq:H-product}; they do not assert a Riemannian product decomposition
of the full Sasakian metric.
\end{definition}

\begin{example}
\begin{enumerate}
\item Let $S_{a}=(\mathbb{S}^{3},\xi _{a},\eta _{S},\Phi _{S},g_{S},\psi
_{a})$, $a_{1},a_{2}>0$, $a_{1}+a_{2}=\lambda /2$, be the three-dimensional
soliton of Example~\ref{ex:S3}, with the normalized potential $\psi
_{a}=ct+c/(e^{c}-1)$ as in Proposition~\ref{prop:S3-family}. Thus
\begin{equation*}
\mathcal{S}_{a}=\mathbb{C}\times _{\mathrm{H}}S_{a}.
\end{equation*}

\item The Heisenberg model is
\begin{equation*}
\mathbb{H}_{2}=\mathbb{C}\times _{\mathrm{H}}\mathbb{H}_{1},
\end{equation*}
where $\mathbb{H}_{1}$ is the three-dimensional Heisenberg group with $\eta
=dt+x\,dy-y\,dx$ and potential $\frac{\lambda }{2}|z|^{2}$.

\item The cylindrical lift of Remark~\ref{rem:cylmodel} is
\begin{equation*}
\mathbb{C}\times _{\mathrm{H}}L(d,1).
\end{equation*}

\item Every Heisenberg product is split, with $V$ the horizontal lift
$\partial _{x}+y\xi $ of $\partial _{x}$. Every Gaussian suspension $M_{A,\ell
}$ is globally split, since its universal cover is $\mathbb{H}_{2}$.
\end{enumerate}
\end{example}

\begin{theorem}
\label{thm:split-structure}Let $(M^{5},\xi ,\eta ,\Phi ,g,f)$ be a connected
complete gradient shrinking Sasaki--Ricci soliton with constant $\lambda $,
normalized by \eqref{eq:pre-normal}.

\begin{enumerate}
\item Let $M$ be split and $V,u,\psi ,X_{1},X_{2}$ be as in
Lemma~\ref{lem:split-affine}. Then $N:=u^{-1}(0)$ is a connected, closed,
totally
geodesic and $\Phi $-invariant three-dimensional submanifold tangent to $\xi
$. With the induced structure and $\psi _{N}:=f|_{N}=\psi |_{N},$ it is a
complete gradient shrinking Sasaki--Ricci soliton with constant $\lambda $,
normalized by \eqref{eq:pre-normal}, and
\begin{equation}
\Theta :\mathbb{C}\times _{\mathrm{H}}N\longrightarrow M,\qquad \Theta
(x+iy,q)=\phi _{x}^{1}\circ \phi _{y}^{2}\circ \varphi _{xy}(q),
\label{eq:split-iso}
\end{equation}
where $\phi ^{j}$ is the flow of $X_{j}$ and $\varphi $ the Reeb flow, is an
isomorphism of solitons with $u\circ \Theta (z,q)=z$.

\item If $M$ is globally split, then its universal cover is isomorphic to
$\mathbb{C}\times _{\mathrm{H}}N$ for a simply connected complete
three-dimensional soliton $(N,\psi _{N})$ as in (1) such that
\begin{equation*}
M\cong (\mathbb{C}\times _{\mathrm{H}}N)/\Gamma
\end{equation*}
with $\Gamma \cong \pi _{1}(M)$ acting freely by automorphisms.

\item Let $(N,\psi _{N})$ be a complete three-dimensional soliton with
$R^{T}\not\equiv 0$. Then the transversely parallel sections of $D$ on
$\mathbb{C}\times _{\mathrm{H}}N$ are the constant linear combinations of $V$
and $\Phi V$ and
\begin{equation*}
\operatorname{Aut}(\mathbb{C}\times _{\mathrm{H}}N,f)=U(1)\times
\operatorname{Aut}(N,\psi
_{N}),
\end{equation*}
with $U(1)$ acting by rotations of $\mathbb{C}$.

\item Under the nonflat-factor hypothesis of (3), let $\Gamma$ be a
 discrete group of automorphisms preserving the normalized potential and
acting freely and properly discontinuously on
$\mathbb C\times_{\mathrm H}N$. Its action has the form $\gamma (z,q)=(\chi (\gamma
)z,\gamma _{N}(q))$, where $\gamma \mapsto \gamma _{N}$ is injective with
image $\Gamma _{N}$ acting freely on $N$ and $\chi :\Gamma \rightarrow U(1)$
is a homomorphism. The quotient is split if and only if $\chi $ is trivial,
and then it is $\mathbb{C}\times _{\mathrm{H}}(N/\Gamma _{N})$.
\end{enumerate}
\end{theorem}

\begin{proof}
See Appendix~\ref{app:B}.
\end{proof}

\subsection{Nonnegative Transverse Holomorphic Bisectional Curvature}

\label{sec:nonneg-bisec}

The remaining factor $N$ is classified as follows. This is where the
positivity of its curvature and the critical points of its potential enter.

\begin{proposition}
\label{prop:3d-shrinkers}Let $(N^{3},\xi ,\eta ,\Phi ,g,\psi )$ be a
connected complete gradient shrinking Sasaki--Ricci soliton normalized by
$R^{T}+|\nabla ^{T}\psi |^{2}=2\lambda \psi $. Let $K=R^{T}/2$ be its
transverse Gaussian curvature. Then there is a constant $\kappa \in \lbrack
0,\lambda /e]$ such that
\begin{equation}
K=\kappa e^{\psi }\qquad \text{and}\qquad |\nabla ^{T}\psi |^{2}=2\lambda \psi -2\kappa e^{\psi }.  \label{eq:3d-identities}
\end{equation}

\begin{enumerate}
\item If $\kappa =0$, then $N$ is transversely flat and noncompact. In
addition, it is isomorphic to $\mathbb{H}_{1}$ if it is simply connected.

\item If $\kappa >0$, then $N$ is compact, $\pi _{1}(N)$ is finite, and the
universal cover of $N$ is isomorphic to $S_{a}$ for a unique $a$ with
$a_{1}\geq a_{2}>0$ and $a_{1}+a_{2}=\lambda /2$.

\begin{enumerate}
\item If $\kappa =\lambda /e$, then $\psi \equiv 1$, $K\equiv \lambda $ and
$a_{1}=a_{2}$.

\item If $\kappa <\lambda /e$, let $\psi _{-}<1<\psi _{+}$ be the solutions
of $\lambda x=\kappa e^{x}$. Then $\psi (N)=[\psi _{-},\psi _{+}]$, $\lambda
\psi _{-}\leq K\leq \lambda \psi _{+}$ and the critical set of $\psi $ is
$\psi ^{-1}(\psi _{-})\cup \psi ^{-1}(\psi _{+})$. Each of these two sets is
a single closed Reeb orbit and $\psi _{+}-\psi _{-}$ is the constant $c>0$
attached to $a$ by \eqref{eq:Hc}, in particular $a_{1}>a_{2}$.
\end{enumerate}

\item $\operatorname{Aut}(S_{a},\psi _{a})$ is the angle torus $T^{2}$ if
$a_{1}\neq
a_{2}$, and $U(2)$ if $a_{1}=a_{2}$.
\end{enumerate}
\end{proposition}

\begin{proof}
See Appendix~\ref{app:B}.
\end{proof}

We can now describe the globally split solitons.

\begin{theorem}
\label{thm:split-class}Let $(M^{5},\xi ,\eta ,\Phi ,g,f)$ be a connected
complete gradient shrinking Sasaki--Ricci soliton with constant $\lambda $,
normalized by \eqref{eq:pre-normal}, which is globally split.

\begin{enumerate}
\item If $R^{T}\equiv 0$, then $M$ is isomorphic to $\mathbb{H}_{2}$ or to a
Gaussian suspension $M_{A,\ell }$. All of these are globally split.

\item Otherwise $R^{T}>0$, and $M\cong \mathcal{S}_{a}/\Gamma $ for a unique
$a$ with $a_{1}\geq a_{2}>0$, $a_{1}+a_{2}=\lambda /2$, and a finite group
$\Gamma $ of automorphisms of $\mathcal{S}_{a}$ acting freely, unique up to
conjugation.

\begin{enumerate}
\item If $a_{1}>a_{2}$, then $\Gamma $ is cyclic of some order $p\geq 1$,
generated by
\begin{equation*}
(z,w_{1},w_{2})\longmapsto (e^{2\pi ij/p}z,e^{2\pi i/p}w_{1},e^{2\pi
iq/p}w_{2}),\qquad \gcd (q,p)=1,\ j\in \mathbb{Z},
\end{equation*}
in the coordinates $w_{1}=\sqrt{1-t}\,e^{i\theta _{1}}$,
$w_{2}=\sqrt{t}\,e^{i\theta _{2}}$ of $\mathbb{S}^{3}$. Here $q$ and $j$
are residue classes modulo $p$, with $q$ invertible modulo $p$, and
\[
L(p,q):=\mathbb S^3/\bigl\langle
(w_1,w_2)\mapsto(e^{2\pi i/p}w_1,e^{2\pi iq/p}w_2)\bigr\rangle.
\]
The symbol $q$ here is the lens-space action weight, not the Euler number
used for circle bundles in Theorem~\ref{thm:intro-topology}; $p=1$ denotes
the trivial group and $L(1,q)=\mathbb S^3$.
The projection $[z,w]\mapsto[w]$ makes $M$ the complex line bundle defined
by the character $\gamma\mapsto e^{2\pi i j/p}$ of this cyclic action.
It is flat as a line bundle, meaning that it carries a flat unitary
connection; the Sasakian metric on its total space is not flat.
The soliton is split if and only if $p\mid j$, in which case
$M=\mathbb{C}\times_{\mathrm H}(S_a/\mathbb Z_p)$.

\item If $a_{1}=a_{2}$, then $M$ is a finite quotient of the regular
cylinder lift with $d=1$.
\end{enumerate}

\item In case (2) let $\psi _{\pm }$ be as in
Proposition~\ref{prop:3d-shrinkers}, so that $\psi _{-}=|c|/(e^{|c|}-1)$ and
$\psi
_{+}=|c|e^{|c|}/(e^{|c|}-1)$, with $\psi _{\pm }=1$ if $a_{1}=a_{2}$. The
transverse curvature operator is nonnegative, the kernel of $\mathrm{Ric}^{T}
$ is at every point the flat factor, and $2\lambda \psi _{-}\leq
R^{T}\leq 2\lambda \psi _{+}$.

\begin{enumerate}
\item If $a_{1}>a_{2}$, the critical set of $f$ consists of two closed Reeb
orbits: the minimum set, at the level $\psi _{-}$ of length $2\pi /(pa_{1})$,
and an orbit of length $2\pi /(pa_{2})$ at the level $\psi _{+}$.

\item If $a_{1}=a_{2}$, the critical set of $f$ is the three-dimensional
minimum set $f^{-1}(1)$.
\end{enumerate}

\item In case (2) the following are equivalent: (a) $a_{1}=a_{2}$; (b) $R^{T}
$ is constant; (c) $\min f=1$; (d) every critical point of $f$ is a
minimum point; (e) $R^{T}=2\lambda $ at some critical point of $f$; (f) the
Reeb foliation is quasi-regular with developable leaf space; (g) $M$ is
\emph{rigid}. That is, its universal cover is a Heisenberg product
\begin{equation*}
\mathbb{C}\times _{\mathrm{H}}N
\end{equation*}
with $\mathrm{Ric}^{T_{N}}=\lambda h_{N}$ and constant potential on $N$. If
they fail, the Reeb foliation is irregular, and then
Theorem~\ref{thm:nondev-structure} applies, or quasi-regular with
nondevelopable leaf
space.

\item A complete five-dimensional gradient shrinking Sasaki--Ricci soliton
whose transverse Ricci tensor is positive definite at some point is not
globally split.
\end{enumerate}
\end{theorem}

\begin{proof}
See Appendix~\ref{app:B}.
\end{proof}

Theorem~\ref{thm:split-class} classifies complete globally split solitons.
For the local surface factors in Theorems~\ref{thm:splitting} and
\ref{thm:irregular-local}, any use of a complete-factor classification
requires the complete smooth soliton extension stated there; a local
transverse product alone does not provide such an extension.

We now characterize the solitons of Theorem~\ref{thm:split-class} by
nonnegative transverse holomorphic bisectional curvature. A globally split
shrinker has nonnegative transverse curvature operator by
Theorem~\ref{thm:split-class}(3), hence nonnegative transverse holomorphic
bisectional
curvature. We prove the converse without any curvature bound. The main step
is the following compactness criterion. For complete gradient shrinking
K\"{a}hler--Ricci solitons it was obtained by Ni~\cite{Ni} and, by the
least-eigenvalue method which Munteanu--Wang~\cite{MWpos} introduced for
shrinkers with nonnegative sectional curvature, by
Wu--Zhang~\cite[Theorem~1.2]{WZh}. We follow \cite{MWpos,WZh}. Their volume
input, the
infinite volume of noncompact shrinkers, cannot be applied to $(M,g)$, which
is not a gradient Ricci soliton and has
$\mathrm{Ric}_{g}|_{D}=\mathrm{Ric}^{T}-2g^{T}$. It is replaced by the
Calabi--Yau type estimate \cite{Cal,Yau}
for tubes about compact Reeb-closure orbits in Lemma~\ref{lem:transverse-CY},
which uses only $\mathrm{Ric}^{T}\geq 0$. The second ingredient of
Theorem~\ref{thm:nonneg-op} is a direct elliptic argument for the kernel of
the
transverse Ricci tensor.

\begin{proposition}
\label{prop:pos-ricci}A connected complete gradient shrinking Sasaki--Ricci
soliton $(M^{5},\xi ,\eta ,\Phi ,g,f)$ with nonnegative transverse
holomorphic bisectional curvature and positive definite transverse Ricci
tensor is compact.
\end{proposition}

The proof is given at the end of this subsection. Its volume input is proved
in Lemma~\ref{lem:transverse-CY} of Appendix~\ref{app:B}.

\begin{theorem}
\label{thm:nonneg-op}Let $(M^{5},\xi ,\eta ,\Phi ,g,f)$ be a connected
complete noncompact gradient shrinking Sasaki--Ricci soliton with constant
$\lambda >0$ and nonnegative transverse holomorphic bisectional curvature.
Then

\begin{enumerate}
\item $M$ is globally split.

\item After normalizing $f$ by \eqref{eq:pre-normal}, $M$ is isomorphic to
$\mathbb{H}_{2}$, to a Gaussian suspension $M_{A,\ell }$, or to a quotient
$\mathcal{S}_{a}/\Gamma $ as in Theorem~\ref{thm:split-class}(2).

\item In particular, its transverse curvature operator $\mathcal{R}^{T}$ is
nonnegative. In the last case $\mathcal{R}^{T}$ has rank one at every point,
$\mathrm{Ric}^{T}$ has real rank two everywhere, and
\begin{equation*}
2\lambda \psi _{-}\leq R^{T}\leq 2\lambda \psi _{+}.
\end{equation*}
\end{enumerate}
\end{theorem}

\begin{proof}
Normalize $f$ by \eqref{eq:pre-normal} and write
\begin{equation*}
Q=\mathrm{Ric}^{T},\quad s=R^{T},\quad J=\Phi |_{D},
\end{equation*}
with the transverse complex structure fixed in the preliminaries.
The bisectional curvature assumption gives $Q\geq0$. By Lemma~\ref{lem:nonneg}, either $s\equiv 0$ or
$s>0$ everywhere. If $s\equiv 0$, then $Q=0$.

Let $v$ be a unit vector and complete it to a unitary real frame $(v,Jv,e,Je)
$. The identity
\begin{equation*}
0=Q(v,v)=\mathcal{B}^{T}(v,v)+\mathcal{B}^{T}(v,e)
\end{equation*}
and nonnegativity imply that every holomorphic sectional curvature
$\mathcal{B}^{T}(v,v)$ vanishes. Polarization of the K\"{a}hler curvature
tensor gives
$\mathrm{Rm}^{T}=0$. Proposition~\ref{prop:transflat} now gives the
Heisenberg model or a Gaussian suspension, which are globally split.

In the rest of the proof we assume $s>0$.

We first record three identities. By
\cite[Proposition~\Iref{prop:tensor}{3.1}]{PartI},
\begin{equation}
(L_{f}Q)_{il}=2\lambda Q_{il}-2\,T_{iplq}Q_{pq},\qquad T=\mathrm{Rm}^{T},
\label{eq:ricci-drift-hbc}
\end{equation}
where $L_{f}$ is the transverse rough drift Laplacian. Let $x\in M$ and let
$(e_{1},Je_{1},e_{2},Je_{2})$ be a unitary frame of $D_{x}$ diagonalizing the
$J$-invariant tensor $Q$, with eigenvalues $\mu _{1}\leq \mu _{2}$, each of
real multiplicity two. For a unit vector $V\in D_{x}$, the K\"{a}hler
identity $T(V,e,V,e)+T(V,Je,V,Je)=\mathcal{B}^{T}(V,e)$ gives the
paired-frame contraction, writing
$(\varepsilon_1,\varepsilon_2,\varepsilon_3,\varepsilon_4)
=(e_1,Je_1,e_2,Je_2)$ for the full real frame,
\begin{equation}
(T\ast Q)(V,V):=\sum_{\alpha,\beta=1}^{4}
T(V,\varepsilon_\alpha,V,\varepsilon_\beta)
Q(\varepsilon_\alpha,\varepsilon_\beta)
=\mu_1\mathcal B^T(V,e_1)+\mu_2\mathcal B^T(V,e_2),
\label{eq:hbc-contraction}
\end{equation}
where $\mathcal{B}^{T}(V,e)=\mathrm{Rm}^{T}(V,JV,e,Je)$ is the transverse
holomorphic bisectional curvature of unit vectors. In particular $(T\ast
Q)(V,V)\geq 0$ when $\mathcal{B}^{T}\geq 0$ and $Q\geq 0$. Let $V$ be a unit
eigenvector of $Q$ at $x$ for the least eigenvalue $\mu (x)=\mu _{1}$,
extended as a unit field on a local leaf space with $\nabla V(x)=0$. The
unit-length condition gives $\langle\Delta V(x),V(x)\rangle=0$; together
with $Q(V,\cdot)=\mu\langle V,\cdot\rangle$ at $x$, this eliminates the
second-derivative terms involving $V$. Thus $q=Q(V,V)$ satisfies $q\geq\mu$
near $x$, $q(x)=\mu(x)$, and by \eqref{eq:ricci-drift-hbc} and
\eqref{eq:hbc-contraction}
\begin{equation}
L_{f}q(x)=(L_{f}Q)(V,V)(x)=2\lambda \mu (x)-2(T\ast Q)(V,V)(x)\leq 2\lambda
\mu (x).  \label{eq:least-ricci-barrier}
\end{equation}

(i) Vanishing of the least eigenvalue. Let $\mu $ be the least eigenvalue of
$Q$ with respect to $g^{T}$. It is continuous and basic.
Equations~\eqref{eq:ricci-drift-hbc}--\eqref{eq:least-ricci-barrier}, which
use only $Q\geq 0$ and $\mathcal{B}^{T}\geq 0$, give
\begin{equation*}
(L_{f}-2\lambda )\mu \leq 0
\end{equation*}
in the upper-support sense described there. If $\mu >0$ everywhere,
Proposition~\ref{prop:pos-ricci} would make $M$ compact. Thus the closed set
$Z_{0}=\{\mu =0\}$ is nonempty.

We show that $Z_{0}$ is open by using only the scalar strong
minimum principle. At any $p\in Z_{0}$, the two eigenvalues of the Hermitian
Ricci endomorphism are $0$ and $s(p)/2>0$ and each occurs twice as a real
eigenvalue. They remain distinct in a sufficiently small transverse chart,
so the smaller eigenvalue $\mu $ is smooth there. Its upper supports imply
\begin{equation*}
(L_{f}-2\lambda )\mu \leq 0
\end{equation*}
in the classical sense on this chart: if $q\geq \mu $ is such a support at $x
$, then
\begin{equation*}
dq(x)=d\mu (x)\quad \text{and}\quad \Delta _{B}\mu (x)\leq \Delta _{B}q(x).
\end{equation*}
Since $\mu \geq 0$, $\mu (p)=0$, and the zeroth-order coefficient is
$-2\lambda <0$, the ordinary strong minimum principle gives $\mu =0$ on a
smaller connected transverse chart. Basicness extends this vanishing to the
corresponding foliated neighborhood of $p$. Thus $Z_{0}$ is open as well as
closed, and connectedness gives $Z_{0}=M$.

(ii) Parallelism. It follows that $Q$ has real rank two everywhere and
\begin{equation*}
E:=\ker Q\subset D
\end{equation*}
is a smooth basic $J$-invariant real rank-two subbundle. We next prove that
$E$ is transversely parallel. Work on a local leaf space and let $V$ be any
smooth section of $E$; we show that $\nabla _{X}V\in E$ for every transverse
vector $X$.

Write $\nabla $ for its Levi-Civita connection and let $(e_{i})_{i=1}^{4}$
be a local orthonormal frame. Differentiating $Q(V,\cdot )=0$ gives
\begin{equation*}
(\nabla _{e_{i}}Q)(V,W)=-Q(\nabla _{e_{i}}V,W)\qquad \hbox{for every }W.
\end{equation*}
Applying the scalar drift Laplacian to $Q(V,V)\equiv 0$ and using this
identity yields
\begin{equation}
0=(L_{f}Q)(V,V)-2\sum_{i=1}^{4}Q(\nabla _{e_{i}}V,\nabla _{e_{i}}V).
\label{eq:ricci-kernel-parallel}
\end{equation}
Here $L_{f}Q$ denotes the transverse rough drift Laplacian as in
\eqref{eq:ricci-drift-hbc}. The terms involving $Q(V,\nabla V)$ or $Q(V,\nabla
^{2}V)$ vanish. On the other hand, the soliton identity and the
paired-frame contraction \eqref{eq:hbc-contraction} give
\begin{equation*}
(L_{f}Q)(V,V)=-2(\mathrm{Rm}^{T}\ast Q)(V,V)\leq 0.
\end{equation*}
Both terms on the right-hand side of \eqref{eq:ricci-kernel-parallel} are
nonpositive. Consequently both vanish. For a positive semidefinite symmetric
form, $Q(W,W)=0$ implies $W\in \ker Q$. Hence $\nabla _{X}V\in E$ for every
transverse vector $X$. The bundle $E$ is therefore transversely parallel,
and so is its orthogonal complement. Since $E$ is basic, it is preserved
also by $\nabla _{\xi }^{T}=[\xi ,\cdot ]$.

(iii) Global splitting. On each sufficiently small local leaf space, the two
parallel orthogonal plane distributions give a local Riemannian product.
They are $J$-invariant, so both factors are K\"{a}hler curves. The factor
tangent to $E$ has zero Ricci tensor, hence zero Gaussian curvature and is
flat. Its local parallel orthonormal frame extends across the other factor
to a parallel frame of $E$. Lifting this frame to a foliated neighborhood
gives basic sections, which are parallel also in the Reeb direction. Thus
the restriction of the full transverse connection to $E\rightarrow M$ is
flat, including curvature components with a Reeb argument. This is a local
statement and requires neither a global leaf space nor completeness of a
transverse chart.

The pullback of this flat metric connection to the simply connected
universal cover $\widehat{M}$ has trivial monodromy. Parallel transport
therefore supplies a global unit section $V$ of $\widehat{E}$ with $\nabla
^{T}V=0$. By Definition~\ref{def:split}, $M$ is globally split.

Theorem~\ref{thm:split-class}(1)--(3) now gives the stated classification
and scalar bounds. In the nonflat case the transverse metric is locally a
product of a flat plane and a surface of positive Gaussian curvature. Its
curvature operator is nonnegative of rank one. In particular curvature
boundedness follows from the classification.
\end{proof}

\begin{corollary}
\label{cor:nonneg-nondev}Let $(M^{5},\xi ,\eta ,\Phi ,g,f)$ be a complete
noncompact gradient shrinking Sasaki--Ricci soliton whose Reeb foliation is
neither quasi-regular nor globally developable. The following are equivalent:

\begin{enumerate}
\item[(a)] the transverse curvature operator of $M$ is nonnegative;

\item[(a$^{\prime }$)] the transverse holomorphic bisectional curvature of $M
$ is nonnegative;

\item[(b)] $M$ is globally split;

\item[(c)]
\begin{equation*}
M\cong \mathcal{S}_{a}/\Gamma
\end{equation*}
for a unique $a$ with $a_{1}>a_{2}$ and $a_{1}/a_{2}\notin \mathbb{Q}$, and
a cyclic group $\Gamma \cong \mathbb{Z}/p\mathbb{Z}$ as in
Theorem~\ref{thm:split-class}(2).
\end{enumerate}

\noindent In this case $M$ is a flat complex line bundle over a lens space
$L(p,q)$, which is split if and only if $\Gamma $ acts trivially on the
factor $\mathbb{C}$. Moreover $k=2,$ $\dim H_{f}=3$ and $R^{T}$ is bounded
with
\begin{equation*}
\sup R^{T}/\inf R^{T}=e^{c(a)}.
\end{equation*}
The critical set of $f$ consists of two closed Reeb orbits, at the levels
$\psi _{-}<1<\psi _{+}$, and $M$ is not rigid. For quasi-regular solitons with
nondevelopable leaf space, (a), (a$^{\prime }$), (b) and (c) are equivalent
if $a_{1}/a_{2}\notin \mathbb{Q}$ is replaced by $a_{1}/a_{2}\in
\mathbb{Q}\setminus \{1\}$ in (c).
\end{corollary}

\begin{proof}
(i) (a) implies (a$^{\prime }$): a nonnegative curvature operator has
nonnegative sectional curvature, and
$\mathcal{B}^{T}(V,e)=\mathrm{Rm}^{T}(V,e,V,e)+\mathrm{Rm}^{T}(V,Je,V,Je)$ for
unit vectors $V,e$.

(ii) (a$^{\prime }$) implies (b): By Theorem~\ref{thm:nonneg-op}.

(iii) (b) implies (c): If (b) holds, then $R^{T}>0$ by
Theorem~\ref{thm:nondev-structure}(1), and Theorem~\ref{thm:split-class}(2)
gives $M\cong \mathcal{S}_{a}/\Gamma $ with $a_{1}\geq a_{2}$ and $\Gamma $
finite.
If $a_{1}/a_{2}$ were rational, $\mathcal{S}_{a}$ would be quasi-regular by
Proposition~\ref{prop:S3-family}(2), and so would be $M$. Hence
$a_{1}/a_{2}\notin \mathbb{Q}$, so $a_{1}>a_{2}$ and $\Gamma $ is cyclic by
Theorem~\ref{thm:split-class}(2). This is (c).

(iv) (c) implies (a): By Theorem~\ref{thm:split-class}(3).

(v) The remaining assertions follow from Theorem~\ref{thm:split-class}(2)--(4)
and Proposition~\ref{prop:S3-family}(5),(6). Since $\Gamma $
commutes with the torus $T^{3}$ of Proposition~\ref{prop:S3-family}(5), $k$
and $\dim H_{f}$ do not change under the finite quotient.

In the quasi-regular case, $R^{T}>0$ by Corollary~\ref{cor:nondev-qr}, and
in (b)$\Rightarrow $(c) the ratio $a_{1}/a_{2}$ is rational, since otherwise
$M$ would be irregular by Proposition~\ref{prop:S3-family}(3). It is
different from one by Theorem~\ref{thm:split-class}(4), because the leaf
space is not developable. Conversely, by the last assertion of
Theorem~\ref{thm:split-class}(4), these quotients are quasi-regular with
nondevelopable
leaf space.
\end{proof}

\begin{corollary}
Let $(M^5,\xi,\eta,\Phi,g,f)$ be a connected complete gradient shrinking
Sasaki--Ricci soliton with constant $\lambda>0$, normalized by
\eqref{eq:pre-normal}, whose Reeb foliation is not globally developable.
If $M$ is globally split,
then
\begin{equation*}
M\cong \mathcal{S}_{a}/\Gamma
\end{equation*}
for a unique $a$ with $a_{1}>a_{2}>0$, $a_{1}+a_{2}=\lambda /2$, and a
finite group $\Gamma $ of automorphisms of $\mathcal{S}_{a}$ acting freely.
Here $\mathcal{S}_{a}$ is quasi-regular if $a_{1}/a_{2}\in \mathbb{Q}$ and
irregular if $a_{1}/a_{2}\notin \mathbb{Q}$.
\end{corollary}

\begin{proof}
By Theorem~\ref{thm:split-class}, $M$ is $\mathbb{H}_{2}$, a Gaussian
suspension, or a quotient $\mathcal{S}_{a}/\Gamma $ with $a_{1}\geq a_{2}$.
The Reeb foliations of $\mathbb{H}_{2}$ and of the Gaussian suspensions are
globally developable, and so is that of a finite quotient of the regular
cylinder lift, which is the case $a_{1}=a_{2}$
(Theorem~\ref{thm:split-class}(2)(b)). Hence $a_{1}>a_{2}$, and the last
assertion is Proposition~\ref{prop:S3-family}(2),(3).
\end{proof}

\begin{remark}

\begin{enumerate}
\item Naber \cite{Nab} proved that a complete noncompact four-dimensional
shrinker with bounded nonnegative curvature operator is a finite quotient of
\begin{equation*}
\mathbb{R}^{4},\quad \mathbb{S}^{3}\times \mathbb{R}\text{ or
}\mathbb{S}^{2}\times \mathbb{R}^{2}.
\end{equation*}
In the K\"{a}hler case only $\mathbb{C}^{2}$ and $\mathbb{C}\times
\mathbb{CP}^{1}$ remain. The solitons $\mathcal{S}_{a}$ have bounded
nonnegative
transverse curvature operator, but for $a_{1}\neq a_{2}$ their transverse
metrics are not locally isometric to these models.

\item This is a transverse version of Naber's theorem \cite{Nab}, in which
the solitons $\mathcal{S}_{a}$ take the place of the round cylinder. The
removal of the curvature bound corresponds to \cite{MWpos}.
\end{enumerate}
\end{remark}

\medskip

It remains to prove Proposition~\ref{prop:pos-ricci}. We keep the notation
of the proof of Theorem~\ref{thm:nonneg-op}: $Q=\mathrm{Ric}^{T}$, $\mu $ is
the least eigenvalue of $Q$ with respect to $g^{T}$, and $L_{f}=\Delta
_{B}-\langle \nabla ^{T}f,\nabla ^{T}\cdot \,\rangle $ is the transverse
drift Laplacian; on a local leaf space it is the drift Laplacian of the
K\"{a}hler metric $g^{T}$. Basic functions are identified with the functions
they induce on local leaf spaces.

\begin{proof}[Proof of Proposition~\protect\ref{prop:pos-ricci}]
Suppose, to the contrary, that $M$ is noncompact. We first justify that the
normalized potential has compact sublevels on the total space. By
Lemma~\ref{lem:reeb}, either the Reeb closure is a compact torus or the Reeb
flow is a
free proper $\mathbb{R}$-action. In the latter case the quotient
$Z=M/\mathbb{R}$ is a complete smooth K\"{a}hler surface. The basic potential
descends, and the transverse soliton equation becomes
\begin{equation*}
\mathrm{Ric}_{h_Z}+\operatorname{Hess}_{h_Z}f_Z=\lambda h_Z.
\end{equation*}
The quotient has nonnegative holomorphic bisectional curvature and positive
definite Ricci tensor, so it is compact by Wu--Zhang~\cite[Theorem~1.2]{WZh}.
On the other hand, the principal $\mathbb{R}$-bundle $M\to Z$ is trivial
because $\mathbb{R}$ is contractible. A smooth section $\sigma$ therefore
gives a one-form $\alpha=\sigma^*\eta$ satisfying $d\alpha=2\omega_Z$. This
is impossible on the compact K\"{a}hler surface $Z$, since Stokes' theorem
would give
\begin{equation*}
0<\int_Z\omega_Z^2=\frac12\int_Zd(\alpha\wedge\omega_Z)=0.
\end{equation*}
Consequently the Reeb closure $H$ is compact. This argument does not use the
conclusion of Theorem~\ref{thm:nonneg-op} or the classification of
noncompact shrinking K\"{a}hler surfaces.

Normalize $f$ by \eqref{eq:pre-normal}. To apply Lemma~\ref{lem:growth} at its
stated scale, put $\tau =\lambda /6$ and use the constant $D$-homothety
\begin{equation*}
\eta _{\tau }=\tau \eta ,\qquad \xi _{\tau }=\tau ^{-1}\xi ,\qquad g_{\tau
}=\tau g^{T}+\tau ^{2}\eta ^{2}.
\end{equation*}
The same function $f$ satisfies the soliton equation with constant $6$ and the
normalized Hamilton identity for $g_{\tau }$. The metrics $g_{\tau }$ and $g$
are uniformly equivalent, and the Reeb orbits and their closure are unchanged.
Lemma~\ref{lem:growth} thus implies that $f$ is proper on $M$. Properness is
a topological property and remains true at the original scale. Write
$Q=\mathrm{Ric}^T$ and $s=R^T$. Since $f$ is basic, $\nabla f=\nabla^T f$ and
$\Delta_g f=\Delta_B f$; the last equality uses $\nabla_\xi\xi=0$. The
soliton identities yield
\begin{equation}
L_f f=4\lambda-2\lambda f,\qquad |\nabla f|^2=2\lambda f-s,\qquad \nabla
s=2Q(\nabla f,\cdot).  \label{eq:posric-identities}
\end{equation}
Here $\nabla s$ is identified with its metric-dual one-form in the last
identity. Positivity of $Q$ gives $s>0$, hence $f>0$. Properness then
implies that $f$ attains a positive minimum, and that $f(x)\to\infty$
whenever $x$ leaves every compact subset of $M$.

Let $\mu(x)>0$ be the least eigenvalue of $Q$ with respect to $g^T$. It is
continuous and basic, but need not be smooth at a repeated eigenvalue. We
use smooth upper supports to avoid differentiating it. On a local transverse
K\"{a}hler chart the tensor soliton identity is
\begin{equation*}
L_f Q=2\lambda Q-2\mathrm{Rm}^T*Q.
\end{equation*}
If $(e_1,Je_1,e_2,Je_2)$ diagonalizes $Q$ at a point, with eigenvalues
$\mu_1,\mu_2>0$, the K\"{a}hler symmetries give, for any unit $v$,
\begin{equation*}
(\mathrm{Rm}^T*Q)(v,v) =\sum_{j=1}^2\mu_j\mathcal{B}^T(v,e_j)\geq0,
\end{equation*}
as in \eqref{eq:hbc-contraction}. Fix a least unit eigenvector at $x_0$ and
extend it to a smooth unit field $V$ in the transverse chart with $\nabla^T
V(x_0)=0$. Then $q=Q(V,V)$ satisfies $q\geq\mu$ and $q(x_0)=\mu(x_0)$. At
$x_0$ the first-derivative terms involving $V$ vanish, while the unit-length
condition gives
\begin{equation*}
\langle\Delta^T V,V\rangle=-|\nabla^T V|^2=0.
\end{equation*}
Because $Q(V,\cdot)=\mu\langle V,\cdot\rangle$ there, the remaining
second-derivative term also vanishes. Thus
\begin{equation}
L_f q(x_0)=(L_f Q)(V,V)(x_0)\leq2\lambda\mu(x_0).
\label{eq:posric-eigenbarrier}
\end{equation}
In particular no assumption that $\Delta^T V(x_0)=0$ is needed. The tensor
identity used here is the local identity of
\cite[Proposition~\Iref{prop:tensor}{3.1}]{PartI}; it requires no global leaf
space and no
curvature bound.

We now obtain a global positive lower bound for $f\mu$. Set $\vartheta
=f^{-1}+4f^{-2}$. The chain rule and \eqref{eq:posric-identities} give the
exact formulas
\begin{align*}
L_f(f^{-1})&=2\lambda f^{-1}-4\lambda f^{-2} +2f^{-3}|\nabla f|^2, \\
L_f(f^{-2})&=4\lambda f^{-2}-8\lambda f^{-3} +6f^{-4}|\nabla f|^2.
\end{align*}
Consequently, on $\{f\geq16\}$,
\begin{align*}
(L_f-2\lambda)\vartheta &=4\lambda f^{-2}-32\lambda f^{-3}
+(2f^{-3}+24f^{-4})|\nabla f|^2 \\
&\geq2\lambda f^{-2}>0.
\end{align*}
Choose a regular value $F>\max\{16,\min f\}$ and then choose $c_{0}>0$ so
small that $w:=\mu-c_{0}\vartheta>0$ on the compact level $\{f=F\}$. Since
$\mu>0$ and $\vartheta\to0$ at infinity, $\liminf_{x\to\infty}w(x)\geq0$. If
$w$ were negative somewhere in $\{f\geq F\}$, its negative infimum would
therefore be achieved on a compact subset. The positive boundary values force
a minimum point $x_0$ to lie in $\{f>F\}$. Since the upper support $q$
constructed above satisfies $q\geq\mu$ near $x_0$ with equality at $x_0$, the
smooth function $q-c_{0}\vartheta$ also has a local minimum at $x_0$. Hence
\begin{equation*}
0\leq L_f(q-c_{0}\vartheta)(x_0)
<2\lambda\bigl(\mu(x_0)-c_{0}\vartheta(x_0)\bigr)<0,
\end{equation*}
a contradiction. It follows that $\mu\geq c_{0}\vartheta\geq c_{0}/f$ on the
exterior. Since $f\mu$ has a positive minimum on the compact sublevel $\{f\leq
F\}$, decreasing $c_{0}$ if necessary yields
\begin{equation}
Q\geq\frac{c_{0}}{f}\,g^T\qquad\text{on }M.  \label{eq:posric-lower}
\end{equation}

Next we prove uniform scalar-curvature growth using the negative gradient
flow. For any $x\in M$, let $x(t)$ solve $x^{\prime }(t)=-\nabla f(x(t))$
with $x(0)=x$. Since
\begin{equation*}
\frac{d}{dt}f(x(t))=-|\nabla f(x(t))|^2,
\end{equation*}
the trajectory remains in the compact set $\{f\leq f(x)\}$ and exists for
every $t\geq0$. Its potential decreases to a limit $\ell\geq\min f>0$, and
\begin{equation*}
\int_0^\infty |\nabla f(x(t))|^2\,dt=f(x)-\ell<\infty.
\end{equation*}
There is a sequence $t_j\to\infty$ for which $|\nabla f(x(t_j))|\to0$.
Passing to a subsequence in the compact sublevel gives $x(t_j)\to x_*$,
where
\begin{equation*}
\nabla f(x_*)=0,\qquad f(x_*)=\ell,\qquad s(x_*)=2\lambda\ell.
\end{equation*}
Convergence of the entire trajectory is not required. By
\eqref{eq:posric-identities} and \eqref{eq:posric-lower},
\begin{equation*}
\frac{d}{dt}\bigl(s(x(t))-2c_{0}\log f(x(t))\bigr) =-2Q(\nabla f,\nabla
f)+\frac{2c_{0}}{f}|\nabla f|^2\leq0.
\end{equation*}
Taking the limit along $t_j$ shows that
\begin{equation}
\begin{split}
s(x)-2c_{0}\log f(x) &\geq 2\lambda\ell-2c_{0}\log\ell \\
&\geq\inf_{u>0}(2\lambda u-2c_{0}\log u)
=2c_{0}\bigl(1-\log(c_{0}/\lambda)\bigr).
\end{split}
\label{eq:posric-scalar-growth}
\end{equation}
The last infimum is attained at $u=c_{0}/\lambda$ and is independent of $x$.
Thus $s(x)\to\infty$ at infinity. This argument allows critical points at
arbitrary levels and does not assume that every gradient trajectory meets
one prescribed level set.

It remains to compare this growth with the integrated trace equation.
Lemma~\ref{lem:transverse-CY}, proved in Appendix~\ref{app:B} using only
$Q\geq0$
and compactness of the Reeb closure, gives $\operatorname{Vol}_g(M)=\infty$.
For a
regular value $t$ let $\Omega_t=\{f<t\}$. Its closure is compact, its
outward unit normal is $\nabla f/|\nabla f|$, and $\Delta_g f+s=4\lambda$.
The divergence theorem therefore gives
\begin{equation*}
\int_{\Omega_t}s\,d\mathrm{vol}_g =4\lambda\operatorname{Vol}_g(\Omega_t)
-\int_{\partial\Omega_t}|\nabla f|\,dA_g
\leq4\lambda\operatorname{Vol}_g(\Omega_t).
\end{equation*}
Choose a fixed compact sublevel $K$ such that $s\geq8\lambda$ outside $K$,
which is possible by \eqref{eq:posric-scalar-growth}. For every sufficiently
large regular $t$, positivity of $s$ implies
\begin{equation*}
8\lambda\bigl(\operatorname{Vol}_g(\Omega_t)-\operatorname{Vol}_g(K)\bigr)
\leq\int_{\Omega_t}s\,d\mathrm{vol}_g
\leq4\lambda\operatorname{Vol}_g(\Omega_t).
\end{equation*}
Thus $\operatorname{Vol}_g(\Omega_t)\leq2\operatorname{Vol}_g(K)$. Choose an
increasing
sequence of regular values tending to infinity; the corresponding sublevels
exhaust $M$. Monotone convergence would give
$\operatorname{Vol}_g(M)\leq2\operatorname{Vol}_g(K)<\infty$, a contradiction.
\end{proof}

\appendix

\section{The Irregular Models \texorpdfstring{$S_{a}$
and
$\mathcal{S}_{a}$}{Sa and script Sa}}

\label{app:A}

\noindent (I) \emph{Computation of the irregular models.}

\noindent (A) \emph{The three-dimensional model $S_{a}$ of
Example~\ref{ex:S3}.} Represent $\mathbb{S}^{3}$ by $(t,\theta _{1},\theta
_{2})$, $0\leq t\leq 1$, with both angles of period $2\pi $, collapsing the
$\theta
_{2}$ circle at $t=0$ and the $\theta _{1}$ circle at $t=1$. Define on the
interior
\begin{equation}
\eta _{S}=\frac{1-t}{a_{1}}d\theta _{1}+\frac{t}{a_{2}}d\theta _{2},\qquad
\beta =\frac{d\theta _{2}}{a_{2}}-\frac{d\theta _{1}}{a_{1}}.
\label{eq:etaS}
\end{equation}
The one-form $\eta _{S}$ extends smoothly across the collapsed circles. It
has $\eta _{S}(\xi _{S})=1$ and $d\eta _{S}=dt\wedge \beta $.

Let $c$ be the unique real solution of
\begin{equation}
\mathcal{H}(c):=\frac{1}{c}-\frac{1}{e^{c}-1}=\frac{a_{2}}{a_{1}+a_{2}}.
\label{eq:Hc}
\end{equation}
Indeed, $\mathcal{H}(0)=1/2$,
\begin{equation*}
\mathcal{H}^{\prime }(c)=-\frac{1}{c^{2}}+\frac{1}{4\sinh ^{2}(c/2)}<0\quad
(c\neq 0),\qquad \mathcal{H}(-\infty )=1,\qquad \mathcal{H}(+\infty )=0.
\end{equation*}
One has $c=0$ if and only if $a_{1}=a_{2}$, and $c>0$ if and only if
$a_{1}>a_{2}$. For $c\neq 0$ put
\begin{equation}
V(t)=\frac{\lambda }{c}\Bigl(t-\frac{e^{ct}-1}{e^{c}-1}\Bigr).  \label{eq:V}
\end{equation}
and for $c=0$ let $V(t)=\frac{\lambda }{2}t(1-t)$, the limit of \eqref{eq:V}
as $c\rightarrow 0$. Direct differentiation and \eqref{eq:Hc} give
\begin{gather}
V(0)=V(1)=0,\qquad V^{\prime }(0)=2a_{2},\qquad V^{\prime }(1)=-2a_{1},
\label{eq:Vboundary} \\
V^{\prime \prime }-cV^{\prime }=-\lambda ,\qquad -V^{\prime \prime
}=\frac{\lambda ce^{ct}}{e^{c}-1}>0.  \label{eq:VODE}
\end{gather}
Strict concavity and the boundary values imply $V>0$ on $(0,1)$. Define
\begin{equation}
h_{S}^{T}=\frac{dt^{2}}{2V(t)}+\frac{V(t)}{2}\beta ^{2},\qquad
g_{S}=h_{S}^{T}+\eta _{S}^{2},\qquad \psi =ct.  \label{eq:S3metric}
\end{equation}
In this calculation $\psi=ct$ is an unnormalized potential. The normalized
potential used for $S_a$ elsewhere is
$\psi_a=ct+c/(e^c-1)$, with the added constant interpreted as $1$ when
$c=0$; this change does not affect its Hessian or soliton equation.
The tensors in \eqref{eq:S3metric} are invariant under the two angular
circles, and $h_{S}^{T}$ annihilates $\xi _{S}$. Locally, $h_{S}^{T}$ is the
momentum construction of Hwang--Singer \cite{HS} with momentum coordinate
$t/2$ and momentum profile $V/2$.

Near $t=0$ put $t=a_{2}\rho ^{2}$. Since $V(t)=2a_{2}t+O(t^{2})$, the radial
part is $(1+O(\rho ^{2}))d\rho ^{2}$, and the angular term is
\begin{equation*}
\rho ^{2}\bigl(d\theta _{2}-(a_{2}/a_{1})d\theta _{1}\bigr)^{2}+O(\rho
^{4})\bigl(d\theta _{2}-(a_{2}/a_{1})d\theta _{1}\bigr)^{2}.
\end{equation*}
The leading radial and collapsing angular coefficients agree, giving a
smooth disk metric in $\rho \cos \theta _{2}$, $\rho \sin \theta _{2}$. The
higher coefficients are smooth functions of $\rho ^{2}$, and the cross terms
with the noncollapsing circle are smooth. The contact form is
\begin{equation*}
\eta _{S}=a_{1}^{-1}d\theta _{1}+\rho ^{2}\bigl(d\theta
_{2}-(a_{2}/a_{1})d\theta _{1}\bigr).
\end{equation*}
Thus both metric and contact form extend smoothly. Near $t=1$ use
$1-t=a_{1}\rho ^{2}$ and exchange the two angles. The second boundary slope
gives the same conclusion. The tensor $g_{S}$ is positive definite
everywhere.

The transverse area form is $\omega _{S}^{T}=\frac{1}{2}dt\wedge \beta
=\frac{1}{2}d\eta _{S}$. The oriented two-dimensional transverse metric
defines an invariant complex structure. In complex dimension one it is
integrable. Locally the construction is the K\"{a}hler contactization of
\cite[Proposition~\Iref{prop:contactization}{2.1}]{PartI}, so the smooth
global structure is Sasakian. In particular $\xi _{S}$ is its unit Killing
Reeb field.

On a local quotient write $du=\beta $. The metric is $dt^{2}/(2V)+(V/2)du^{2}
$. Introducing arc length $d\varrho =dt/\sqrt{2V}$
and the warping function $b=\sqrt{V/2}$ gives $db/d\varrho =V^{\prime }/2$,
hence
\begin{equation}
K=-\frac{1}{b}\frac{d^{2}b}{d\varrho ^{2}}=-V^{\prime \prime }=\frac{\lambda
ce^{ct}}{e^{c}-1}>0.  \label{eq:KS3}
\end{equation}
For $\psi =ct$ the Hessian in the unit radial and angular directions has the
same value $cV^{\prime }$. Its mixed component vanishes. Therefore
\begin{equation*}
\operatorname{Hess}_{S}^{T}\psi =cV^{\prime }h_{S}^{T},\qquad
\mathrm{Ric}_{S}^{T}+\operatorname{Hess}_{S}^{T}\psi =(-V^{\prime \prime
}+cV^{\prime })h_{S}^{T}=\lambda
h_{S}^{T}.
\end{equation*}
The curvature is positive, and nonconstant unless $c=0$. Finally the Reeb
flow is $(\theta _{1},\theta _{2})\mapsto (\theta _{1}+a_{1}s,\theta
_{2}+a_{2}s)$. If $a_{1}/a_{2}\notin \mathbb{Q}$, it is dense on every
interior two-torus; there are two closed endpoint orbits, but the generic
orbits are not closed, and the structure is irregular. If $a_{1}/a_{2}\in
\mathbb{Q}$, all orbits are closed (Proposition~\ref{prop:S3-family}(2)).

\noindent (B) \emph{The five-dimensional model $\mathcal{S}_{a}$ of
Example~\ref{ex:counterexample}.} Locally the Reeb quotient is the product of
the
plane and a transverse chart of $\mathbb{S}^{3}$. Its K\"{a}hler metric is
$h_{\mathrm{E}}+h_{S}^{T}$ and its K\"{a}hler form is $\omega
_{\mathrm{E}}+\omega _{S}^{T}$. Since $d\eta =2\omega _{\mathrm{E}}+d\eta
_{S}$, the
local contactization of
\cite[Proposition~\Iref{prop:contactization}{2.1}]{PartI} proves the Sasakian
condition. The Reeb field remains $\xi _{S}$, so
the Reeb type is that of $S_{a}$.

The transverse Ricci tensor is zero on the plane and $Kh_{S}^{T}$ on the
surface. The Hessian of $\lambda |z|^{2}/2$ is $\lambda h_{\mathrm{E}}$,
while Example~\ref{ex:S3} provides the other factor. This proves the full
real soliton equation. The potential is basic and proper because
$\mathbb{S}^{3}$ is compact. To check completeness, note that any
finite-length
escaping curve projects to a Cauchy curve in the Euclidean plane, because
$g\geq dx^{2}+dy^{2}$. Its plane coordinates therefore remain bounded, and
its $\mathbb{S}^{3}$ coordinate lies in a compact space. It remains in a
compact subset on which $g$ is a smooth positive metric, contradicting
finite-length escape.

For the normalization, on the three-dimensional factor $|\nabla ^{T}\psi
|^{2}=2c^{2}V$, and
\begin{equation*}
K+c^{2}V-\lambda ct=\frac{\lambda c}{e^{c}-1}.
\end{equation*}
The latter identity follows directly from \eqref{eq:V} and \eqref{eq:KS3}.
Adding the Gaussian plane term gives exactly the Hamilton identity for the
constant in \eqref{eq:ctr-g}. Formula \eqref{eq:ctr-RT} is immediate.

We next prove the last assertion of Example~\ref{ex:counterexample}. Let
$a_{1}\neq a_{2}$, so that $c\neq 0$. Let $q_{\ast }\in \mathbb{S}^{3}$ be a
point of the endpoint orbit $\{t=1\}$ if $c>0$, resp.\ $\{t=0\}$ if $c<0$,
and put $p_{r}=(r,q_{\ast })\in \mathbb{C}\times \mathbb{S}^{3}$ for $r\geq
0 $. Denote the complete Reeb flow of $\mathbb{S}^{3}$ by $\varphi _{s}$.
For $v=(v_{1},v_{2})\in \mathbb{R}^{2}$ put
\begin{equation*}
b_{v}(x,y)=v_{1}y-v_{2}x,\qquad I_{v}(z,q)=\bigl(z+v,\varphi
_{-b_{v}(z)}q\bigr).
\end{equation*}
Since translation changes $\alpha $ by $db_{v}$, while the variable Reeb
shift changes $\eta _{S}$ by $-db_{v}$, one has $I_{v}^{\ast }\eta =\eta $.
Reeb invariance and horizontality of $h_{S}^{T}$ give $I_{v}^{\ast
}h_{S}^{T}=h_{S}^{T}$. Thus $I_{v}$ preserves $g,\eta ,\xi $ and the
transverse complex structure. For $v=(-r,0)$, $b_{v}(r,0)=0$, hence
$I_{v}(p_{r})=p_{0}$. These are exact pointed Sasakian isometries. They do
not preserve the uncentered potential: $I_{v}^{\ast }f-f$ is an affine
function of $z$, so the pulled-back potentials still satisfy the local
soliton equation \eqref{LS}. This diverging affine term is the familiar
reason to distinguish a splitting function from a limiting shrinking
potential.

At $q_{\ast }$ one has $\nabla ^{T}\psi =0$, since $V$ vanishes at the
endpoints. Along the positive real axis $\alpha $ vanishes on the radial
direction, so the full transverse gradient at $p_{r}$ is $\lambda r\partial
_{r}$. Thus $r\mapsto p_{r}$ parametrizes an escaping gradient trajectory,
up to the reparametrization $r=e^{\lambda s}r_{0}$, and
$R^{T}(p_{r})=2K(q_{\ast })$ for all $r$.

It remains to check the strict inequality and the local obstruction. By
\eqref{eq:Vboundary},
\begin{equation*}
\int_{0}^{1}K(t)\,dt=V^{\prime }(0)-V^{\prime }(1)=2(a_{1}+a_{2})=\lambda .
\end{equation*}
Also $K^{\prime }=cK$, so the maximum of $K$ occurs at the same endpoint as
the maximum of $ct$. Since $K$ is nonconstant, $K(q_{\ast })>\lambda $.
Every neighborhood of that endpoint in a local transverse disk contains
different $t$ values and thus different $K$ values. The product transverse
metric therefore has nonconstant scalar curvature on every neighborhood of
the basepoint. A round K\"{a}hler cylinder has constant scalar curvature
$2\lambda $, so no such neighborhood can be transversely isometric to it.

\noindent (II) \emph{Boundary slopes and the scalar ODE.}

Differentiating \eqref{eq:V} gives
\begin{equation*}
V^{\prime }(t)=\frac{\lambda }{c}\Bigl(1-\frac{ce^{ct}}{e^{c}-1}\Bigr),\qquad
V^{\prime \prime }(t)=-\frac{\lambda ce^{ct}}{e^{c}-1}.
\end{equation*}
Therefore, by \eqref{eq:Hc} and $\lambda =2(a_{1}+a_{2})$,
\begin{equation*}
V^{\prime }(0)=\lambda \Bigl(\frac{1}{c}-\frac{1}{e^{c}-1}\Bigr)=\lambda
\mathcal{H}(c)=2a_{2},
\end{equation*}
and
\begin{equation*}
V^{\prime }(1)=\lambda \Bigl(\frac{1}{c}-\frac{e^{c}}{e^{c}-1}\Bigr)=\lambda
\bigl(\mathcal{H}(c)-1\bigr)=-2a_{1}.
\end{equation*}
Also $V^{\prime \prime }-cV^{\prime }=-\lambda $. When $c=0$ the continuous
limit is $V=\frac{\lambda }{2}t(1-t)$, and then $K\equiv \lambda $. This is
precisely the symmetric case $a_{1}=a_{2}$.

\section{Proofs of the Global Splitting Results}

\label{app:B}

\begin{lemma}
\label{lem:split-affine}Let $(M^{5},\xi ,\eta ,\Phi ,g,f)$ be a connected
complete gradient shrinking Sasaki--Ricci soliton with constant $\lambda $,
normalized by \eqref{eq:pre-normal}, and let $V$ be a section of $D$ with
$|V|=1$ and
\begin{equation*}
\nabla ^{T}V=0.
\end{equation*}
Put
\begin{equation*}
u_{1}=\lambda ^{-1}Vf,\qquad u_{2}=\lambda ^{-1}(\Phi V)f,\qquad
u=u_{1}+iu_{2},\qquad \psi =f-\tfrac{\lambda }{2}|u|^{2}.
\end{equation*}

\begin{enumerate}
\item $\Phi V$ is transversely parallel,
\begin{equation*}
\mathrm{Rm}^{T}(V,\cdot ,\cdot ,\cdot )=\mathrm{Rm}^{T}(\Phi V,\cdot ,\cdot
,\cdot )=0,
\end{equation*}
and
\begin{equation*}
\mathrm{Ric}^{T}(V,\cdot )=\mathrm{Ric}^{T}(\Phi V,\cdot )=0.
\end{equation*}

\item $u_{1},u_{2},\psi $ are basic with
\begin{equation*}
\nabla u_{1}=V,\quad \nabla u_{2}=\Phi V,\quad
\operatorname{Hess}^{T}u_{1}=\operatorname{Hess}^{T}u_{2}=0,\quad
\text{and}\quad V\psi =(\Phi V)\psi =0.
\end{equation*}

\item The fields $X_{1}=V-2u_{2}\xi $ and $X_{2}=\Phi V+2u_{1}\xi $ are
complete Killing fields which preserve $\xi ,\eta ,\Phi ,\psi ,V$ and $\Phi
V $, and
\begin{equation*}
X_{1}u=1,\qquad X_{2}u=i,\qquad \lbrack X_{1},X_{2}]=2\xi ,\qquad \lbrack
\xi ,X_{1}]=[\xi ,X_{2}]=0.
\end{equation*}
\end{enumerate}
\end{lemma}

\begin{proof}
(1) Since $\nabla ^{T}\Phi |_{D}=0$, $\Phi V$ is parallel, and
$\mathrm{Rm}^{T}(X,Y)V=\mathrm{Rm}^{T}(X,Y)\Phi V=0$. The symmetries of
$\mathrm{Rm}^{T}$
give the remaining identities.

(2) Since $\eta ([\xi ,V])=-d\eta (\xi ,V)=0$, one has $\nabla _{\xi
}^{T}V=[\xi ,V]=0$. Hence $\xi (Vf)=V(\xi f)=0$ and $u_{1}$ is basic. For
$X\in D$, by \eqref{eq:soliton} and (1),
\begin{equation*}
X(Vf)=\operatorname{Hess}^{T}f(X,V)+df(\nabla _{X}^{T}V)=\lambda
h(X,V)-\mathrm{Ric}^{T}(X,V)=\lambda h(X,V).
\end{equation*}
Thus $\nabla u_{1}=\nabla ^{T}u_{1}=V$, and $\operatorname{Hess}^{T}u_{1}=\nabla
^{T}V=0$. Likewise for $u_{2}$. Then $V\psi =\lambda u_{1}-\lambda
(u_{1}Vu_{1}+u_{2}Vu_{2})=0$, because $Vu_{1}=1$ and $Vu_{2}=h(\Phi V,V)=0.$
Similarly $(\Phi V)\psi =0$.

(3) As is standard (cf.~\cite{BG}), the contact vector field with basic
Hamiltonian $F$ is $F\xi +\frac{1}{2}\Phi \nabla ^{T}F$. For $F=-2u_{2}$ and
$F=2u_{1}$, this gives $X_{1}$ and $X_{2}$ due to $\Phi ^{2}V=-V$. Since
$\operatorname{Hess}^{T}u_{j}=0$ is $J$-invariant, the argument in the proof of
Proposition~\ref{prop:soliton-torus}(1)
shows
that $X_{1},X_{2}$ preserve $\xi ,\eta ,\Phi $ and $g$. Killing fields on a
complete manifold are complete. By (2), $X_{1}u=Vu=1$, $X_{2}u=(\Phi V)u=i$
and $X_{j}\psi =0$. With $\Phi =\nabla \xi $ and (1), $\nabla _{V}\Phi
V=-g(\Phi V,\Phi V)\xi =-\xi $ and $\nabla _{\Phi V}V=\xi $, so $[V,\Phi
V]=-2\xi $. Together with $[\xi ,V]=[\xi ,\Phi V]=0$ and the derivatives of
$u_{1},u_{2}$ in (2), this gives $[X_{j},V]=[X_{j},\Phi V]=0$,
$[X_{1},X_{2}]=2\xi $ and $[\xi ,X_{j}]=0$.
\end{proof}

\begin{proof}[Proof of Theorem~\protect\ref{thm:split-structure}]
(1) By Lemma~\ref{lem:split-affine}(2), $u:M\rightarrow \mathbb{C}$ is a
submersion with orthonormal gradients $V,\Phi V$, so $N$ is a closed
embedded submanifold with $TN=\mathbb{R}\xi \oplus E$, $E=\{V,\Phi
V\}^{\perp }\cap D$, and $\Phi E=E$. For a basic function $w,$ one has
$\operatorname{Hess}_{g}w=\operatorname{Hess}^{T}w$ on $D\times D$,
$\operatorname{Hess}_{g}w(\xi
,\xi )=0$ and $\operatorname{Hess}_{g}w(\xi ,X)=-(\Phi X)w$ for $X\in D$
(since $\xi
w=0$ and $\nabla _{X}\xi =\Phi X$). For $w=u_{1},u_{2}$ and $X\in E$ all
these terms vanish, since $\operatorname{Hess}^{T}u_{j}=0$ and $\Phi X\perp
V,\Phi V$. Hence the second fundamental form of $N$, which is
$-\sum_{j}\operatorname{Hess}_{g}u_{j}|_{TN}\,\nabla u_{j}$, vanishes. Thus
$N$ is totally geodesic and
complete, $\xi $ restricts to a unit Killing field with $\nabla ^{N}\xi
=\Phi |_{TN}$, and $N$ is Sasakian with $D_{N}=E$. Since $V,\Phi V$ are
parallel, $\nabla ^{T}$ preserves $E$ and induces the transverse connection
of $N$. By Lemma~\ref{lem:split-affine}(1),
$\mathrm{Ric}^{T}|_{E}=\mathrm{Ric}^{T_{N}}$ and $R^{T}=R^{T_{N}}$. As $\nabla
^{T}\psi \in E$ and $\operatorname{Hess}^{T}\psi
=\operatorname{Hess}^{T}f-\lambda \sum_{j}du_{j}^{2}$, the soliton
equation of $M$ restricts to that of $(N,\psi _{N})$, and on $N$, where $u=0$,
the identity \eqref{eq:pre-normal} becomes
\begin{equation*}
R^{T_{N}}+|\nabla \psi _{N}|^{2}=2\lambda \psi _{N}.
\end{equation*}
By Lemma~\ref{lem:split-affine}(3) the maps $\phi _{x}^{1}$, $\phi _{y}^{2}$,
$\varphi _{s}$ are automorphisms of $(\xi ,\eta ,\Phi ,g)$ preserving $\psi
,V,\Phi V$ and $E$, $\varphi $ commutes with $\phi ^{1}$ and $\phi ^{2}$, and
$\phi _{x}^{1},\phi _{y}^{2}$ add $x$, resp.\ $iy$, to $u$. Hence $u\circ
\Theta (z,q)=z$, and $\Theta $ is a diffeomorphism with inverse $p\mapsto
\bigl(u(p),\varphi _{-u_{1}u_{2}}\circ \phi _{-u_{2}}^{2}\circ \phi
_{-u_{1}}^{1}(p)\bigr)$, $u_{j}=u_{j}(p)$. In particular $N$ is connected.
Moreover $\Theta ^{\ast }\psi =\psi _{N}$, so $\Theta ^{\ast }f=\frac{\lambda
}{2}|z|^{2}+\psi _{N}$, and $\Theta _{\ast }\xi _{N}=\xi $. Let $p=\Theta
(z,q)$. From $\Theta (x+s,y,q)=\phi _{s}^{1}\varphi _{sy}\Theta
(z,q),$ we get $d\Theta (\partial _{x})=X_{1}+y\xi =V-y\xi $ at $p$. From
$\Theta (x,y+s,q)=\phi _{x}^{1}\phi _{s}^{2}\varphi _{xs}(\phi
_{x}^{1})^{-1}\Theta (z,q),$ we get $d\Theta (\partial _{y})=d\phi
_{x}^{1}(X_{2}+x\xi )$, evaluated at $p^{\prime }=\phi _{-x}^{1}(p)$, where
$u_{1}(p^{\prime })=0$. Hence $\eta (d\Theta \partial _{y})=2u_{1}(p^{\prime
})+x=x$, and the $D$-component of $d\Theta \partial _{y}$ is $\Phi V(p)$.
For $w\in T_{q}N$, $d\Theta (w)$ is the image of $w$ under the automorphism
$\phi _{x}^{1}\phi _{y}^{2}\varphi _{xy}$, which maps $E_{q}$ isometrically
and $\Phi $-linearly onto $E_{p}$. Therefore
\begin{equation*}
\Theta ^{\ast }\eta =\eta _{N}+x\,dy-y\,dx,\qquad \Theta ^{\ast }(g-\eta
^{2})=dx^{2}+dy^{2}+h_{N},
\end{equation*}
and $\Theta $ maps the product transverse complex structure to that of $M$.
Since $g$ and $\Phi $ are determined by $\eta $, $\xi $ and the transverse
K\"{a}hler structure, $\Theta $ is an isomorphism of solitons.

(2) Apply (1) to the universal cover. The factor $N$ is simply connected
because $\mathbb{C}\times N$ is.

(3) A transversely parallel section $W$ of $D$ has parallel component in $E$.
If this component were nonzero somewhere, it would be a nonzero parallel
section of $D_{N}$ along $N$, and by Lemma~\ref{lem:split-affine}(1) applied
to it the Gaussian curvature of $h_{N}$ would vanish identically,
contradicting $R^{T_{N}}\not\equiv 0$. Hence $W=\alpha V+\beta \Phi V$ with
constants $\alpha ,\beta $.

Let $\gamma \in \operatorname{Aut}(\mathbb{C}\times _{\mathrm{H}}N,f)$. Then
$\gamma
_{\ast }$ preserves this two-dimensional space of parallel sections,
commutes with $\Phi $ and preserves $h$, so it acts on it by an element of
$U(1)$. Since $\gamma ^{\ast }f=f$, the definition of $u$ gives $u\circ
\gamma =\rho \,u$ for some $\rho \in U(1)$. Composing with the rotation by
$\rho ^{-1}$ we obtain an automorphism $\gamma ^{\prime }$ with $u\circ
\gamma ^{\prime }=u$. It preserves $u_{1},u_{2}$ and then $V,\Phi V$,
$X_{1},X_{2}$ as well as their flows, and it commutes with $\varphi $. By
\eqref{eq:split-iso}, $\gamma ^{\prime }(z,q)=(z,\gamma _{N}(q))$ with $\gamma
_{N}=\gamma ^{\prime }|_{\{0\}\times N}\in \operatorname{Aut}(N,\psi _{N})$.
Conversely both factors act by automorphisms. An element $(\rho ,\gamma
_{N}) $ fixes a point of $\{0\}\times N$ whenever $\gamma _{N}$ has a fixed
point; hence freeness of $\Gamma $ implies that $\gamma \mapsto \gamma _{N}$
is injective and that $\Gamma _{N}$ acts freely. Finally a parallel section on the
quotient lifts to a $\Gamma $-invariant $\alpha V+\beta \Phi V$ which forces
$\chi \equiv 1$. If $\chi \equiv 1$, then $V$ descends and the quotient is
$\mathbb{C}\times _{\mathrm{H}}(N/\Gamma _{N})$.
\end{proof}

\begin{proof}[Proof of Proposition~\protect\ref{prop:3d-shrinkers}]
By \eqref{eq:pre-identities} with $n=1$, $\nabla
^{T}R^{T}=2\mathrm{Ric}^{T}(\nabla ^{T}\psi ,\cdot )$ and $dK=K\,d\psi $, as
in \eqref{eq:dK}.
Hence $d(Ke^{-\psi })=0$, so $K=\kappa e^{\psi }$ with a constant $\kappa $.
If $N$ is noncompact, the dimension-general scalar-nonnegativity lemma
\cite[Lemma~\Iref{lem:nonneg}{3.3}]{PartI}, applied with $n=1$ to the complete
shrinker $N$, gives $K\geq 0$. If $N$ is compact, apply the minimum
principle to $L_{\psi }R^{T}=2\lambda R^{T}-2|\mathrm{Ric}^{T}|^{2}$ at a
minimum of $R^{T}$ and a negative minimum is impossible. Thus $\kappa \geq 0$
in either case, and the normalization gives the second identity. If $\kappa
>0$, the nonnegativity of $|\nabla ^{T}\psi |^{2}$ means $\lambda \psi
e^{-\psi }\geq \kappa $. Since $\max_{x}xe^{-x}=1/e$ is attained only at $x=1
$, this gives $\kappa \leq \lambda /e$ and $\psi (N)\subset \lbrack \psi
_{-},\psi _{+}]$, with $\psi _{\pm }=1$ if $\kappa =\lambda /e$.

(1) If $\kappa =0$, then $K\equiv 0$ and then $\mathrm{Rm}^{T}\equiv 0$. If
$N$ were compact, integration of $R^{T}+\Delta _{B}\psi =2\lambda $ would
give $0=2\lambda \operatorname{Vol}(N)$. Steps~1 and~2 of the proof of
Proposition~\ref{prop:transflat}, with $\mathbb{C}^{2}$ replaced by
$\mathbb{C}$,
identify the simply connected case with $\mathbb{H}_{1}$.

(2) By the above, $\mathrm{Ric}^{T}=Kh\geq \lambda \psi _{-}h$ with $\lambda
\psi _{-}>0$. For $0<\tau <\min \{1,\lambda \psi _{-}/2\}$, the $D$-homothetic
metric $g_{\tau }=\tau g^{T}+\tau ^{2}\eta ^{2}\geq \tau ^{2}g$
is complete, and \eqref{eqn-SR2} and \eqref{eq:pre-ric}, applied to the
structure $g_{\tau }$, give $\mathrm{Ric}_{g_{\tau
}}|_{D}=\mathrm{Ric}^{T}-2\tau g^{T}\geq (\lambda \psi _{-}-2\tau )g^{T}$,
$\mathrm{Ric}_{g_{\tau }}(\xi _{\tau },\xi _{\tau })=2$ and
$\mathrm{Ric}_{g_{\tau }}(\xi
_{\tau },D)=0$. Thus $\mathrm{Ric}_{g_{\tau }}\geq \delta g_{\tau }$ for
some $\delta >0$, and by Myers' theorem $N$ is compact with finite
fundamental group. The lifted soliton on the compact universal cover
satisfies \eqref{eq:3d-identities} with the same $\kappa $. We may assume
from now on that $N$ is simply connected.

(i) $\kappa =\lambda /e$. Then $\psi \equiv 1$ and $K\equiv \lambda $. After
the $D$-homothety which changes $\lambda $ into $4$, \eqref{eqn-SR2} and
\eqref{eq:pre-ric} give $\mathrm{Ric}_{g}=2g$, so $g$ has constant curvature
one and $N$ is the round sphere. Its Reeb field is a unit Killing field
$x\mapsto Bx$, $B\in \mathfrak{so}(4)$, so $B^{2}=-I$. All such $B$ are
conjugate under $O(4)$, and $\eta =g(\xi ,\cdot )$, $\Phi =\nabla \xi $.
Hence $N$ is the Hopf structure, which is $S_{(1,1)}$ at the normalization
$\lambda =4$ (Example~\ref{ex:S3}). That is $S_{a}$ with $a_{1}=a_{2}=\lambda
/4$ at the original normalization.

(ii) $\kappa <\lambda /e$. If $\psi $ were constant, then
$\operatorname{Hess}^{T}\psi =0$ would give $K\equiv \lambda $ and $\psi
\equiv 1$, hence $\kappa =\lambda /e$. So $\psi $ is not constant. At its
extrema $\nabla
^{T}\psi =0$, so $\psi \in \{\psi _{-},\psi _{+}\}$ there. Hence $\psi
(N)=[\psi _{-},\psi _{+}]$ and the critical set is $\psi ^{-1}(\psi
_{-})\cup \psi ^{-1}(\psi _{+})$ on whose normal bundle $D$ the Hessian
$(\lambda -K)h=\lambda (1-\psi _{\pm })h$ is definite. As in the proof of
Proposition~\ref{prop:min-set}, $\psi ^{-1}(\psi _{\pm })$ is the zero set
of the Killing field $\zeta _{\psi }-\psi _{\pm }\xi $, where $\zeta _{\psi
}=\psi \xi +\frac{1}{2}\Phi \nabla ^{T}\psi $ is the contact vector field
with basic Hamiltonian $\psi $, and its components are closed Reeb orbits.
By Proposition~\ref{prop:soliton-torus}(2), the closure $T$ of the group
generated by the flows of $\xi $ and $\zeta _{\psi }$ is a torus of
dimension at most two. Since $\zeta _{\psi }-\psi \xi \neq 0$ on $N^{\circ
}=\psi ^{-1}((\psi _{-},\psi _{+}))$, $T\cong T^{2}$ and its Lie algebra
$\mathfrak{t}$ is spanned by $\xi ,\zeta _{\psi }$. In particular $\eta (X)$
is an affine function of $\psi $ for every $X\in \mathfrak{t}$. On $N^{\circ
}$ the vectors $\xi ,\nabla ^{T}\psi ,\Phi \nabla ^{T}\psi $ span the
tangent space and are fixed by the isotropy, so $T$ acts freely there, and
$\psi :N^{\circ }\rightarrow (\psi _{-},\psi _{+})$ is a principal $T$-bundle.
Each component of $N^{\circ }$ is a product $(\psi _{-},\psi
_{+})\times T^{2}$. Since the normal Hessians are definite, a punctured
tubular neighborhood of a critical orbit is connected, so the closure of
such a component is open and closed in $N$ and contains exactly one orbit of
each critical level. As $N$ is connected, $N^{\circ }$ is connected and $\psi
^{-1}(\psi _{\pm })$ is a single closed Reeb orbit $\gamma _{\pm }$.
The critical levels are the images of these orbits on a quotient of $N$. The
isotropy group of $\gamma _{\pm }$ is a circle, generated by a primitive
element $v_{\pm}$ of the lattice
$\Lambda=\{X\in\mathfrak t:\exp(2\pi X)=\mathrm{id}\}$.
This normalization makes a primitive lattice generator have flow period
$2\pi$, as do the angle fields below. It rotates the normal planes of
$\gamma_{\pm}$ with
weight $\pm 1$ because the action is effective. By the van Kampen theorem
applied to invariant tubular neighborhoods of $\gamma _{\pm }$, $\pi
_{1}(N)=\Lambda /(\mathbb{Z}v_{-}+\mathbb{Z}v_{+})$ (compare \cite{Ler}), so
$\{v_{+},v_{-}\}$ is a basis of $\Lambda $. Let $\partial _{\theta _{1}}=\pm
v_{+}$ and $\partial _{\theta _{2}}=\pm v_{-}$, with signs such that $\eta
(\partial _{\theta _{j}})>0$ on $N^{\circ }$. This is possible because $\eta
(\partial _{\theta _{2}})$ is affine in $\psi $, vanishes on $\gamma _{-}$
and not on $\gamma _{+}$, and similarly for $\partial _{\theta _{1}}$.

Put $c=\psi _{+}-\psi _{-}>0$ and $t=(\psi -\psi _{-})/c$, and write $\xi
=a_{1}\partial _{\theta _{1}}+a_{2}\partial _{\theta _{2}}$. As $\eta (\xi
)=1$, the affine functions $\eta (\partial _{\theta _{1}})$, $\eta (\partial
_{\theta _{2}})$ are $(1-t)/a_{1}$ and $t/a_{2}$, with $a_{1},a_{2}>0$. The
$T$-invariant field $\partial _{t}=c\nabla ^{T}\psi /|\nabla ^{T}\psi |^{2}$
lies in $D$ and satisfies $dt(\partial _{t})=1$. All together with the
$T$-action it gives coordinates $(t,\theta _{1},\theta _{2})$ on $N^{\circ
}\cong (0,1)\times T^{2}$ in which $\eta $ is given by \eqref{eq:etaS}. As
in Example~\ref{ex:S3}, the $T$-invariant transverse metric compatible with
$\frac{1}{2}d\eta =\frac{1}{2}dt\wedge \beta $ and orthogonal to $\partial
_{t}$ on the orbits is $h=dt^{2}/(2V)+\frac{V}{2}\beta ^{2}$, with $V=|\nabla
^{T}t|^{2}/2$. Since $\lambda \psi _{\pm }=\kappa e^{\psi _{\pm }}$, one has
$\psi _{+}/\psi _{-}=e^{c}$. Hence $\psi _{-}=c/(e^{c}-1)$ and $\psi
=ct+c/(e^{c}-1)$, and \eqref{eq:3d-identities} becomes
\begin{equation*}
V(t)=\frac{|\nabla ^{T}\psi |^{2}}{2c^{2}}=\frac{\lambda
}{c}\Bigl(t-\frac{e^{ct}-1}{e^{c}-1}\Bigr),
\end{equation*}
which is \eqref{eq:V}. The flow of $\partial _{\theta _{2}}$ rotates the
normal planes of $\gamma _{-}$ with weight $\pm 1$, so $|\partial _{\theta
_{2}}|_{g}$ equals the distance to $\gamma _{-}$ up to second order with
$|\partial _{\theta _{2}}|_{g}^{2}=V/(2a_{2}^{2})+t^{2}/a_{2}^{2}.$ This
gives $V^{\prime }(0)=2a_{2}$, and likewise $V^{\prime }(1)=-2a_{1}$. These
are the boundary values \eqref{eq:Vboundary}. They give $\lambda
=2(a_{1}+a_{2})$ and $\mathcal{H}(c)=a_{2}/(a_{1}+a_{2})$, so that $c=c(a)$
as $c>0$, $a_{1}>a_{2}$. Consequently, the coordinates $(t,\theta
_{1},\theta _{2})$ identify $N^{\circ }$ with the corresponding dense open
subset of $S_{a}$ of Example~\ref{ex:S3}. The complements are unions of two
closed curves, so this identification is an isometry for the restricted
distances and extends to an isometry of the completions $N\rightarrow
\mathbb{S}^{3}$, which is smooth \cite{MyS} and, by continuity, an
isomorphism onto $S_{a}$. Uniqueness of $a$ follows from the invariant $\max
R^{T}/\min R^{T}=e^{|c|}$, as in the proof of
Proposition~\ref{prop:S3-family}(6).

(3) An automorphism $\gamma $ of $(S_{a},\psi _{a})$ preserves $\zeta _{\psi
}$, which is determined by $\psi $ and the Sasakian structure, and commutes
with $T$. If $a_{1}\neq a_{2}$, it preserves $\psi $ and $\nabla ^{T}\psi $,
so it acts on each free orbit $\psi ^{-1}(s)$ by a translation $\sigma
(s)\in T$ and commutes with the flow of $\partial _{t}$. Thus $\sigma $ is
constant and $\gamma \in T$. If $a_{1}=a_{2}$, the automorphisms are the
isometries of the round sphere commuting with the Hopf field; they form
$U(2)$.
\end{proof}

\begin{proof}[Proof of Theorem~\protect\ref{thm:split-class}]
By Theorem~\ref{thm:split-structure}(2), the universal cover is
$\mathbb{C}\times _{\mathrm{H}}N$ with $N$ simply connected, and
$R^{T}=R^{T_{N}}$.

(1) If $R^{T}\equiv 0$, Proposition~\ref{prop:transflat} gives (1), and
$\mathbb{H}_{2}=\mathbb{C}\times _{\mathrm{H}}\mathbb{H}_{1}$ is split.

(2) Otherwise $R^{T}>0$ by Lemma~\ref{lem:nonneg}, and
Proposition~\ref{prop:3d-shrinkers}(2) gives $N\cong S_{a}$, so the universal
cover is $\mathcal{S}_{a}$, with $a$ unique by
Proposition~\ref{prop:S3-family}(6).
Uniqueness of $\Gamma $ up to conjugation is the usual lifting of
isomorphisms to universal covers. By Theorem~\ref{thm:split-structure}(3)
and Proposition~\ref{prop:3d-shrinkers}(3), $\Gamma $ embeds into
$\operatorname{Aut}(S_{a},\psi _{a})$ with image $\Gamma _{N}$ acting freely
on the compact $\mathbb{S}^{3}$, so $\Gamma $ is finite.

(i) If $a_{1}>a_{2}$, then $\Gamma _{N}\subset T^{2}$, and freeness on the
circles $\{w_{2}=0\}$ and $\{w_{1}=0\}$, whose isotropy groups are the two
coordinate circles of $T^{2}$, means that both projections $T^{2}\rightarrow
\mathbb{S}^{1}$ are injective on $\Gamma _{N}$. Hence $\Gamma _{N}$ is cyclic,
generated by $(e^{2\pi i/p},e^{2\pi iq/p})$ with $\gcd (q,p)=1$, and
$\mathbb{S}^{3}/\Gamma _{N}=L(p,q).$ Now with $\chi $ as in
Theorem~\ref{thm:split-structure}(4), the projection $[z,w]\mapsto[w]$
exhibits $M$ as the quotient line bundle defined by the character $\chi$.
The product connection on the trivial bundle upstairs descends to a flat
unitary connection on this quotient line bundle, and $M$ is split
if and
only if $\chi $ is trivial.

(ii) If $a_{1}=a_{2}$, $\mathcal{S}_{a}$ is the cylinder lift with $d=1$ by
Proposition~\ref{prop:S3-family}(2).

This proves (2).

(3) The transverse structure is the product of the flat plane and the
surface of curvature $K=\kappa e^{\psi _{N}}>0$, which gives the assertions
on the curvature operator, on $\mathrm{Ric}^{T}$ and, by
Proposition~\ref{prop:3d-shrinkers}(2), on $R^{T}=2K$. In the model
\eqref{eq:H-product}, $\nabla ^{T}f$ is the sum of $\lambda $ times the
horizontal lift of the
position vector $z$ and $\nabla ^{T}\psi _{N}$, so the critical set of $f$
is $\{0\}\times \operatorname{Crit}(\psi _{N})$. It is $\Gamma $-invariant,
and each
of its components is preserved by $\Gamma $, which preserves $f$.
Proposition~\ref{prop:3d-shrinkers}(2) and Proposition~\ref{prop:S3-family}(5)
identify these components, their levels and the lengths of their images
in $M$. The Hessian of $f$ along $\{0\}\times \gamma _{+}$ is $\lambda $ on
the plane and $\lambda (1-\psi _{+})<0$ on $D_{N}$.

(4) By \eqref{eq:3d-identities}, $R^{T}=2\kappa e^{\psi _{N}}$ is constant
if and only if $\psi _{N}$ is constant, and then $\kappa =\lambda /e$.
Moreover $\min f=\psi _{-}$, which is smaller than one unless $\kappa
=\lambda /e$. So (a), (b) and (c) are equivalent, and (3) gives the
equivalence with (d). At a critical point, \eqref{eq:pre-normal} gives
$R^{T}=2\lambda f\in \{2\lambda \psi _{-},2\lambda \psi _{+}\}$, which is
$2\lambda $ only if $\psi _{\pm }=1 $. So (e) is equivalent to (a). If
$a_{1}=a_{2}$, $M$ is a finite quotient of a regular lift, so (f) holds. If
$a_{1}\neq a_{2}$, the Reeb foliation of the simply connected manifold
$\mathcal{S}_{a}$ has nonclosed leaves or leaves with nontrivial holonomy
(Proposition~\ref{prop:S3-family}(2),(3)), so it is not globally
developable, and neither is that of $M$.
Moreover (f) fails. This is clear if $a_{1}/a_{2}\notin \mathbb{Q}$, since
then $M$ is irregular. If $a_{1}/a_{2}=m_{1}/m_{2}\in \mathbb{Q}$ with
$(m_{1},m_{2})\neq (1,1)$, then $\Gamma $ commutes with the circle action, so
the leaf space $\mathbb{C}\times \mathbb{CP}_{[m_{1},m_{2}]}^{1}$ of
$\mathcal{S}_{a}$ is a finite orbifold cover of the leaf space of $M$, as in
the proof of Theorem~\ref{thm:globdev}. An orbifold cover of a developable
orbifold is developable, whereas $\mathbb{CP}_{[m_{1},m_{2}]}^{1}$ is a bad
orbifold. Hence the leaf space of $M$ is not developable. If $a_{1}=a_{2} $,
$N=S_{a}$ is round with $K\equiv \lambda $ and $\psi _{a}\equiv 1$, so (g)
holds. Conversely (g) gives $R^{T}=R^{T_{N}}\equiv 2\lambda $, which is (b).
The last assertion follows from Proposition~\ref{prop:S3-family}(2),(3) and
Theorem~\ref{thm:pi1}(2).

(5) By Lemma~\ref{lem:split-affine}(1), $\mathrm{Ric}^{T}(V,V)=0$ on the
universal cover of a globally split soliton, at every point.
\end{proof}

\begin{lemma}
\label{lem:transverse-CY}Let $(M^{5},\xi ,\eta ,\Phi ,g)$ be a connected
complete noncompact Sasakian manifold such that the closure $H$ of the Reeb
flow in $\operatorname{Isom}(M,g)$ is a compact torus, and assume that
\begin{equation*}
\mathrm{Ric}^{T}\geq 0.
\end{equation*}
For an $H$-orbit $P$, put $\rho _{P}=d_{g}(P,\cdot )$ and
$V_{P}(r)=\operatorname{Vol}_{g}\{\rho _{P}<r\}$. This is the
five-dimensional Riemannian volume of a tube in $M$, not a volume on a
transverse quotient. The exponent $4$ below belongs to the comparison
estimate and does not change the measure being used. Then, for every
$H$-orbit $P$:

\begin{enumerate}
\item $r\mapsto r^{-4}V_{P}(r)$ is nonincreasing on $(0,\infty )$;

\item the volume grows at least linearly:
\begin{equation}
V_{P}(2r+3)\geq \tfrac{r}{8}\,V_{P}(1)\qquad (r\geq 1).  \label{eq:CY}
\end{equation}
In particular, $\operatorname{Vol}_{g}(M)=\infty $.
\end{enumerate}
\end{lemma}

\begin{proof}
Every element of $H$ preserves $\eta$. If $X,Y$ are the fundamental vector
fields of elements of its abelian Lie algebra, then
\begin{equation*}
X(\eta(Y))=(\mathcal{L}_X\eta)(Y)+\eta([X,Y])=0, \qquad d\eta(X,Y)=0.
\end{equation*}
Since $\xi$ is tangent to every orbit, the projection of $TP$ to $D$ is an
isotropic subspace of the four-dimensional symplectic space $(D,d\eta|_D)$
and has dimension $\dim P-1$. Consequently $1\leq d:=\dim P\leq3$.

Fix an orbit $P$. It is a compact embedded submanifold. Its distance
function $u=\rho_P$ is proper and $H$-invariant, hence basic. Since $\xi$ is
tangent to $P$, the normal spaces of $P$ lie in $D$; in particular a
minimizing geodesic leaving $P$ is initially perpendicular to $\xi$. Since
$\xi$ is Killing, $g(\xi,\dot\gamma)$ is constant along the geodesic, so it
remains horizontal. Write $\nu^1P$ for the unit normal bundle, $d\sigma$ for
its natural measure, and $c(v)$ for the normal cut distance. Outside $P$ and
its cut locus, the normal exponential map gives coordinates $E(r,v)=\exp(rv)$,
$0<r<c(v)$, with volume density $\mathcal{J}(r,v)>0$ satisfying
\begin{equation*}
E^*dV_g=\mathcal{J}(r,v)\,dr\,d\sigma(v),\qquad \partial_r\log
\mathcal{J}=(\Delta_g u)\circ E.
\end{equation*}
The excluded sets have Riemannian measure zero.

We prove $\Delta_g u\leq3/r$ along each such normal geodesic. Put $W=\nabla_g
u=\nabla^T u$ and $S=\operatorname{Hess}^T u|_{D\cap W^\perp}$. At a
regular point, $|W|=1$, $\operatorname{Hess}^T u(W,\cdot)=0$, and
\begin{equation*}
\operatorname{Hess}_g u(\xi,\xi)=0,\qquad \operatorname{Hess}_g u(X,\xi)=g(\Phi
W,X)\quad(X\in D).
\end{equation*}
On $D\times D$ the ambient and transverse Hessians agree. Thus, with
$m=\Delta_g u=\operatorname{tr}S$,
\begin{equation*}
|\operatorname{Hess}_g u|^2=|S|^2+2,\qquad
\mathrm{Ric}_g(W,W)=\mathrm{Ric}^T(W,W)-2.
\end{equation*}
The two constant terms cancel in the ordinary Bochner formula for the unit
gradient $W$. Along $\gamma_v(r)=E(r,v)$ it follows that
\begin{equation*}
m^{\prime }=-|S|^{2}-\mathrm{Ric}^T(W,W) \leq-|S|^2\leq-\frac{m^2}{3}.
\end{equation*}
Here the first inequality uses $\mathrm{Ric}^{T}\geq 0$, and the last one uses
the dimension three of $D\cap W^\perp$. This
calculation does not impose a minimality or mean-curvature condition on $P$.

The normal-tube expansion gives
\begin{equation*}
\mathcal{J}(r,v)=r^{4-d}(1+O(r)),\qquad
m(r)=\frac{4-d}{r}+O(1)\longrightarrow+\infty \quad(r\downarrow0).
\end{equation*}
While $m>0$, the Riccati inequality implies
\begin{equation*}
\left(\frac1m\right)^{\prime }=-\frac{m^{\prime }}{m^2}\geq\frac13.
\end{equation*}
Integrating from $\varepsilon$ to $r$ and letting $\varepsilon\downarrow0$
yields $m(r)\leq3/r$. Moreover $m^{\prime }\leq0$, so if $m$ reaches zero it
stays nonpositive thereafter. The same bound therefore holds throughout
$0<r<c(v)$.

It follows that $\mathcal{J}(r,v)/r^3$ is nonincreasing. Extend
$\mathcal{J}(r,v)$ by zero for $r\geq c(v)$ and denote the resulting density
by $\bar{\mathcal{J}}(r,v)$. Then $\bar{\mathcal{J}}(r,v)/r^3$ is still
nonincreasing, and the normal-coordinate volume formula gives
\begin{equation*}
\frac{V_P(r)}{r^4} =\int_{\nu^1P}\int_0^1
\frac{\bar{\mathcal{J}}(r\theta,v)}{(r\theta)^3}\,\theta^3\,d\theta\,d\sigma(v).
\end{equation*}
This proves (1).

For (2), fix $r\geq1$. Since $M$ is connected and noncompact and $\rho_P$ is
proper, $\rho_P$ takes every value in $[0,\infty)$; choose $x$ with
$\rho_P(x)=r+1$. Set $P^{\prime }=Hx$. Invariance under $H$ gives $\rho_P=r+1$
on $P^{\prime }$ and $\rho_{P^{\prime }}=r+1$ on $P$. The triangle
inequality therefore gives
\begin{equation*}
\{\rho_P<1\}\subset\{r<\rho_{P^{\prime }}<r+2\},\qquad \{\rho_{P^{\prime
}}<r+2\}\subset\{\rho_P<2r+3\}.
\end{equation*}
By (1) and $(1+2/r)^{-4}\geq1-8/r$,
\begin{equation*}
\begin{split}
V_P(1)&\leq V_{P^{\prime }}(r+2)-V_{P^{\prime }}(r) \\
&\leq\left[1-\left(\frac{r}{r+2}\right)^4\right]V_{P^{\prime }}(r+2)
\leq\frac8r V_P(2r+3).
\end{split}
\end{equation*}
This is \eqref{eq:CY}. Since $V_P(1)>0$, letting $r\to\infty$ proves
$\operatorname{Vol}_g(M)=\infty$.
\end{proof}

\vskip12pt

\noindent{\bf Acknowledgements} Shu-Cheng Chang is partially supported by the Startup Foundation
for Advanced Talents of the Shanghai Institute for Mathematics and
Interdisciplinary Sciences (No.~2302-SRFP-2024-0049). Chien
Lin is partially supported by the NSTC (No.~114-2115-M-003-002-MY2).
Hongbing Qiu is partially supported by the NSFC
(No.~12471050).

\vskip12pt

\noindent{\bf Declaration of AI use.} The authors used ChatGPT (OpenAI) 
during the
preparation of the manuscript for language editing and exploratory calculations.
All mathematical arguments and computations were independently checked by the
authors.


\vskip24pt

\begin{thebibliography}{CHLW3}
\bibitem[BCCD]{BCCD} R. H. Bamler, C. Cifarelli, R. J. Conlon and A.
Deruelle, \textit{A new complete two-dimensional shrinking gradient
K\"{a}hler--Ricci soliton}, Geom. Funct. Anal. 34 (2024), no. 2, 377--392;
arXiv:2206.10785.

\bibitem[BG]{BG} C. P. Boyer and K. Galicki, \textit{Sasakian Geometry},
Oxford Mathematical Monographs, Oxford University Press, Oxford, 2008.

\bibitem[Bir]{Bir} C. Birkar, \textit{Singularities of linear systems and
boundedness of Fano varieties}, Ann. of Math. (2) 193 (2021), 347--405.

\bibitem[BM]{BM} J. Bernstein and T. Mettler, \textit{Two-dimensional
gradient Ricci solitons revisited}, Int. Math. Res. Not. IMRN 2015, no. 1,
78--98.

\bibitem[Cal]{Cal} E. Calabi, \textit{On manifolds with non-negative Ricci
curvature II}, Notices Amer. Math. Soc. 22 (1975), A205.

\bibitem[Cao1]{Cao96} H.-D. Cao, \textit{Existence of gradient
K\"{a}hler--Ricci solitons}, in: Elliptic and Parabolic Methods in Geometry
(Minneapolis, MN, 1994), A K Peters, Wellesley, MA, 1996, 1--16.

\bibitem[Cao2]{Cao} H.-D. Cao, \textit{Recent progress on Ricci solitons},
in: Recent Advances in Geometric Analysis, Adv. Lect. Math. (ALM) 11, Int.
Press, Somerville, MA, 2010, 1--38.

\bibitem[CaZ]{CaZ} H.-D. Cao and D. Zhou, \textit{On complete gradient
shrinking Ricci solitons}, J. Differential Geom. 85 (2010), no. 2, 175--185;
arXiv:0903.3932.

\bibitem[CCD]{CCD} C. Cifarelli, R. J. Conlon and A. Deruelle, \textit{On
finite time Type I singularities of the K\"{a}hler--Ricci flow on compact
K\"{a}hler surfaces}, J. Eur. Math. Soc. 28 (2026), no. 2, 463--504;
arXiv:2203.04380.

\bibitem[CDS]{CDS} R. J. Conlon, A. Deruelle and S. Sun,
\textit{Classification results for expanding and shrinking gradient
K\"{a}hler--Ricci solitons}, Geom. Topol. 28 (2024), no. 1, 267--351;
arXiv:1904.00147.

\bibitem[CGT]{CGT} J. Cheeger, M. Gromov and M. Taylor, \textit{Finite
propagation speed, kernel estimates for functions of the Laplace operator,
and the geometry of complete Riemannian manifolds}, J. Differential Geom. 17
(1982), 15--53.

\bibitem[CHLW]{CHLWHT} S.-C. Chang, Y.-B. Han, C. Lin and C.-T. Wu,
\textit{On the Hamilton--Tian conjecture in a compact transverse Fano Sasakian
5-manifold}, arXiv:2605.20852.

\bibitem[CHW]{CHW} S.-C. Chang, Y.-B. Han and C.-T. Wu, \textit{Geometry and
topology of gradient shrinking Sasaki--Ricci solitons}, arXiv:2508.13495.

\bibitem[CLL]{CLL} S.-C. Chang, F. Li and C. Lin, \textit{Geometry of
shrinking Sasaki--Ricci solitons I: fundamental equations and
characterization of rigidity}, arXiv:2502.16148.

\bibitem[CLQ]{PartI} S.-C. Chang, C. Lin and H. Qiu, \textit{On the
structure of gradient shrinking Sasaki--Ricci solitons of dimension five I:
curvature estimates and local transverse splitting}, preprint, 2026.

\bibitem[Co]{Co1} T. C. Collins, \textit{Stability and convergence of the
Sasaki--Ricci flow}, J. Reine Angew. Math. 716 (2016), 1--27.

\bibitem[CTZ]{CTZ} S.-C. Chang, J. Tie and J. Zhou, \textit{Complete gradient
steady Sasaki--Ricci solitons}, preprint, 2026.

\bibitem[FIK]{FIK} M. Feldman, T. Ilmanen and D. Knopf, \textit{Rotationally
symmetric shrinking and expanding gradient K\"{a}hler--Ricci solitons}, J.
Differential Geom. 65 (2003), no. 2, 169--209.

\bibitem[FOW]{FOW} A. Futaki, H. Ono and G. Wang, \textit{Transverse
K\"{a}hler geometry of Sasaki manifolds and toric Sasaki--Einstein manifolds},
J.
Differential Geom. 83 (2009), no. 3, 585--635.

\bibitem[GPSS]{GPSS} B. Guo, D. H. Phong, J. Song and J. Sturm,
\textit{Compactness of K\"{a}hler--Ricci solitons on Fano manifolds}, Pure
Appl.
Math. Q. 18 (2022), no. 1, 305--316.

\bibitem[Ham]{Ham} R. S. Hamilton, \textit{The formation of singularities in
the Ricci flow}, in: Surveys in Differential Geometry, Vol. II, Int. Press,
Cambridge, MA, 1995, 7--136.

\bibitem[Hat]{Hat} A. Hatcher, \textit{Algebraic Topology}, Cambridge
University Press, Cambridge, 2002.

\bibitem[He]{He} W. He, \textit{The Sasaki--Ricci flow and compact Sasaki
manifolds of positive transverse holomorphic bisectional curvature}, J.
Geom. Anal. 23 (2013), 1876--1931.

\bibitem[Her]{Her} R. Hermann, \textit{A sufficient condition that a mapping
of Riemannian manifolds be a fibre bundle}, Proc. Amer. Math. Soc. 11
(1960), 236--242.

\bibitem[HS]{HS} A. D. Hwang and M. A. Singer, \textit{A momentum
construction for circle-invariant K\"{a}hler metrics}, Trans. Amer. Math.
Soc. 354 (2002), no. 6, 2285--2325.

\bibitem[Koi]{Koi} N. Koiso, \textit{On rotationally symmetric Hamilton's
equation for K\"{a}hler--Einstein metrics}, in: Recent Topics in Differential
and Analytic Geometry, Adv. Stud. Pure Math. 18-I, Academic Press, Boston, MA,
1990, 327--337.

\bibitem[Kol1]{Kol1} J. Koll\'{a}r, \textit{Seifert $G_{m}$-bundles},
preprint, 2004; arXiv:math/0404386.

\bibitem[Kol2]{Kol2} J. Koll\'{a}r, \textit{Einstein metrics on
five-dimensional Seifert bundles}, J. Geom. Anal. 15 (2005), no. 3, 445--476.

\bibitem[Ler]{Ler} E. Lerman, \textit{Contact toric manifolds}, J.
Symplectic Geom. 1 (2003), no. 4, 785--828.

\bibitem[LW]{LW} Y. Li and B. Wang, \textit{On K\"{a}hler Ricci shrinker
surfaces}, Acta Math. 236 (2026), no. 1, 1--50; arXiv:2301.09784.

\bibitem[LZ]{LZ} Y. Li and J. Zhang, \textit{Scalar curvature bounds and
tangent cones for K\"{a}hler--Ricci shrinkers}, arXiv:2610.05643.

\bibitem[Mol]{Mol} P. Molino, \textit{Riemannian Foliations}, Progress in
Mathematics 73, Birkh\"{a}user, Boston, 1988.

\bibitem[MW1]{MW1} O. Munteanu and J. Wang, \textit{Geometry of shrinking
Ricci solitons}, Compos. Math. 151 (2015), no. 12, 2273--2300.

\bibitem[MW2]{MW2} O. Munteanu and J. Wang, \textit{Structure at infinity
for shrinking Ricci solitons}, Ann. Sci. \'{E}c. Norm. Sup\'{e}r. (4) 52
(2019), no. 4, 891--925.

\bibitem[MW3]{MWpos} O. Munteanu and J. Wang, \textit{Positively curved
shrinking Ricci solitons are compact}, J. Differential Geom. 106 (2017),
499--505.

\bibitem[MyS]{MyS} S. B. Myers and N. E. Steenrod, \textit{The group of
isometries of a Riemannian manifold}, Ann. of Math. (2) 40 (1939), 400--416.

\bibitem[Nab]{Nab} A. Naber, \textit{Noncompact shrinking four solitons with
nonnegative curvature}, J. Reine Angew. Math. 645 (2010), 125--153;
arXiv:0710.5579.

\bibitem[Ni]{Ni} L. Ni, \textit{Ancient solutions to K\"{a}hler--Ricci flow},
Math. Res. Lett. 12 (2005), no. 5--6, 633--653.

\bibitem[Ru]{Ru} P. Rukimbira, \textit{Chern--Hamilton's conjecture and
K-contactness}, Houston J. Math. 21 (1995), no. 4, 709--718.

\bibitem[Shi]{Shi} W.-X. Shi, \textit{Deforming the metric on complete
Riemannian manifolds}, J. Differential Geom. 30 (1989), no. 1, 223--301.

\bibitem[Ste]{Ste} N. Steenrod, \textit{The Topology of Fibre Bundles},
Princeton Mathematical Series 14, Princeton University Press, Princeton,
1951.

\bibitem[SWZ]{SWZ} K. Smoczyk, G. Wang and Y. Zhang, \textit{The
Sasaki--Ricci flow}, Internat. J. Math. 21 (2010), no. 7, 951--969.

\bibitem[Tan]{Tan} S. Tanno, \textit{Sasakian manifolds with constant $\phi
$-holomorphic sectional curvature}, T\^{o}hoku Math. J. (2) 21 (1969), no. 3,
501--507.

\bibitem[Wad]{Wad} A. W. Wadsley, \textit{Geodesic foliations by circles},
J. Differential Geom. 10 (1975), no. 4, 541--549.

\bibitem[WZh]{WZh} G. Wu and S. Zhang, \textit{Remarks on shrinking gradient
K\"{a}hler--Ricci solitons with positive bisectional curvature}, C. R. Math.
Acad. Sci. Paris 354 (2016), no. 7, 713--716; arXiv:1906.00412.

\bibitem[Wu]{Wu} L.-F. Wu, \textit{The Ricci flow on $2$-orbifolds with
positive curvature}, J. Differential Geom. 33 (1991), 575--596.

\bibitem[WZ]{WZ} X.-J. Wang and X. Zhu, \textit{K\"{a}hler--Ricci solitons
on toric manifolds with positive first Chern class}, Adv. Math. 188 (2004),
no. 1, 87--103.

\bibitem[WZ2]{WZ2} G. Wang and Y. Zhang, \textit{The Sasaki--Ricci flow on
Sasakian 3-spheres}, Commun. Math. Stat. 1 (2013), no. 1, 43--71.

\bibitem[Yau]{Yau} S.-T. Yau, \textit{Some function-theoretic properties of
complete Riemannian manifolds and their applications to geometry}, Indiana
Univ. Math. J. 25 (1976), 659--670.
\end{thebibliography}
\end{document}